\documentclass[a4paper,twoside,reqno,tbtags]{article}

\usepackage{amssymb,amsmath,amsthm,mathrsfs}
\usepackage{graphicx}
\usepackage{times}
\usepackage{natbib}

\leftmargini=\baselineskip

\newtheoremstyle{localthm}
	{5pt} 
	{5pt} 
	{\sl} 
	{} 
	{\bf} 
	{{\rm.}} 
	{.7em} 
	{} 

\theoremstyle{localthm}
\newtheorem{theorem}{Theorem}[section]
\newtheorem{corollary}[theorem]{Corollary}
\newtheorem{fact}[theorem]{Fact}
\newtheorem{proposition}[theorem]{Proposition}
\newtheorem{lemma}[theorem]{Lemma}

\newtheoremstyle{localrem}
	{5pt} 
	{5pt} 
	{\rm} 
	{} 
	{\bf} 
	{{\rm.}} 
	{.7em} 
	{} 

\theoremstyle{localrem}
\newtheorem{condition}[theorem]{Condition}
\newtheorem{remark}[theorem]{Remark}
\newtheorem{example}[theorem]{Example}

\newcommand{\FF}{{\mathbb F}}
\newcommand{\GG}{{\mathbb G}}

\newcommand{\RR}{{\mathbb R}}
\newcommand{\UU}{{\mathbb U}}

\newcommand{\WW}{{\mathbb W}}

\newcommand{\ZZ}{{\mathbb Z}}

\newcommand{\logit}{\mathop{\mathrm{logit}}\nolimits}
\newcommand{\TT}{\mathcal{T}}

\DeclareMathOperator{\argmax}{argmax}

\newcommand{\eps}{\varepsilon}

\begin{document}

\title{A New Approach to Tests and Confidence Bands for Distribution Functions}

\author{Lutz D\"umbgen and Jon A.\ Wellner\\
	University of Bern and University of Washington, Seattle}

\date{October 2022}

\maketitle

\begin{abstract}
We introduce new goodness-of-fit tests and corresponding confidence bands for distribution functions. They are inspired by multi-scale methods of testing and based on refined laws of the iterated logarithm for the normalized uniform empirical process $\UU_n (t)/\sqrt{t(1-t)}$ and its natural limiting process, the normalized Brownian bridge process $\UU (t)/\sqrt{t(1-t)}$.  The new tests and confidence bands refine the procedures of Berk and Jones (1979) and Owen (1995). Roughly speaking, the high power and accuracy of the latter methods in the tail regions of distributions are essentially preserved while gaining considerably in the central region. 
The goodness-of-fit tests perform well in signal detection problems involving sparsity, as in Ingster (1997), Donoho and Jin (2004) and Jager and Wellner (2007), but also under contiguous alternatives. Our analysis of the confidence bands sheds new light on the influence of the underlying $\phi$-divergences.
\end{abstract}

\paragraph{AMS subject classifications.}
60E10, 60F10 (primary); 62D99 (secondary).

\paragraph{Key words.}
Confidence band, goodness-of-fit, law of the iterated logarithm, limit distribution, multi-scale test statistics

\tableofcontents

\section{Introduction and motivations}
\label{sec:intro}

\subsection{Some well-known facts}
\label{subsec:EDFknownFacts} 
Let $\FF_n$ be the empirical distribution function of independent random variables $X_1, X_2, \ldots, X_n$ with unknown distribution function $F$ on the real line. The main topic of the present paper is to construct a confidence band $(A_{n,\alpha},B_{n,\alpha})$ for $F$ with given confidence level $1 - \alpha \in (0,1)$. That is, $A_{n,\alpha} = A_{n,\alpha}(\cdot, (X_i)_{i=1}^n)$ and $B_{n,\alpha} = B_{n,\alpha}(\cdot, (X_i)_{i=1}^n)$ are data-driven functions on the real line such that for any true distribution function $F$,
\begin{equation}
\label{ineq:conf.level}
	P_F(A_{n,\alpha} \le F \le B_{n,\alpha} \ \text{on} \ \RR) \ge 1 - \alpha .
\end{equation}
Let us recall some well-known facts about $\FF_n$ (cf.\ \cite{Shorack_Wellner_1986, Shorack_Wellner_2009}). The stochastic process $\bigl( \FF_n(x) \bigr)_{x \in \RR}$ has the same distribution as $\bigl( \GG_n(F(x)) \bigr)_{x \in \RR}$, where $\GG_n$ is the empirical distribution of independent random variables $\xi_1, \xi_2, \ldots, \xi_n$ with uniform distribution on $[0,1]$. This enables the well-known Kolmogorov--Smirnov confidence bands: let
\[
	\UU_n(t) := \sqrt{n} (\GG_n(t) - t) ,
\]
and let $\kappa_{n,\alpha}^{\rm KS}$ be the $(1 - \alpha)$-quantile of the supremum norm $\|\UU_n\|_\infty := \sup_{t \in [0,1]} |\UU_n(t)|$. Then the confidence band $(A_{n,\alpha}^{\rm KS}, B_{n,\alpha}^{\rm KS})$ with $A_{n,\alpha}^{\rm KS} := \max(\FF_n - n^{-1/2} \kappa_{n,\alpha}^{\rm KS}, 0)$ and $B_{n,\alpha}^{\rm KS} := \min(\FF_n + n^{-1/2}\kappa_{n,\alpha}^{\rm KS}, 1)$ satisfies \eqref{ineq:conf.level} with equality if $F$ is continuous. 
Since $\UU_n$ converges in distribution in $\ell^\infty([0,1])$ to standard Brownian bridge $\UU$, $\kappa_{n,\alpha}^{\rm KS}$ converges to the $(1 - \alpha)$-quantile $\kappa_\alpha^{\rm KS}$ of $\|\UU\|_\infty$. In particular, the width $B_{n,\alpha}^{\rm KS} - A_{n,\alpha}^{\rm KS}$ of the Kolmogorov--Smirnov band is bounded uniformly by $2 n^{-1/2} \kappa_{n,\alpha}^{\rm KS} = O(n^{-1/2})$. (Throughout this paper, asymptotic statements refer to $n \to \infty$, unless stated otherwise.) On the other hand, it is well-known that Kolmogorov-Smirnov confidence bands give little or no information in the tails of the distribution $F$; see e.g.\ \cite{Milbrodt_Strasser_1990}, \cite{Janssen_1995}, and \cite{Lehmann_Romano_2005}, chapter 14, for a useful summary.

\subsection{Confidence bands by inversion of tests}
\label{subsec:CBinGeneralAndKS}  
In general, confidence bands can be obtained by inverting goodness-of-fit tests. For a given continuous distribution function $F_0$, let $T_n(F_0) = T_n(F_0, (X_i)_{i=1}^n)$ be some test statistic for the null hypothesis that $F \equiv F_0$. Suppose that for any test level $\alpha \in (0,1)$, the $(1 - \alpha)$-quantile $\kappa_{n,\alpha}$ of $T_n(F_0)$ under the null hypothesis does not depend on $F_0$. Then a $(1 - \alpha)$-confidence band $(A_{n,\alpha},B_{n,\alpha})$ for a continuous distribution function $F$ is given by
\[
	A_{n,\alpha}(x) := \inf \bigl\{ F(x) \colon T_n(F) \le \kappa_{n,\alpha} \bigr\} , \quad
	B_{n,\alpha}(x) := \sup \bigl\{ F(x) \colon T_n(F) \le \kappa_{n,\alpha} \bigr\} .
\]
Depending on the specific choice of $T_n$, these functions $A_{n,\alpha}$ and $B_{n,\alpha}$ can be computed explicitly, and the constraint \eqref{ineq:conf.level} is even satisfied for arbitrary, possibly noncontinuous distribution functions $F$; see Section~S.6 for further details.

Since $(A_{n,\alpha}^{\rm KS}, B_{n,\alpha}^{\rm KS})$ corresponds to $T_n^{\rm KS}(F_0) := \sqrt{n} \|\FF_n - F_0\|_\infty$, one possibility to enhance precision in the tails is to consider weighted supremum norms such as
\begin{equation}
\label{eq:weighted.sup-norm.1}
	T_n(F_0) := \ \sup_{x\colon 0 < F_0(x) < 1}
		\frac{\sqrt{n} |\FF_n - F_0|}{w(F_0)}(x)
\end{equation}
or
\begin{equation}
\label{eq:weighted.sup-norm.2}
	T_n(F_0) := \ \sup_{x \in [X_{n:1}, X_{n:n})}
		\frac{\sqrt{n} |\FF_n - F_0|}{w(\FF_n)}(x) ,
\end{equation}
where $X_{n:1} \le X_{n:2} \le \cdots \le X_{n:n}$ are the order statistics of $X_1, X_2, \ldots, X_n$. Here, $w : (0,1) \to (0,\infty)$ is some continuous weight function such that $w(1-t) = w(t)$ for $0 < t < 1$ and $w(t) \to 0$ as $t \to 0$. Specific proposals include
\[
	w(t) := \sqrt{t(1 - t) h(t)} ,
\]
where $h \equiv 1$, see \cite{Jaeschke_1979} and \cite{Eicker_1979}, or 
$h(t) \to \infty$ sufficiently fast as $t \to 0$, see \cite{OReilly_1974} 
or \cite{Csorgo_etal_1986}. Specifically, Stepanova and Pavlenko \cite{Stepanova_Pavlenko_2018} propose to construct confidence bands with the test statistic 
\eqref{eq:weighted.sup-norm.2} and $h(t) := \log \log(1/[t(1-t)])$. The latter choice is motivated by the law of the iterated logarithm (LIL) for the Brownian bridge process $\mathbb{U}$, stating that
\begin{equation}
\label{eq:LIL.BB}
	\limsup_{t \searrow 0} \frac{\mathbb{U}(t)}{\sqrt{2 t \log\log(1/t)}}
	\ = \ \limsup_{t \nearrow 1} \frac{\mathbb{U}(t)}{\sqrt{2 (1-t) \log\log(1/(1-t))}}
	\ = \ 1
\end{equation}
almost surely. 

\subsection{The tests of Berk and Jones and Owen's bands}
\label{BerkJonesTestandOwenBands} 

Another goodness-of-fit test, proposed by Berk and Jones \cite{Berk_Jones_1979}, uses the test statistic
\begin{equation}
\label{eq:SupFirst}
	T_n^{\rm BJ}(F_0) := n \sup_{x \colon 0 < F_0(x) < 1} K(\FF_n(x),F_0(x)) ,
\end{equation}
where
\[
	K(u,t) := u \log \Bigl( \frac{u}{t} \Bigr)
		+ (1 - u) \log \Bigl( \frac{1-u}{1-t} \Bigr) 
\]
for $u \in [0,1]$ and $t \in (0,1)$. Note that $K(u,t)$ is the Kullback-Leibler divergence between the $\mathrm{Bernoulli}(u)$ and $\mathrm{Bernoulli}(t)$ distributions. Owen \cite{Owen_1995} proposed and analyzed confidence bands for $F$ based on this test statistic. As noted by \cite{Jager_Wellner_2007}, the test statistic $T_n^{\rm BJ}(F_0)$ can be embedded into a general family of test statistics $T_{n,s}^{\rm BJ}(F_0)$, $s \in \RR$. Let
\begin{equation}
\label{eq:SupFirst.S}
	T_{n,s}^{\rm BJ}(F_0) := \begin{cases}
		\displaystyle
		\sup_{x \colon 0 < F_0(x) < 1} n K_s(\FF_n(x),F_0(x))
		& \text{if} \ s > 0 , \\
		\displaystyle
		\sup_{x \in [X_{n:1},X_{n:n})} n K_s(\FF_n(x),F_0(x))
		& \text{if} \ s \le 0 ,
	\end{cases}
\end{equation}
with the following divergence function $K_s$: for $t, u \in (0,1)$,
\begin{equation}
\label{TheKFunctions}
	K_s(u,t) = \begin{cases}
		\bigl( t (u/t)^s + (1 - t)[(1 - u)/(1 - t)]^s - 1 \bigr)/[s(s-1)] ,
			& s \ne 0,1, \\
		u \log(u/t) + (1 - u) \log[(1-u)/(1-t)] ,
			& s = 1, \\
		t \log(t/u) + (1 - t) \log[(1 - t)/(1 - u)] ,
			& s = 0.
	\end{cases}		
\end{equation}
(An alternative representation of $K_s$ is given in \eqref{eq:Ks.via.phis}.) Moreover, for fixed $t \in (0,1)$ and $u \in \{0,1\}$, the limit $K(u,t) := \lim_{u' \to u} K_s(u',t)$ equals $\infty$ if $s \le 0$ and exists in $(0,\infty)$ otherwise. A detailed discussion of these divergences is given in Section~S.3 of the online supplement. At present it suffices to note that for any fixed $t \in (0,1)$, $K_s(u,t)$ is strictly convex in $u$ with unique minimum $0$ at $u = t$ and second derivative $[t(1 - t)]^{-1}$ there. Interesting special cases are $K = K_1$, $K_{1/2}(u,t) = 4 \bigl( 1 - \sqrt{ut} - \sqrt{(1 - u)(1 - t)} \bigr)$ and
\[
	K_2(u,t) = \frac{(u - t)^2}{2t(1-t)} , \quad
	K_{-1}(u,t) = \frac{(u - t)^2}{2u(1-u)} .
\]
Consequently, if $w(t) := \sqrt{t(1 - t)}$, then the test statistic $T_{n,2}^{\rm BJ}(F_0)$ coincides with $0.5$ times the square of $T_n(F_0)$ in \eqref{eq:weighted.sup-norm.1}, and $T_{n,-1}^{\rm BJ}(F_0)$ equals $0.5$ times the square of \eqref{eq:weighted.sup-norm.2}. As shown by \cite{Jager_Wellner_2007}, for any $s \in [-1,2]$, the null distribution of $T_{n,s}^{\rm BJ}(F_0)$ has the same asymptotic behavior, and the corresponding $(1 - \alpha)$-quantiles $\kappa_{n,s,\alpha}^{\rm BJ}$ satisfy
\begin{equation}
\label{eq:CenterFunctionOfnForBJ}
	\kappa_{n,s,\alpha}^{\rm BJ} \
	= \ \log\log n + 2^{-1} \log \log \log n + O(1) .
\end{equation}
From this one can deduce that the resulting confidence band $(A_{n,s,\alpha}^{\rm BJO}, B_{n,s,\alpha}^{\rm BJO})$ for $F$ satisfies
\[
	B_{n,s,\alpha}^{\rm BJO}(x) - A_{n,s,\alpha}^{\rm BJO}(x)
	\le 2 \sqrt{ 2 \gamma_n \, \FF_n(1 - \FF_n)(x) } + 4 \gamma_n
\]
where $\gamma_n := n^{-1} \kappa_{n,s,\alpha}^{\rm BJ} = (1 + o(1)) n^{-1} \log \log n$; see Lemma~S.12 in Section~S.3. Hence the band $(A_{n,s,\alpha}^{\rm BJO}, B_{n,s,\alpha}^{\rm BJO})$ is substantially more accurate than $(A_{n\alpha}^{\rm KS}, B_{n,\alpha}^{\rm KS})$ in the tail regions. But in the central region, i.e.\ when $\FF_n(x)$ is bounded away from $0$ and $1$, they are of width $O(n^{-1/2} (\log \log n)^{1/2})$ rather than $O(n^{-1/2})$.

\subsection{Goals revisited}
\label{GoalsAgain}
The goal of Berk and Jones \cite{Berk_Jones_1979} was to find goodness-of-fit tests with optimal Bahadur efficiencies. They interpret their test statistic $T_n^{\rm BJ}(F_0)$ also as a union-intersection test statistic, where $n K(\FF_n(x), F_0(x))$ is the negative likelihood ratio statistic for the null hypothesis that $F(x) = F_0(x)$, based on the binomial distribution of $n \FF_n(x)$. The union-intersection and related paradigms for the present goodness-of-fit testing problem have been treated in more generality by \cite{Gontscharuk_etal_2016}.

In view of the previous considerations, the confidence band $(A_{n,\alpha}^{\rm SP}, B_{n,\alpha}^{\rm SP})$ of \cite{Stepanova_Pavlenko_2018}, based on the test statistic
\begin{equation}
\label{eq:TnSP}
	T_n^{\rm SP}(F_0) := \sup_{x \in [X_{n:1},X_{n:n})}
		\frac{\sqrt{n} |\FF_n - F_0|}{\sqrt{\FF_n(1 - \FF_n) h(\FF_n)}}(x)
\end{equation}
with $h(t) := \log \log(1/[t(1-t)])$, provides a trade-off 
between tail behavior and behavior in the center 
of the distribution. Previous proposals for the same 
purpose include \cite{Mason_Schuenemeyer_1983} 
and \cite{Revesz_1982}. But we shall demonstrate later 
that with purely multiplicative correction factors as in 
\eqref{eq:TnSP}, the tail regions are asymptotically 
underemphasized in comparison with the new methods presented here.

\subsection{Our new test statistics and confidence bands} 
\label{NewTestsNewCBs} 

To obtain a better compromise between the Kolmogorov--Smirnov and Berk--Jones tests, we propose a refined adjustment of $\FF_n(x)$ involving a pointwise standardization together with a pointwise additive correction, where the latter takes into account whether $x$ is in the center or in the tails of $F_0$ or $\FF_n$. Only after standardization and additive correction, we take a supremum over $x$. This approach of pointwise standardization plus additive correction before taking a supremum has been developed in the context of multi-scale testing and has proved quite successful there; see e.g.\ \cite{Duembgen_Spokoiny_2001}, \cite{Duembgen_Walther_2008}, \cite{Schmidt-Hieber_etal_2013} and \cite{Rohde_Duembgen_2013}. In the present setting, pointwise standardization means that we consider $n K_s(\FF_n(x),F_0(x))$, which behaves asymptotically like $\UU(F_0(x))^2/[2 F_0(x)(1 - F_0(x))]$ under the null hypothesis, that is, a squared standard Gaussian random variable times $0.5$. To identify an appropriate additive correction term, we utilize a refinement of the LIL \eqref{eq:LIL.BB}, based on Kolmogorov's upper class test (cf.\ \cite{Erdos_1942}, or \cite{Ito_McKean_1974}, Chapter~1.8). For $t \in (0,1)$ define
\begin{align*}
	C(t)
	&:= \log \log \frac{e}{4t(1-t)}
		= \log \bigl( 1 - \log(1 - (2t-1)^2) \bigr) \ \ge \ 0 , \\
	D(t)
	&:= \log(1 + C(t)^2) \in \bigl[ 0, \min\{C(t), C(t)^2\} \bigr] .
\end{align*}
Then for any fixed $\nu > 3/4$,
\begin{equation}
\label{eq:additive.LIL.BB}
	T_\nu := \sup_{t \in (0,1)}
		\Bigl( \frac{\mathbb{U}(t)^2}{2t(1-t)} - C_\nu(t) \Bigr)
	< \infty
\end{equation}
almost surely, where $C_\nu := C + \nu D$. Note that $C(t) = C(1 - t)$, $D(t) = D(1 - t)$, and, as $t \searrow 0$,
\begin{align*}
	C(t)
	&= \log \log(1/t) + O \bigl( (\log(1/t) )^{-1} \bigr) , \\
	D(t)
	&= 2 \log \log \log(1/t) + O \bigl( (\log \log(1/t))^{-1} \bigr) .
\end{align*}
This indicates why \eqref{eq:additive.LIL.BB} follows from Kolmogorov's test (see Section~S.1), and shows the connection between \eqref{eq:additive.LIL.BB} and \eqref{eq:LIL.BB}. On $(0,1/2]$, both functions $C$ and $D$ are decreasing with $C(1/2) = D(1/2) = 0$ and
\[
	\lim_{t \to 1/2} \frac{C(t)}{(2t-1)^2}
	\ = \ \lim_{t \to 1/2} \frac{D(t)}{(2t-1)^4}
	\ = \ 1 .
\]
Consequently, we propose the following test statistics:
\begin{equation}
\label{eq:GoF.statistic}
	T_{n,s,\nu}(F_0)
	:= \begin{cases}
		\displaystyle
		\sup_{x\colon 0 < F_0(x) < 1}
			\bigl[ n K_s(\FF_n(x),F_0(x)) - C_\nu(\FF_n(x),F_0(x)) \bigr]
				& \text{if} \ s > 0 , \\
		\displaystyle
		\sup_{x \in [X_{n:1},X_{n:n})}
			\bigl[ n K_s(\FF_n(x),F_0(x)) - C_\nu(\FF_n(x),F_0(x)) \bigr]
				& \text{if} \ s \le 0 ,
	\end{cases}
\end{equation}
where for $t,u \in [0,1]$,
\[
	C_\nu(u,t)
	:= \min_{\min(u,t) \le v \le \max(u,t)} C_\nu(v)
	= \begin{cases}
		C_\nu(\min(u,t)) & \text{if} \ \min(u,t) > 1/2 , \\
		C_\nu(\max(u,t)) & \text{if} \ \max(u,t) < 1/2 , \\
		0 & \text{else} ,
	\end{cases}
\]
with $C(0), C(1), D(0), D(1) := \infty$. As seen later, using this bivariate version $C_\nu(\FF_n(x),F_0(x))$ instead of $C_\nu(F_0(x))$ or $C_\nu(\FF_n(x))$ has computational advantages and increases power. The additive correction term $C_\nu(\FF_n(x),F_0(x))$ is large only if $x$ is far in the tails of $\FF_n$ \textsl{and} of $F_0$. 
\smallskip

The remainder of this paper is organized as follows. 

\begin{itemize} 
\item 
In Section~\ref{sec:LimitDistributionsUnif01} we show that under the null hypothesis, the test statistics $T_{n,s,\nu}(F_0)$ in \eqref{eq:GoF.statistic} converge in distribution to $T_\nu$ in \eqref{eq:additive.LIL.BB} for any fixed value of $s \in \RR$. 
\item 
Section~\ref{sec:Implications} discusses statistical implications of this finding. As explained in Section~\ref{subsec:GoF.tests}, goodness-of-fit tests based on  $T_{n,s,\nu}(F_0)$ have desirable asymptotic power. In particular, they are shown to attain a detection boundary of Ingster \cite{Ingster_1997} for Gaussian mixture models. 
Moreover, even under contiguous alternatives they have nontrivial asymptotic power, as opposed to goodness-of-fit tests based on $T_{n,s}^{\rm BJ}$ in \eqref{eq:SupFirst.S}. 
\item
In Section~\ref{subsec:ConfBands}, we analyze the confidence bands $(A_{n,s,\nu,\alpha}, B_{n,s,\nu,\alpha})$ resulting from inversion of the tests $T_{n,s,\nu}(\cdot)$. It will be shown that these bands have similar accuracy as those of Owen \cite{Owen_1995} and the bands $(A_{n,s,\alpha}^{\rm BJO}, B_{n,s,\alpha}^{\rm BJO})$ based on $T_{n,s}^{\rm BJ}(\cdot)$ in the tail regions while achieving the usual root-$n$ consistency everywhere. In addition, we compare our bands with the confidence bands of \cite{Stepanova_Pavlenko_2018}, confirming our claim that a purely multiplicative adjustment of $\FF_n - F_0$ is necessarily suboptimal in the tail regions.
\item 
Our results for the confidence bands elucidate the impact of the parameter $s$ on these bands for large sample sizes. These considerations are based on new inequalities and expansions for the divergences $K_s$ which are of independent interest.
\end{itemize}

All proofs and auxiliary results are deferred to Sections~\ref{sec:ProofsSec2}, \ref{sec:ProofsSec3} and an online supplement. References to the latter start with `S.' or `(S.'. Essential ingredients for the proofs in Section~\ref{sec:ProofsSec2} are tools and techniques of Cs{\"o}rg{\H{o}} et al.\ \cite{Csorgo_etal_1986}. A first version of this paper used a different, more self-contained approach which is probably of independent interest and outlined in Section~S.2. This also includes an alternative proof of \eqref{eq:additive.LIL.BB}.

\section{Limit distributions under the null hypothesis}
\label{sec:LimitDistributionsUnif01} 

Recall the uniform empirical process $\GG_n$ mentioned in the introduction. 
Under the null hypothesis that $F \equiv F_0$, the test statistic $T_{n,s,\nu}(F_0)$ has the same distribution as
\begin{equation}
\label{eq:GoF.statistic.0}
	T_{n,s,\nu} := \begin{cases}
		\displaystyle
		\sup_{t \in (0,1)}
			\bigl[ n K_s (\GG_n(t),t) - C_\nu(\GG_n(t),t) \bigr]
			& \text{if} \ s > 0 , \\
		\displaystyle
		\sup_{t \in [\xi_{n:1},\xi_{n:n})}
			\bigl[ n K_s (\GG_n(t),t) - C_\nu(\GG_n(t),t) \bigr]
			& \text{if} \ s \le 0 ,
	\end{cases}
\end{equation}
where $\xi_{n:1} < \cdots < \xi_{n:n}$ are the order statistics of the uniform sample $\xi_1,\ldots,\xi_n$. In particular, the $(1 - \alpha)$-quantile of $T_{n,s,\nu}(F_0)$ under the null hypothesis coincides with the $(1 - \alpha)$-quantile $\kappa_{n,s,\nu,\alpha}$ of $T_{n,s,\nu}$. Here is our main result for $T_{n,s,\nu}$ and $\kappa_{n,s,\nu,\alpha}$.

\begin{theorem}
\label{thm:PhiDivergenceNull}
For all $\nu > 3/4$ and $s \in \RR$,
\[
	T_{n,s,\nu} \rightarrow_d T_{\nu} .
\]
Moreover, $\kappa_{n,s,\nu,\alpha} \to \kappa_{\nu,\alpha} > 0$ for any fixed test level $\alpha \in (0,1)$, where $\kappa_{\nu,\alpha}$ is the $(1 - \alpha)$-quantile of $T_\nu$.
\end{theorem}

A key step along the way to proving Theorem~\ref{thm:PhiDivergenceNull} will be to consider the case $s=2$ and prove the following theorem for the uniform empirical process $\UU_n = \sqrt{n} (\GG_n - I)$, where $I$ denotes the distribution function of the uniform distribution on $[0,1]$.

\begin{theorem}
\label{IntermediateThmPhiDiv}
For all $\nu > 3/4$,
\[
	\tilde{T}_{n,\nu} := \sup_{t \in (0,1)}
		\left ( \frac{\UU_n(t)^2}{2 t(1-t)} - C_\nu(t) \right )
	\rightarrow_d T_{\nu} .
\]
\end{theorem}

\begin{remark}[The impact of $s$ and the definition of $T_{n,s,\nu}$]
Note that the parameter $s$ could be an arbitrary real number. However, numerical experiments indicate that the convergence to the asymptotic distribution is very slow if, say, $s < -0.5$ or $s > 1.5$. More precisely, Monte Carlo experiments show that for parameters $s \not \in [-0.5,1.5]$, the test statistcs $T_{n,s,\nu}$ are mainly influenced by just a few very small or very large order statistics. Moreover, if $s \in (0,0.5]$, one should redefine $T_{n,s,\nu}$ as a supremum over $[\xi_{n:1}, \xi_{n:n})$ rather than $(0,1)$. As shown in our proof of Theorem~\ref{thm:PhiDivergenceNull}, this modification does not alter the asymptotic distribution, but for realistic sample sizes $n$, taking the supremum over the full set $(0,1)$ for small parameters $s > 0$ leads to distributions which are mainly influenced by $\xi_{n:1}$.

Tables~S.1 and S.2 provide exact critical values $\kappa_{n,s,\nu,\alpha}$ for various sample sizes $n$, $s \in \{j/10: -10 \le j \le 20\}$, $\nu = 1$ and $\alpha = 0.5, 0.1, 0.05, 0.01$.

Similar discrepancies between asymptotic theory and finite sample behaviour can be observed for the Berk-Jones quantiles $\kappa_{n,s,\alpha}^{\rm BJ}$ if $s \not\in [-0.5,1.5]$, see Tables~S.3 and S.4.
\end{remark}

\section{Statistical implications}
\label{sec:Implications}

\subsection{Goodness-of-fit tests}
\label{subsec:GoF.tests}

As explained in the introduction, we can reject the null hypothesis that $F$ is a given continuous distribution function $F_0$ at level $\alpha$ if the test statistic $T_{n,s,\nu}(F_0)$, defined in \eqref{eq:GoF.statistic}, exceeds the $(1 - \alpha)$-quantile $\kappa_{n,s,\nu,\alpha}$ of $T_{n,s,\nu}$. The test statistics $T_{n,s,\nu}$ and $T_{n,s,\nu}(F_0)$ can be represented as the maximum of at most $2n$ terms: with $u_{n,i} := i/n$, the statistic $T_{n,s,\nu}$ equals
\[
	\max_{1 \le i \le n} \max \bigl\{
		n K_s(u_{n,i-1},\xi_{n:i}) - C_\nu(u_{n,i-1},\xi_{n:i}),
		n K_s(u_{n,i},\xi_{n:i}) - C_\nu(u_{n,i},\xi_{n:i}) \bigr\}
\]
if $s > 0$, and
\[
	\max_{1 \le i < n} \max \bigl\{
		n K_s(u_{n,i},\xi_{n:i}) - C_\nu(u_{n,i},\xi_{n:i}),
		n K_s(u_{n,i},\xi_{n:i+1}) - C_\nu(u_{n,i},\xi_{n:i+1}) \bigr\}
\]
if $s \le 0$. The statistic $T_{n,s,\nu}(F_0)$ can be represented analogously with $F_0(X_{n:i})$ in place of $\xi_{n:i}$. These formulae follow from the fact that for fixed $u \in (0,1)$, the function $t \mapsto n K_s(u,t) - C_\nu(u,t)$ is continuous on $(0,1)$, increasing on $[u,1)$ and decreasing on $(0,u]$. For $K_s(u,t) = K_{1-s}(t,u)$ is convex in $t$ with minimum at $t = u$, see (S.12) in Section~S.3, and $C_\nu(u,t)$ is increasing in $t \in (0,u]$ and decreasing in $t \in [u,1)$. If $s > 0$, these monotonicities are also true for $u \in \{0,1\}$, precisely,
\[
	C_\nu(0,t) = C_\nu(\min(t,1/2))
	\quad\text{and}\quad
	K_s(0,t) = \begin{cases}
		- \log(1 - t) & \text{if} \ s = 1 , \\
		((1 - t)^{1-s} - 1)/(s(s-1)) & \text{if} \ s \ne 1 ,
	\end{cases}
\]
while $C_\nu(1,t) = C_\nu(0,1-t)$ and $K_s(1,t) = K_s(0,1-t)$.

\subsubsection{Non-contiguous alternatives}

Now suppose that the true distribution function of the observations $X_i$ is a continuous distribution function $F_n$ such that $\{x \in \RR : 0 < F_n(x) < 1\} \subset \{x \in \RR : 0 < F_0(x) < 1\}$. A first question is: under what conditions on the sequence $(F_n)_n$ does our goodness-of-fit test have asymptotic power one for any fixed test level $\alpha \in (0,1)$. Since $\kappa_{n, s,\nu,\alpha} \rightarrow \kappa_{\nu,\alpha} < \infty$, this goal is equivalent to 
\begin{eqnarray}
	\label{PowerConsistencyGOF}
	P_{F_n}(T_{n,s,\nu}(F_0) > \kappa) \rightarrow 1 \quad \text{for any fixed} \ \kappa > 0 .
\end{eqnarray}
To verify this property, the following function $\Delta_n : \RR \to [0,\infty)$ plays a key role:
\[
	\Delta_n := \frac{\sqrt{n} |F_n - F_0|}{\min\{H_n(F_n),H_n(F_0)\}}
	\quad\text{with}\quad
	H_n(t) := \sqrt{(1 + C(t)) t(1 - t)} + \frac{1 + C(t)}{\sqrt{n}}
\]
for $t \in [0,1]$ with the conventions $C(t) := \infty$ and $C(t) t(1 - t) := 0$ for $t \in \{0,1\}$.

\begin{theorem}
\label{thm:GOF-PwrOneCondition}
Suppose that the sequence $(F_n)_n$ satisfies the condition
\begin{equation}
\label{eq:GOF-PwrOneCondition}
	\sup_{x \in \RR} \Delta_n(x) \to \infty .
\end{equation} 
Then \eqref{PowerConsistencyGOF} holds true for any $s \in [-1,2]$.
\end{theorem}

It follows immediately from this theorem that \eqref{PowerConsistencyGOF} is satisfied whenever $F_n \equiv F_*$ 
for all sample sizes $n$, where $F_* \neq F_0$.

As a litmus test for our procedures and Theorem~\ref{thm:GOF-PwrOneCondition}, we consider a testing problem studied in detail by \cite{Ingster_1997}. The null hypothesis is given by $F_0 = \Phi$, the standard Gaussian distribution function, whereas
\[
	F_n(x) := (1 - \epsilon_n) \Phi(x) + \epsilon_n \Phi (x - \mu_n) .
\]
for certain numbers $\epsilon_n \in (0,1)$ and $\mu_n > 0$. By means of Theorem~\ref{thm:GOF-PwrOneCondition} one can derive the following result.

\begin{corollary}
\label{cor:GaussMixtureConsistencyCriterion} 
\textbf{(a)}  Suppose that $\epsilon_n = n^{-\beta + o(1)}$ for some fixed $\beta \in (1/2,1)$.  
Furthermore let $\mu_n = \sqrt{2r \log n}$ for some $r \in (0,1)$.  
Then \eqref{PowerConsistencyGOF} is satisfied for any $s \in [-1,2]$ if 
\[
	r > \begin{cases}
		\beta - 1/2
			& \text{if} \ \beta \in (1/2, 3/4] , \\
		(1 - \sqrt{1-\beta})^2
			& \text{if} \ \beta \in [3/4, 1) .
	\end{cases}
\]
\textbf{(b)}  Suppose that $\epsilon_n = n^{-1/2 + o(1)} $ such that $\pi_n := \sqrt{n} \epsilon_n \rightarrow 0$.
Then \eqref{PowerConsistencyGOF} is satisfied for any $s \in [-1,2]$ if $\mu_n = \sqrt{ 2 \lambda \log (1 / \pi_n )} $ for some $\lambda > 1$.
\end{corollary} 

As explained by \cite{Ingster_1997}, any goodness-of-fit test at fixed level $\alpha \in (0,1)$ has trivial asymptotic power $\alpha$ whenever $\epsilon_n = n^{-\beta}$ for some $\beta \in (1/2,1)$ and $\mu_n = \sqrt{2 r \log n}$ with
\[
	r < \begin{cases}
		\beta-1/2
			& \text{if} \ \beta \in (1/2, 3/4], \\
		(1 - \sqrt{1-\beta})^2
			& \text{if} \ \beta \in [3/4, 1) .
	\end{cases}
\]
Thus part~(a) of the previous corollary shows that our new family of tests achieves this detection boundary, as do the goodness-of-fit tests of \cite{Donoho_Jin_2004}, \cite{Jager_Wellner_2007} and \cite{Gontscharuk_etal_2016}.

A connection between parts~(a) and (b) of Corollary~\ref{cor:GaussMixtureConsistencyCriterion} can be seen as follows: let $\epsilon_n = n^{-\beta}$ for some fixed $\beta \in (1/2, 3/4]$, and $\mu_n = \sqrt{2r \log (n)}$ for some $r > \beta-1/2$. Then $r = \lambda (\beta - 1/2)$ for some $\lambda > 1$, and with $\pi_n = \sqrt{n} \epsilon_n = n^{1/2 - \beta}$, we may write $\sqrt{2r \log(n)} = \sqrt{2 \lambda \log(1/\pi_n)}$.

\subsubsection{Contiguous alternatives}

Suppose that the distribution functions $F_0$ and $F_n$ have densities $f_0$ and $f_n$, respectively, with respect to some continuous measure $\Lambda$ on $\RR$ such that for some function $a$,
\begin{equation}
\label{eq:contiguous.alternative}
	\sqrt{n} (f_n^{1/2} - f_0^{1/2}) \rightarrow 2^{-1} a f_0^{1/2}
	\quad\text{in} \ L_2(\Lambda) .
\end{equation}
Then it follows easily that $a \in L_2(F_0)$, $\int a \, dF_0 = 0$ and
\[
	\sqrt{n} (F_n - F_0 )(t) \to A(t) := \int_{-\infty}^t  a \, dF_0
	\quad\text{uniformly in } \ t \in \RR .
\]
Furthermore, since $\int_{-\infty}^t a \, dF_0 = \int_{\RR} (1_{[x \le t]} - F_0(t)) a(x) \, d F_0(x)$, the Cauchy-Schwarz inequality yields that
\begin{equation}
	|A(t)| \leq \sqrt{F_0(t) (1 - F_0(t))} \, \| a \|_{L_2 (F_0)} .
\label{TheCapABound}
\end{equation}

\begin{lemma}[Power of ``tail-dominated'' tests under contiguous alternatives]
\label{lem:power.contiguous.1}
Let $(\varphi_n)_n$ be a sequence of tests with the following two properties:

\noindent
\textbf{(i)} For a fixed level $\alpha \in (0,1)$,
\[
	E_{F_0} \varphi_n(X_1,\ldots,X_n) \to \alpha .
\]

\noindent
\textbf{(ii)} For any fixed $0 < \rho < 1/2$ and $x_{\rho} := F_0^{-1} (\rho)$, $y_{\rho} := F_0^{-1} (1-\rho)$, there exists a test $\varphi_{n,\rho}$ depending only on $(\FF_n(x))_{x \not\in [x_\rho,y_\rho]}$ such that
\[
	P_{F_0}(\varphi_n \ne \varphi_{n,\rho}) \to 0 .
\]

\noindent
Then under assumption \eqref{eq:contiguous.alternative},
\[
	\limsup_{n \to \infty} E_{F_n} \varphi_n(X_1,\ldots,X_n)
	\le \alpha .
\]
\end{lemma}

Note that the Berk-Jones tests with $T_{n,s}^{\rm BJ}(F_0)$ satisfy the assumptions of Lemma~\ref{lem:power.contiguous.1}, if tuned to have asymptotic level $\alpha$. For all of them involve a test statistic of the type
\[
	T_n(F_0) = \sup_{x \in \RR} \Gamma_n(\FF_n(x))
\]
with a function $\Gamma_n : \RR \to [0,\infty]$ such that under the null hypothesis,
\[
	\sup_{x \in \RR} \Gamma_n(\FF_n(x)) \to_p \infty ,
\]
but for any $0 < \rho < 1/2$,
\[
	\sup_{x \in [x_{\rho},y_{\rho}]} \Gamma_n(\FF_n(x)) = O_p(1) .
\]
Hence $T_n(F_0)$ equals
\[
	T_n^{(\rho)}(F_0) := \sup_{x \not\in [x_{\rho},y_{\rho}]} \Gamma_n(\FF_n(x))
\]
with asymptotic probability one. Thus we may replace the test statistic $T_n(F_0)$ with $T_n^{(\rho)}(F_0)$ while keeping the critical value.

By way of contrast, the goodness-of-fit test based on $T_{n,s,\nu}(F_0)$ has nontrivial asymptotic power in the present setting.

\begin{theorem}[Power of new tests under contiguous alternatives]
\label{thm:power.contiguous.2}
In the setting~\eqref{eq:contiguous.alternative}, the test statistic $T_{n,s,\nu}(F_0)$ converges in distribution to
\[
	T_\nu(A) := \sup_{t \in (0,1)}
		\Bigl( \frac{\bigl( \UU(t) + A(F_0^{-1}(t)) \bigr)^2}{2t(1-t)} - C_\nu(t) \Bigr) .
\]
In particular,
\[
	P_{F_n} \bigl[ T_{n,s,\nu}(F_0) \ge \kappa_{n,s,\nu,\alpha} \bigr]
	\to P[T_\nu(A) \ge \kappa_{\nu,\alpha}] \ge \alpha .
\]
Concerning the impact of $A$,
\[
	P[T_\nu(A) \ge \kappa_{\nu,\alpha}] \to 1
	\quad\text{as}\quad
	\sup_{t \in (0,1)} \Bigl( \frac{|A(F_0^{-1}(t))|}{\sqrt{2 t(1 - t)}} - \sqrt{C(t)} \Bigr)
		\to \infty .
\]
\end{theorem}

\subsection{Confidence bands}
\label{subsec:ConfBands}

The confidence bands of \cite{Owen_1995}, defined in terms of $K = K_1$, may be generalized to arbitrary fixed $s \in [-1,2]$, but we restrict our attention to $s \in (0,2]$, because for $s \le 0$ and a large range of sample sizes $n$, the resulting bands would focus mainly on small regions in the tails and be rather wide elsewhere. With confidence $1 - \alpha$ we may claim that $\sup_{x\colon 0 < F(x) < 1} n K_s(\FF_n(x),F(x))$ does not exceed the $(1 - \alpha)$-quantile $\kappa_{n,s,\alpha}^{\rm BJ}$ of $\sup_{t\in (0,1)} n K_s(\GG_n(t),t)$. As explained in Section~S.6, inverting the inequality $nK_s(\FF_n(x),F(x)) \le \kappa_{n,s,\alpha}^{\rm BJ}$ for fixed $x$ with respect to $F(x)$ reveals that for $0\le i \le n$ and $X_{n:i} \le x < X_{n:i+1}$,
\[
	F(x) \in \bigl[ A_{n,s,\alpha}^{\rm BJO}(x), B_{n,s,\alpha}^{\rm BJO}(x) \bigr]
	= [a_{n,s,\alpha,i}^{\rm BJO}, b_{n,s,\alpha,i}^{\rm BJO}] ,
\]
where $a_{n,s,\alpha,i}^{\rm BJO} \le u_{n,i} \le b_{n,s,\alpha,i}^{\rm BJO}$ are given by $a_{n,s,\alpha,0}^{\rm BJO} := 0$, $b_{n,s,\alpha,n}^{\rm BJO} := 1$ and for $0 \le i < n$,
\begin{align*}
	b_{n,s,\alpha,i}^{\rm BJO}
	&:= \max \bigl\{ t \in (u_{n,i}, 1]:
		n K_s(u_{n,i}, t) \le \kappa_{n,s,\alpha}^{\rm BJ} \bigr\} , \\
	a_{n,s,\alpha,n-i}^{\rm BJO}
	&:= 1 - b_{n,s,\alpha,i}^{\rm BJO} .
\end{align*}
Thus, computing the confidence band $(A_{n,s,\alpha}^{\rm BJO}, B_{n,s,\alpha}^{\rm BJO})$ boils down to determining the $2(n+1)$ numbers $a_{n,s,\alpha,i}^{\rm BJO}$ and $b_{n,s,\alpha,i}^{\rm BJO}$, $0 \le i \le n$.

Our new method is analogous: with confidence $1-\alpha$, for $0 \le i \le n$ 
and $X_{n:i} \le x < X_{n:i+1}$, the value $F(x)$ is contained in
\[
	\bigl[ A_{n,s,\nu,\alpha}(x), B_{n,s,\nu,\alpha}(x) \bigr]
	= [a_{n,s,\nu,\alpha,i}, b_{n,s,\nu,\alpha,i}] ,
\]
where $a_{n,s,\nu,\alpha,0} := 0$, $b_{n,s,\nu,\alpha,n} := 1$ and for $0 \le i < n$,
\begin{align*}
	b_{n,s,\nu,\alpha,i}
	&:= \max \bigl\{ t \in (u_{n,i} , 1] :
		n K(u_{n,i}, t) - C_\nu(u_{n,i},t) \le \kappa_{n,s,\nu,\alpha} \bigr\} , \\
	a_{n,s,\nu,\alpha,n-i}
	&:= 1 - b_{n,s,\nu,\alpha,i} .
\end{align*}

\begin{figure}
\centering
\includegraphics[width=0.7\textwidth]{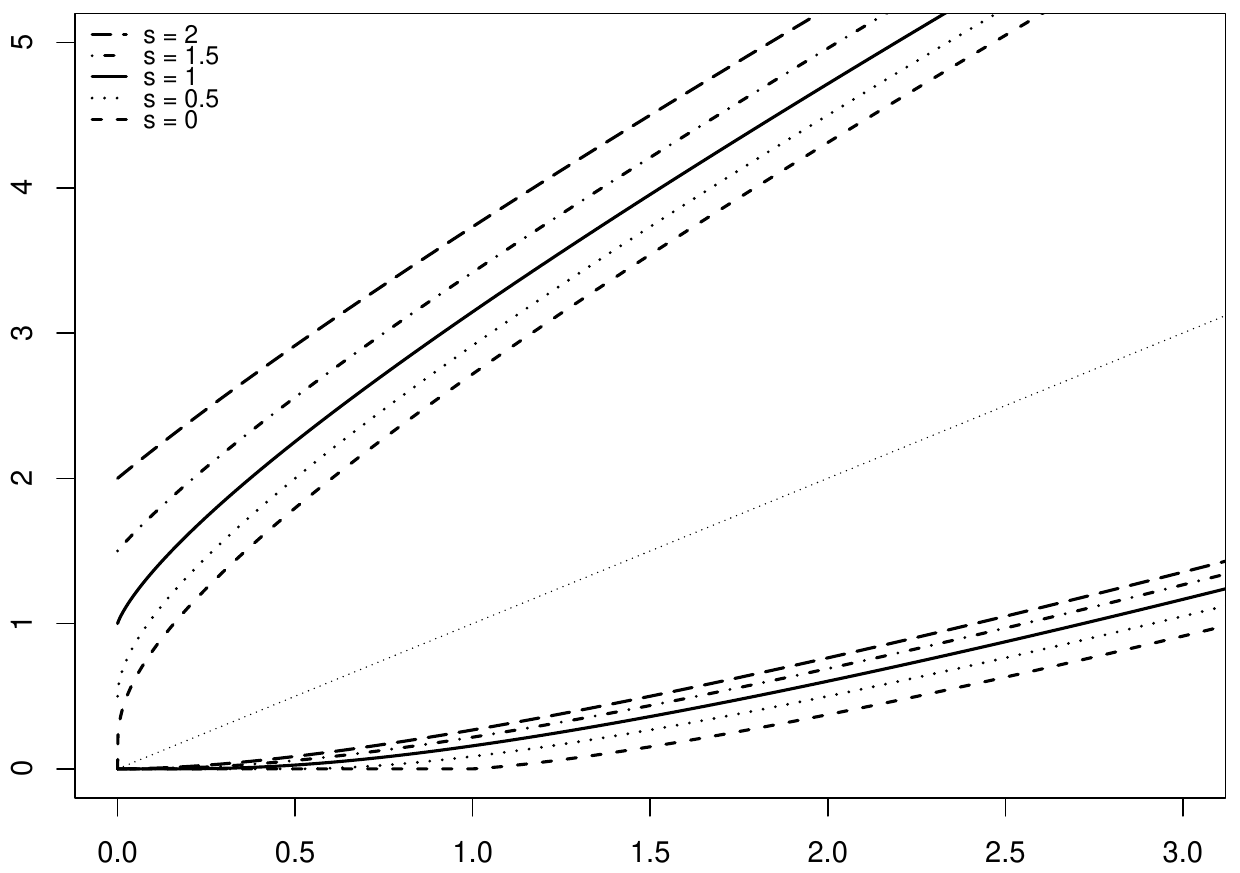}
\caption{The auxiliary functions $a_s$ (below diagonal), $b_s$ (above diagonal) for $s \in \{0, 0.5, 1, 1.5, 2\}$.}
\label{fig:AB}
\end{figure}

To understand the asymptotic performance of these confidence bands properly, we need auxiliary functions $a_s, b_s : [0,\infty) \to [0,\infty)$. Note first that for any $s \in [-1,2]$, $K_s(u,t)$ in \eqref{TheKFunctions} may be represented as
\begin{equation}
\label{eq:Ks.via.phis}
	K_s(u,t)
	= t \phi_s (u/t) + (1-t) \phi_s [(1-u)/(1-t)]
\end{equation}
where
\begin{equation}
	\label{ThePhiFunctions}
	\phi_s(x) = \begin{cases}
		(x^s - sx + s - 1)/[s(s-1)],
			& s \neq 0,1 , \\
		x \log x - x + 1,
			& s = 1, \\
		x - 1 - \log x,
			& s = 0 ,
	\end{cases}
\end{equation}
for $x \in (0,\infty)$, and $\phi_s(0) := \lim_{x \searrow 0} \phi_s(x)$ equals $1/s^+$. If $u$ and $t$ are close to $0$, one may approximate $K_s(u,t)$ by
\[
	H_s(u,t) := t \phi_s(u/t) .
\]
The properties of $H_s : [0,\infty) \times (0,\infty) \to [0,\infty]$ are treated in Lemma~S.13. In particular, it is shown that
\begin{align*}
	a_s(x) &:= \begin{cases}
		0
			& \text{if} \ x = 0 , \\
		\inf \{ y \in (0,x) : H_s(x,y) \le 1 \}
			& \text{else} ,
		\end{cases} \\
	b_s(x) &:= \begin{cases}
		s^+
			& \text{if} \ x = 0 , \\
		\max\{ y > x : H_s(x,y) \le 1 \}
			& \text{else} ,
		\end{cases}
\end{align*}
define continuous functions $a_s, b_s : [0,\infty) \to [0,\infty)$, where $a_s$ is convex with $a_s(0) = 0 = a_s'(0)$, $a_s(x) = 0$ if and only if $x \le (1 - s)^+$, and $b_s$ is concave. Moreover, $a_s(x) = x - \sqrt{2x} + O(1)$ and $b_s(x) = x + \sqrt{2x} + O(1)$ as $x \to \infty$. Finally, for fixed $x > 0$, $a_s(x)$ and $b_s(x)$ are non-decreasing in $s \in [-1,2]$ with $a_s(x) < x < b_s(x)$. Figure~\ref{fig:AB} depicts these functions $a_s, b_s$ on the interval $[0,3]$ for $s \in \{0, 0.5, 1, 1.5, 2\}$.

Our first result shows that the confidence bands $(A_{n,s,\alpha}^{\rm BJO}, B_{n,s,\alpha}^{\rm BJO})$ and $(A_{n,s,\nu,\alpha}, B_{n,s,\nu,\alpha})$ are asymptotically equivalent in the tail regions, that is, for $\FF_n(x)$ close to zero or close to one. Moreover, the test level $\alpha$ is asymptotically irrelevant there, but the parameter $s$ does play a role when $\min\{\FF_n(x), 1 - \FF_n(x)\} \le O(n^{-1} \log\log n)$.

\begin{theorem}
\label{thm:conf.bands.tails}
Let $\gamma_n := n^{-1} \log\log n$. For any fixed $s \in (0,2]$, $\nu > 3/4$ and $\delta \in (0,1)$,
\begin{align*}
	\left. \begin{array}{c}
		u_{n,i} - a_{n,s,\alpha,i}^{\rm BJO} \\
		u_{n,i} - a_{n,s,\nu,\alpha,i} \\[0.5ex]
		b_{n,s,\alpha,n-i}^{\rm BJO} - u_{n,n-i} \\
		b_{n,s,\nu,\alpha,n-i} - u_{n,n-i} \\
	\end{array} \right\}
	&= \gamma_n \bigl( i/\log\log n - a_s(i / \log\log n) \bigr) (1 + o(1))
\intertext{and}
	\left. \begin{array}{c}
		b_{n,s,\alpha,i}^{\rm BJO} - u_{n,i} \\
		b_{n,s,\nu,\alpha,i} - u_{n,i} \\[0.5ex]
		u_{n,n-i} - a_{n,s,\alpha,n-i}^{\rm BJO} \\
		u_{n,n-i} - a_{n,s,\nu,\alpha,n-i} \\
	\end{array} \right\}
	&= \gamma_n \bigl( b_s(i/\log\log n) - i/\log\log n \bigr) (1 + o(1)) ,
\end{align*}
uniformly in $i \in \{0,1,\ldots,n\} \cap [0, n^\delta]$.
\end{theorem}

\begin{remark}[Choice of $s$]
Concerning the choice of $s$, Theorem~\ref{thm:conf.bands.tails} shows that smaller (resp.\ larger) values of $s$ lead to better upper (resp.\ lower) and worse lower (resp. upper) bounds for $F(x)$ in the left tail and better lower (resp.\ upper) and worse upper (resp.\ lower bounds) for $F(x)$ in the right tail. The choice $s = 1$ seems to be a good compromise, see also the numerical examples later.
\end{remark}

The next result shows that in the central region, the parameter $s$ is asymptotically irrelevant, and the width of the band $(A_{n,s,\nu,\alpha}, B_{n,s,\nu,\alpha})$ is of smaller order than the width of $(A_{n,s,\alpha}^{\rm BJO}, B_{n,s,\alpha}^{\rm BJO})$.

\begin{theorem}
\label{thm:conf.bands.center}
For any fixed $s \in (0,2]$, $\nu > 3/4$ and $\delta \in (0,1)$,
\begin{align*}
	\left. \begin{array}{c}
		u_{n,i} - a_{n,s,\alpha,i}^{\rm BJO} \\[0.5ex]
		b_{n,s,\alpha,i}^{\rm BJO} - u_{n,i}
	\end{array} \right\}
	&= \sqrt{ 2 \gamma_n \, u_{n,i}(1 - u_{n,i}) } \,
		(1 + o(1)) , \\
	\left. \begin{array}{c}
		u_{n,i} - a_{n,s,\nu,\alpha,i} \\
		b_{n,s,\nu,\alpha,i} - u_{n,i}
	\end{array} \right\}
	&= \sqrt{ 2 \gamma_{n,\nu,\alpha}(u_{n,i}) \, u_{n,i}(1 - u_{n,i}) } \,
		(1 + o(1)) ,
\end{align*}
uniformly in $i \in \{0,1,\ldots,n\} \cap [n^\delta, n - n^\delta]$, where $\gamma_n = n^{-1} \log\log n$ and $\gamma_{n,\nu,\alpha}(u) := n^{-1} \bigl( C_{\nu}(u) + \kappa_{\nu,\alpha} \bigr)$.
\end{theorem}

Note that $(C_\nu(u) + \kappa_{\nu,\alpha}) u(1 - u) \to 0$ as $u \to \{0,1\}$. Thus one can deduce from Theorems~\ref{thm:conf.bands.tails} and \ref{thm:conf.bands.center} that
\begin{align*}
	\max_{i=0,1, \ldots, n} (b_{n,i}^{\rm BJO} - u_{n,i})
	= \max_{i=0,1, \ldots, n}  (u_{n,i} - a _{n,i}^{\rm BJO})
	&= \sqrt{\gamma_n/2} (1 + o(1)) , \\
	\max_{i=0,1, \ldots, n} (b_{n,i} - u_{n,i})
	= \max_{i=0,1, \ldots , n}  (u_{n,i} - a_{n,i})
	&= O(n^{-1/2}) .
\end{align*}

\begin{remark}[Comparison with Stepanova--Pavlenko \cite{Stepanova_Pavlenko_2018}]
\label{rem:Comparison.with.SP2018}
The confidence band $(A_{n,\alpha}^{\rm SP}, B_{n,\alpha}^{\rm SP})$ of \cite{Stepanova_Pavlenko_2018} with the test statistic $T_n^{\rm SP}(\cdot)$ in \eqref{eq:TnSP} can be represented as follows: for $0 \le i \le n$ and $X_{n:i} \le x < X_{n:i+1}$,
\[
	\bigl[ A_{n,\alpha}^{\rm SP}(x), B_{n,\alpha}^{\rm SP}(x) \bigr]
	= [a_{n,\alpha,i}^{\rm SP}, b_{n,\alpha,i}^{\rm SP}] ,
\]
where $a_{n,\alpha,0}^{\rm SP} = 0$, $b_{n,\alpha,0}^{\rm SP} = b_{n,\alpha,1}^{\rm SP}$, $a_{n,\alpha,n}^{\rm SP} = a_{n,\alpha,n-1}^{\rm SP}$, $b_{n,\alpha,n}^{\rm SP} = 1$, and for $1 \le i < n$,
\[
	[a_{n,\alpha,i}^{\rm SP}, b_{n,\alpha,i}^{\rm SP}]
	= \Bigl[ u_{n,i} \pm n^{-1/2} \kappa_{n,\alpha}^{\rm SP}
		\sqrt{ u_{ni}(1 - u_{n,i}) h(u_{n,i})} \Bigr]
		\cap [0,1] .
\]
Recall that $h(t) = \log\log(1/[t(1-t)])$. 
Here $\kappa_{n,\alpha}^{\rm SP}$ is the $(1 - \alpha)$-quantile of 
$T_{n,\alpha}^{\rm SP}(F_0)$ in case of $F \equiv F_0$, and it 
converges to the $(1 - \alpha)$-quantile $\kappa_\alpha^{\rm SP}$ of
\[
	\sup_{t \in (0,1)} \frac{|\UU(t)|}{\sqrt{t(1 - t) h(t)}} .
\]
Consequently, for fixed $s \in (0,2]$, $\nu > 3/4$ and $\delta \in (0,1)$,
\[
	\frac{b_{n,\alpha,i}^{\rm SP} - u_{n,i}}{b_{n,s,\nu,\alpha,i} - u_{n,i}} ,
	\frac{u_{n,i} - a_{n,\alpha,i}^{\rm SP}}{u_{n,i} - a_{n,s,\nu,\alpha,i}}
	= \frac{\kappa_\alpha^{\rm SP}\sqrt{h(u_{n,i})}}
		   {\sqrt{2 (C_\nu(u_{n,i}) + \kappa_{\nu,\alpha})}} (1 + o(1))
\]
uniformly in $i \in \{0,1,\ldots,n\} \cap [n^\delta, n - n^\delta]$. But
\[
	\lim_{u \to \{0,1\}}
	\frac{\kappa_\alpha^{\rm SP}\sqrt{h(u)}}
		 {\sqrt{2 (C_\nu(u) + \kappa_{\nu,\alpha})}}
	 = \frac{\kappa_\alpha^{\rm SP}}{\sqrt{2}}
	 \begin{cases}
	 	\ge 1 , \\
		\to \infty & \text{as} \ \alpha \searrow 0 ,
	\end{cases}
\]
because $h(t) / \log\log(1/t)$ and $C_\nu(t) / \log\log(1/t)$ converge to $1$ as $t \searrow 0$. Thus, the confidence band $(A_{n,\alpha}^{\rm SP}, B_{n,\alpha}^{\rm SP})$ is asymptotically wider than $(A_{n,s,\nu,\alpha}, B_{n,s,\nu,\alpha})$ in the tail regions for sufficiently small $\alpha$.

Note that these considerations apply to any choice of the continuous function $h : (0,1) \to (0,\infty)$ in \eqref{eq:TnSP} as long as $h(t) / \log\log(1/t) \to 1$ as $t \searrow 0$.
\end{remark}

\begin{remark}[Bahadur and Savage \cite{Bahadur_Savage_1956} revisited]
On $(-\infty,X_{n:1}]$, the upper confidence bounds for $F$ are constant $b_{n,s,\alpha,1}^{\rm BJO}$ or $b_{n,s,\nu,\alpha,1}$, and this is of order $O(n^{-1} \log\log n)$. Likewise, on $(X_{n:n},\infty)$, the lower confidence bounds for $F$ are constant $1 - b_{n,s,\alpha,1}^{\rm BJO}$ or $1 - b_{n,s,\nu,\alpha,1}$. Interestingly, for any $(1 - \alpha)$-confidence band for a continuous distribution function $F$, the upper bound has to be greater than $c/n$ with asymptotic probability at least $e^c \alpha$, and the lower bound has to be smaller than $1 - c/n$ with asymptotic probability at least $e^c \alpha$. This follows from a quantitative version of Theorem~2 of \cite{Bahadur_Savage_1956}, stated as Theorem~\ref{thm:Bahadur-Savage} below.

It is also instructive to consider Daniels' lower confidence bound for a continuous distribution function $F$, namely
\[
	P_F ( \alpha \FF_n (x) \le F(x) \ \text{for all} \ x \in \RR) = 1 - \alpha .
\]
\end{remark}

\begin{theorem}
\label{thm:Bahadur-Savage}
Let $\mathcal{F}$ be a family of continuous distribution functions which is convex and closed under translations, that is, $F(\cdot - \mu) \in \mathcal{F}$ for all $F \in \mathcal{F}$ and $\mu \in \RR$. Let $(A_n, B_n)$ be a $(1-\alpha)$-confidence band for $F \in \mathcal{F}$. Then for any $F \in \mathcal{F}$ and $\epsilon \in (0,1)$,
\[
	P_F \Bigl( \inf_{x \in \RR} B_n (x) < \epsilon \Bigr) \le (1 - \eps)^{-n} \alpha
	\quad\text{and}\quad
	P_F \Bigl( \sup_{x \in \RR} A_n (x) > 1- \epsilon \Bigr)
	\le (1-\epsilon)^{-n} \alpha .
\]
\end{theorem}

In our context, $\mathcal{F}$ would be the family of all continuous distribution functions. But the precision bounds in Theorem~\ref{thm:Bahadur-Savage} apply to much smaller families $\mathcal{F}$ already, for instance, the family of all convex combinations of $F_o(\cdot - \mu)$, $\mu \in \RR$, where $F_o$ is an arbitrary continuous distribution function. For the reader's convenience, a proof of Theorem~\ref{thm:Bahadur-Savage} is provided in Section~S.5.

\begin{figure}
\includegraphics[width=0.49\textwidth]{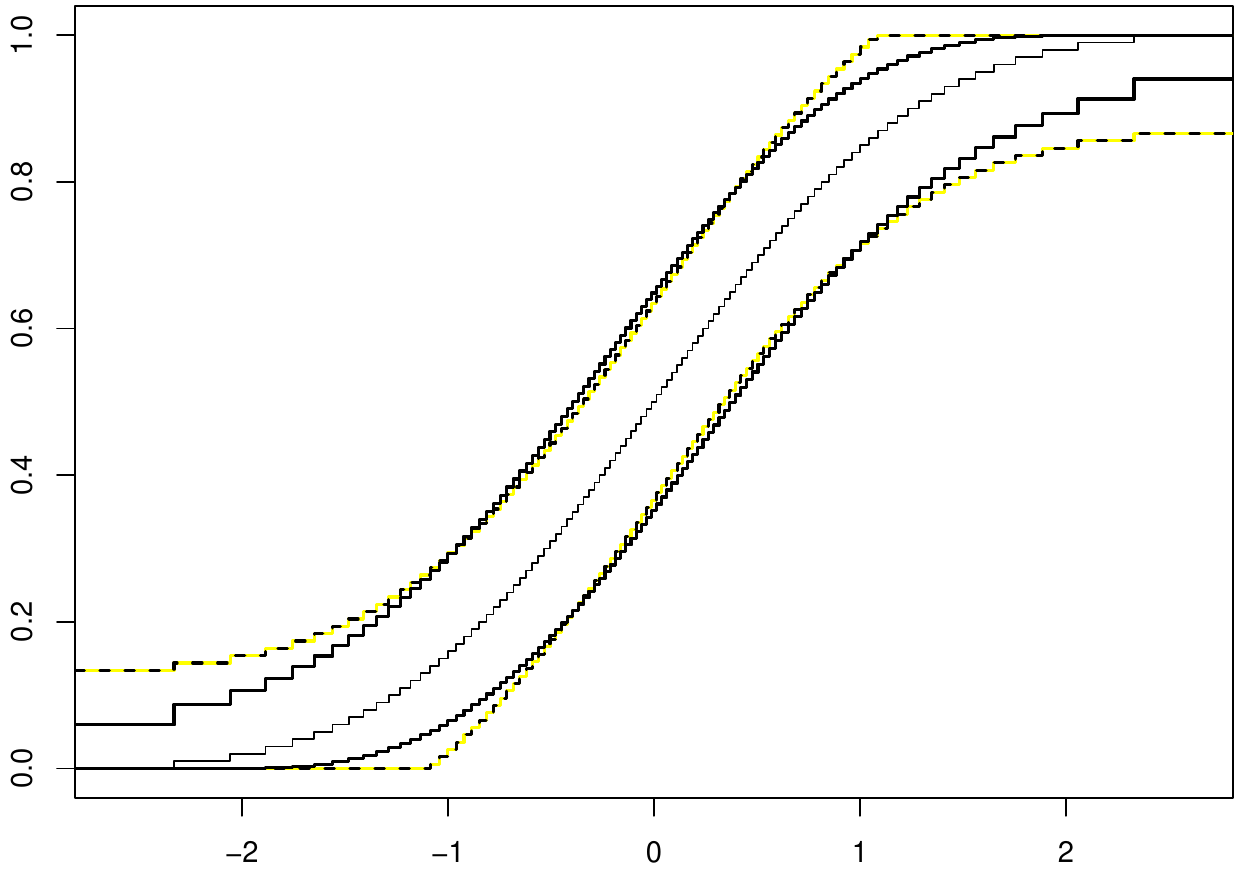}
\hfill
\includegraphics[width=0.49\textwidth]{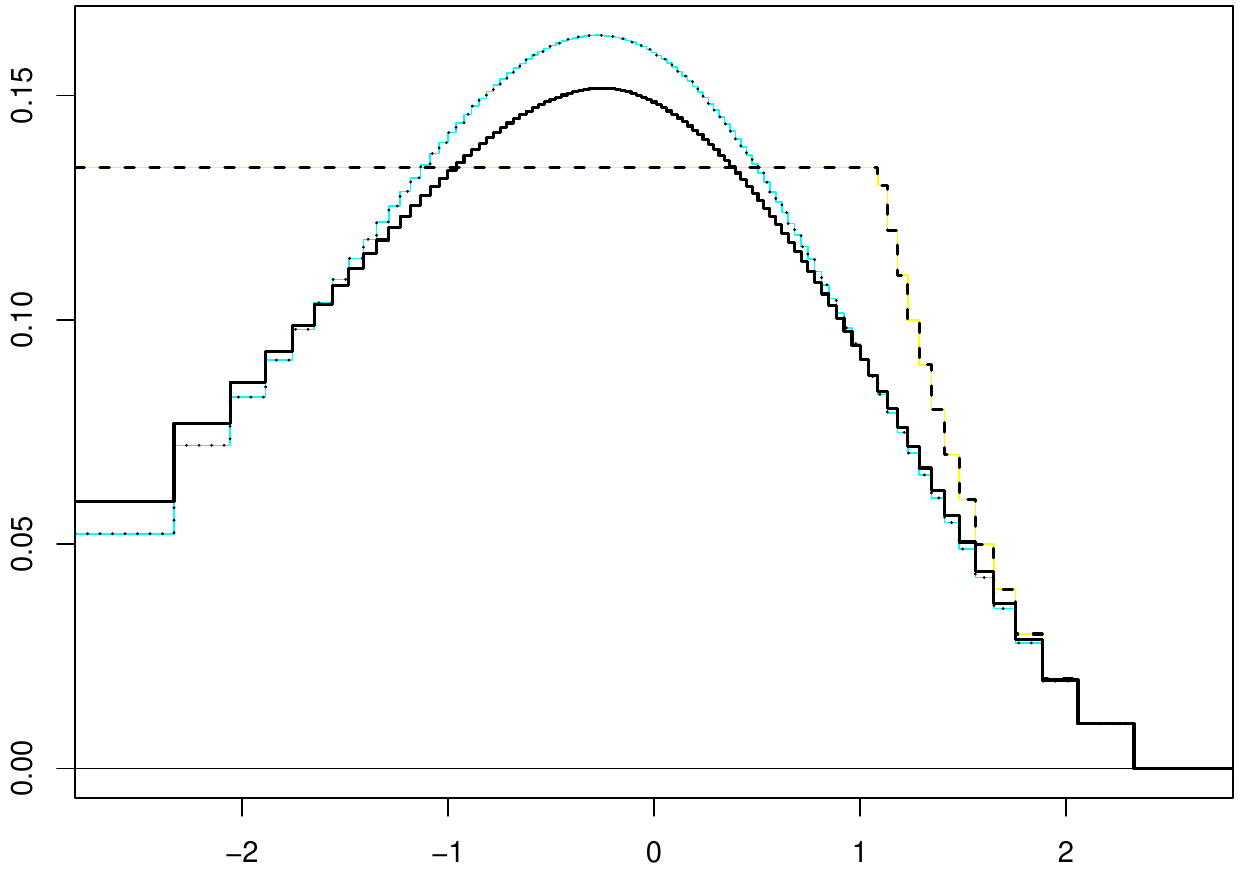}
\caption{$95\%$-confidence bands for $n = 100$. Left panel: $(A_{n,1,1,\alpha}, B_{n,1,1,\alpha})$ (solid) and $(A_{n,\alpha}^{\rm KS}, B_{n,\alpha}^{\rm KS})$ (dashed). Right panel: centered upper bounds $B_{n,1,1,\alpha} - \FF_{n}$ (solid), $B_{n,1,\alpha}^{\rm BJO} - \FF_{n}$ (dotted) and $B_{n,\alpha}^{\rm KS} - \FF_{n}$ (dashed).}
\label{fig:AsEfficiency.1}
\end{figure}

\begin{figure}
\includegraphics[width=0.49\textwidth]{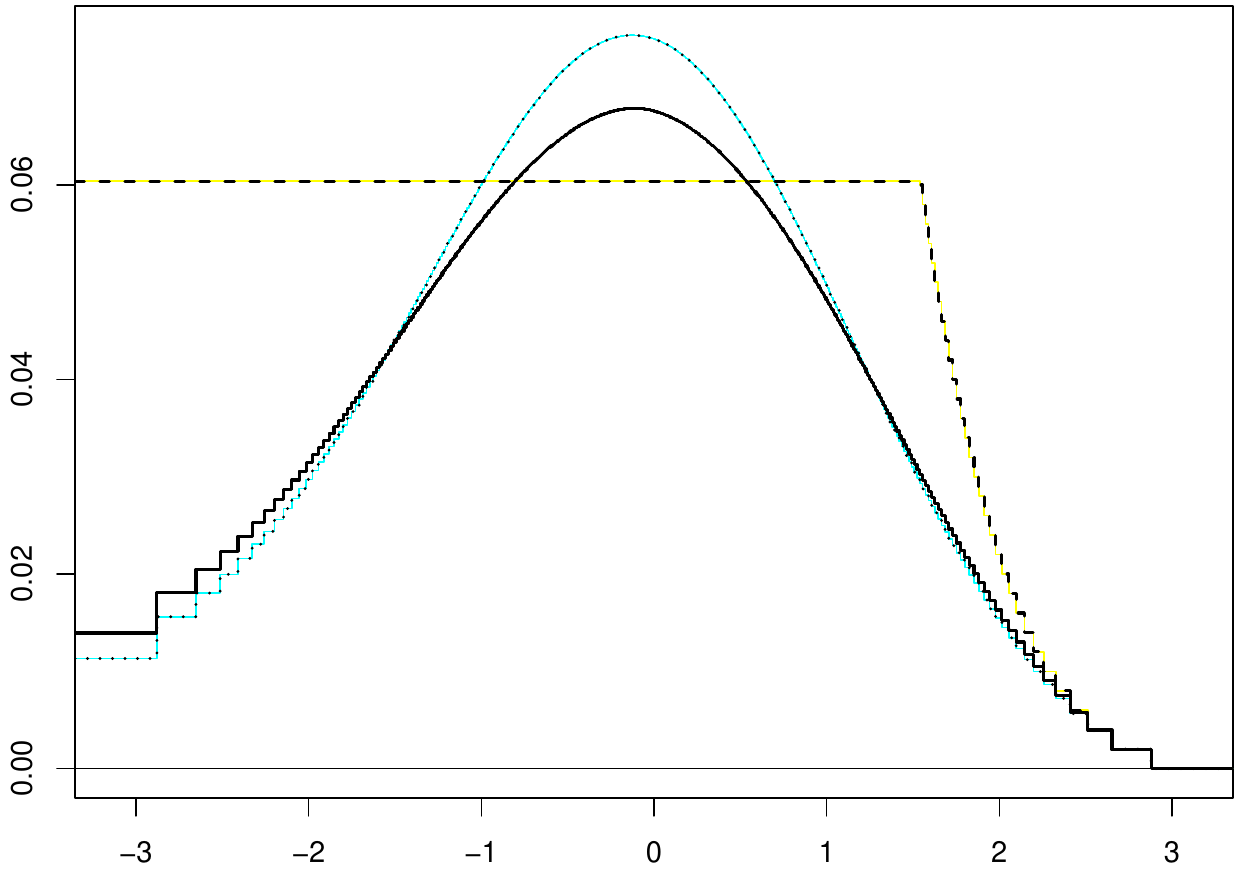}
\hfill
\includegraphics[width=0.49\textwidth]{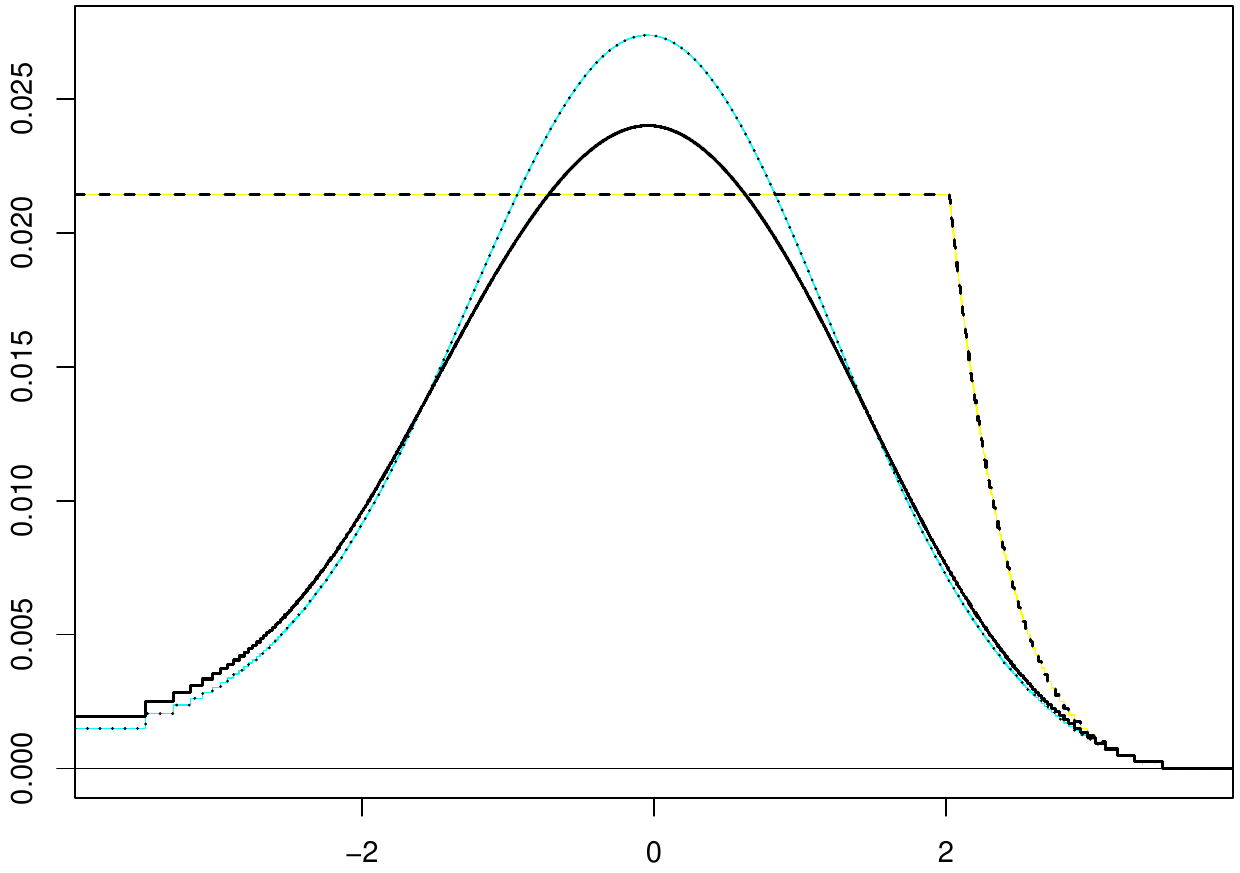}
\caption{Centered upper $95\%$-confidence bounds $B_{n,1,1,\alpha} - \FF_{n}$ (solid), $B_{n,1,\alpha}^{\rm BJO} - \FF_{n}$ (dotted) and $B_{n,\alpha}^{\rm KS} - \FF_{n}$ (dashed) for $n = 500$ (left panel) and $n = 4000$ (right panel).}
\label{fig:AsEfficiency.2}
\end{figure}

\begin{example}[$s = 1$]
The left panel in Figure~\ref{fig:AsEfficiency.1} depicts, for $n=100$, the $95\%$-confidence band $(A_{n,1,1,\alpha}, B_{n,1,1,\alpha})$ in case of an idealized standard Gaussian sample with order statistics
\[
	X_{n:i} = \Phi^{-1}(i/(n+1)) .
\]
In addition, one sees the Kolmogorov--Smirnov $95\%$-confidence band $(A_{n,\alpha}^{\rm KS}, B_{n,\alpha}^{\rm KS})$. In the right panel, one sees for the same setting the centered upper bounds $B_{n,1,1,\alpha} - \FF_n$, $B_{n,1,\alpha}^{\rm BJO} - \FF_n$ and $B_{n,\alpha}^{\rm KS} - \FF_n$. Note that a plot of the centered lower bounds $A_{n,1,1,\alpha} - \FF_n$, $A_{n,1,\alpha}^{\rm BJO} - \FF_n$ and $A_{n,\alpha}^{\rm KS} - \FF_n$ would be the reflection of the plots for the centered upper bounds with respect to the point $(0,0)$. The corresponding critical values $\kappa_{n,1,1,\alpha}$, $\kappa_{n,1,\alpha}^{\rm BJ}$ and $\kappa_{n,\alpha}^{\rm KS}$ have been computed numerically, see Section~S.7.

Figure~\ref{fig:AsEfficiency.2} shows the same as the right panel in Figure~\ref{fig:AsEfficiency.1}, but with sample sizes $n=500$ and $n=4000$ in the left and right panel, respectively.

In the online supplement, these bands $(A_{n,1,1,\alpha},B_{n,1,1,\alpha})$ are also compared with the confidence bands of \cite{Stepanova_Pavlenko_2018}, confirming the purely asymptotic result in
Remark~\ref{rem:Comparison.with.SP2018}.

\end{example}

\begin{figure}
\includegraphics[width=0.49\textwidth]{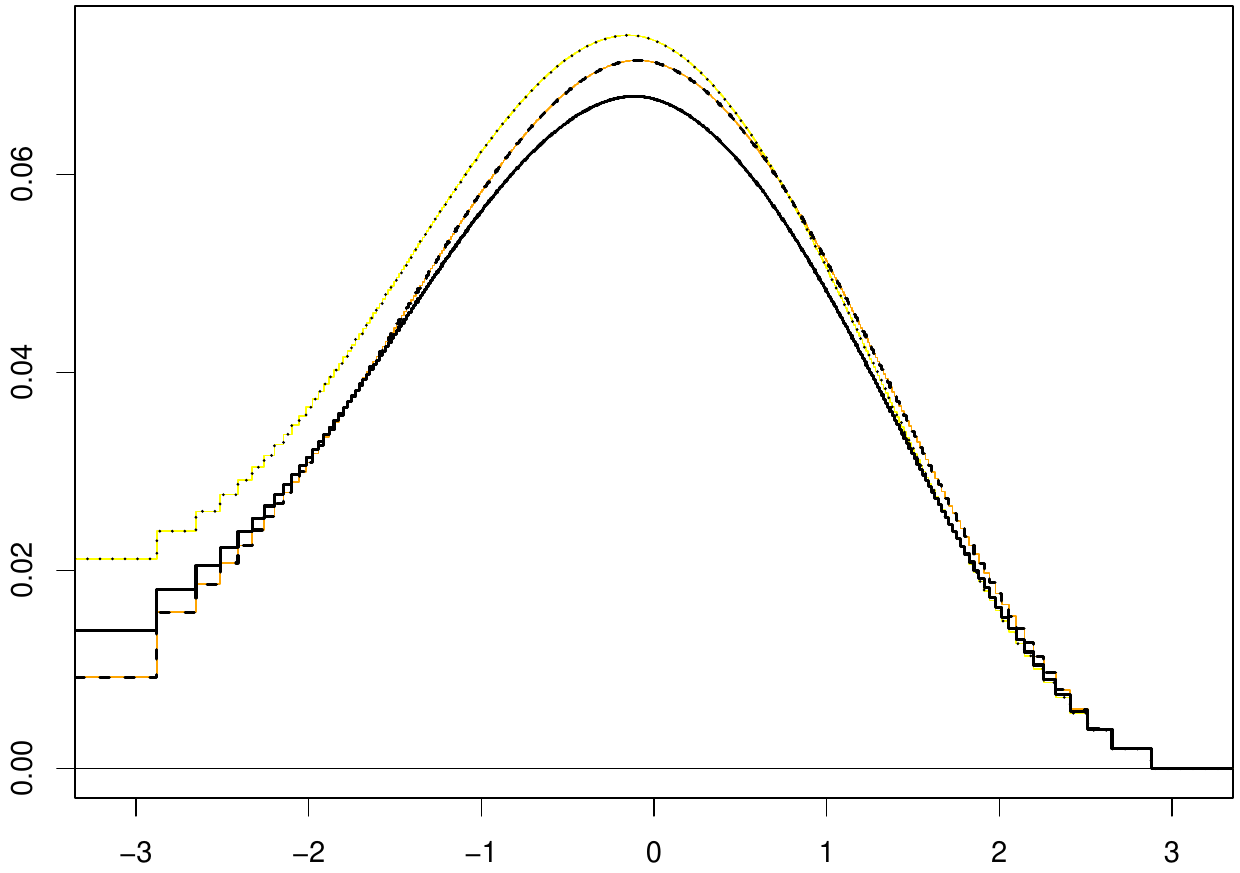}
\hfill
\includegraphics[width=0.49\textwidth]{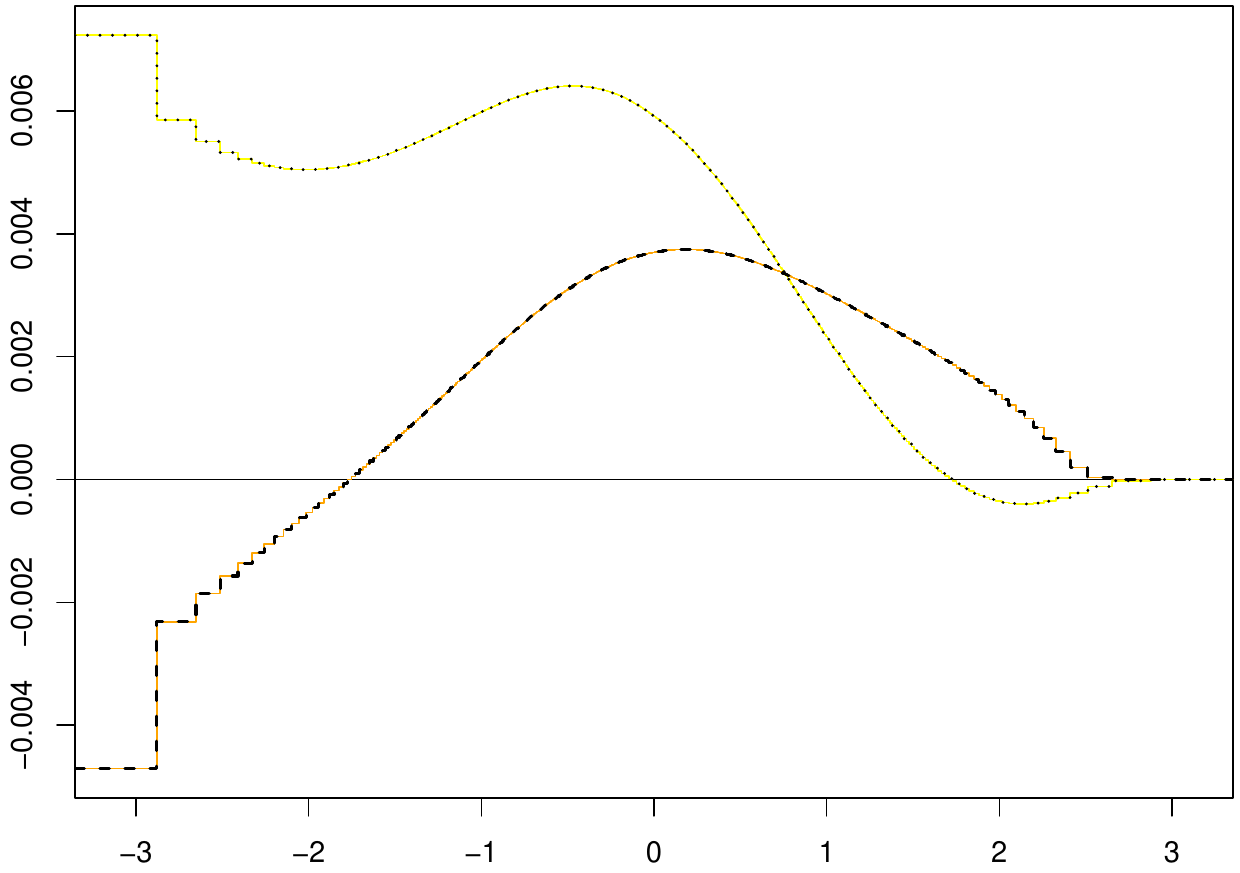}
\caption{Upper $95\%$-confidence bounds for $n = 500$. Left panel: centered bounds $B_{n,s,1,\alpha} - \FF_n$ for $s = 0.6$ (dashed), $s = 1.0$ (solid) and $s = 1.4$ (dotted). Right panel: differences $B_{n,s,1,\alpha} - B_{n,1,1,\alpha}$ for $s = 0.6$ (dashed) and $s = 1.4$ (dotted).}
\label{fig:UBS}
\end{figure}

\begin{example}[The impact of $s$]
Figure~\ref{fig:UBS} shows for an idealized Gaussian sample of size $n = 500$, the centered upper $95\%$-confidence bounds $B_{n,s,1,\alpha} - \FF_n$ for $s = 0.6, 1, 1.4$ (left panel) as well as the differences $B_{n,s,1,\alpha} - B_{n,1,1,\alpha}$ for $s = 0.6, 1.4$, right panel. As predicted by Theorem~\ref{thm:conf.bands.tails}, the upper bounds $B_{n,s,1}(x)$ are increasing in $s$ for small values of $x$ and decreasing in $s$ for large values of $x$. The online supplement contains further plots illustrating the impact of $s$ on our bands. These plots support our claim that choosing $s$ close to $1$ is preferable. Other values of $s$ increase the bands' precision somewhere in the tails, but lead to a substantial loss of precision in the central region.
\end{example}

\begin{remark}[Discontinuous distribution functions]
In the previous considerations, we focused on continuous distribution functions~$F$, and all confidence bands $(A_{n,\alpha},B_{n,\alpha})$ for $F$ we considered are of the form
\[
	\bigl[ A_{n,\alpha}(x), B_{n,\alpha}(x) \bigr]
	= [a_{n,\alpha,i}, b_{n,\alpha,i}]
	\quad\text{for} \ x \in [X_{n:i},X_{n:i+1}) \ \text{and} \ 0 \le i \le n
\]
with certain numbers $a_{n,\alpha,i}, b_{n,\alpha,i} \in [0,1]$. Interestingly, such a band has coverage probability at least $1 - \alpha$ for arbitrary, not necessarily continuous distribution functions $F$; see Section~S.6.
\end{remark}

\section{Proofs for Section~\ref{sec:LimitDistributionsUnif01}} 
\label{sec:ProofsSec2}

\subsection{Proof of Theorem~\ref{IntermediateThmPhiDiv}}

The following three facts are our essential ingredients.

\begin{fact}[\cite{Csorgo_etal_1986}, Theorem~2.2 and Corollary~2.1]
\label{fact1CCHM}
There exist on a common probability space a sequence of i.i.d.\ $U(0,1)$ random 
variables $\xi_1, \xi_2, \xi_3, \ldots$ and a sequence of Brownian bridge processes 
$\UU^{(1)}, \UU^{(2)}, \UU^{(3)}, \ldots $ such that, for all $0 \le \delta < 1/4$,
\[
	\sup_{t \in [1/n,1 - 1/n]}
		\frac{n^{\delta} \bigl| \UU_n(t) - \UU^{(n)}(t) \bigr|}{(t (1-t))^{1/2-\delta}}
	= O_p(1) .
\]
\end{fact}

\begin{fact}[\cite{Csorgo_etal_1986}, Theorem 4.4.1]
\label{fact2CCHM}
\[
	\sup_{t \in (0,1)} \frac{\UU_n(t)^2}{ 2 t (1-t) \log \log n}
	\rightarrow_p 1 .
\]
\end{fact}

\begin{fact}[\cite{Csorgo_etal_1986}, Lemma 4.4.4]
\label{fact3CCHM} 
For any $1 \le d_n \le n$ such that $d_n/n \rightarrow 0$ and $d_n \rightarrow \infty$,
\[
	\sup_{t \in (0,d_n/n]} \frac{\UU_n(t)^2}{ 2 t (1-t) \log \log d_n}
	\rightarrow_p 1 .
\]
The same holds with the supremum over $[1-d_n/n, 1)$.
\end{fact}

The asymptotic distribution of $\tilde{T}_{n,\nu}$ will be derived from the subsequent Lemmas~\ref{lemmaOne}, \ref{lemmaTwo} and \ref{lemmaThree}.

\begin{lemma}  
\label{lemmaOne}
For any sequence of constants $1 \le d_n \le n$ such that $d_n/n \rightarrow 0$ and 
$d_n \rightarrow \infty$ and any choice of $0 < \delta < 1/4$, 
\[
	\sup_{t \in [d_n/n, 1 - d_n/n]}
		\frac{\bigl| \UU_n(t)^2  - \UU^{(n)}(t)^2 \bigr|}{t(1-t)}
	= O_p \bigl( d_n^{-\delta} (\log\log n)^{1/2} \bigr) .
\]
\end{lemma}

\begin{proof}
By Fact~\ref{fact1CCHM}, for $0 < \delta < 1/4$, 
\[
	\sup_{t \in [d_n/n,1 - d_n/n]}
		\frac{|\UU_n(t)  - \UU^{(n)}(t)|}{(t(1-t))^{1/2}}
	\le O(d_n^{-\delta})
		\sup_{t \in [1/n,1 - 1/n]}
			\frac{n^{\delta} |\UU_n(t) - \UU^{(n)}(t)|}{(t(1-t))^{1/2-\delta}}
	= O_p (d_n^{-\delta}) .
\]
Together with Fact~\ref{fact2CCHM} and \eqref{eq:LIL.BB} this implies that
\begin{align*}
\lefteqn{\sup_{t \in [d_n/n,1 - d_n/n]}
		\frac{\bigl| \UU_n(t)^2 - \UU^{(n)}(t)^2 \bigr|}{t(1-t)} }\\
	&\le \sup_{t \in [d_n/n,1 - d_n/n]}
		\frac{|\UU_n(t)  - \UU^{(n)}(t)|}{(t(1-t))^{1/2}}
			\cdot \left( \frac{|\UU_n(t)|}{(t(1-t))^{1/2}} + \frac{|\UU^{(n)}(t)|}{(t(1-t))^{1/2}} \right) \\
	&= O_p \bigl( d_n^{-\delta} (\log\log n)^{1/2} \bigr) .
\end{align*}\\[-5ex]
\end{proof}

\begin{lemma}
\label{lemmaTwo}
For all $\nu \ge 0$,
\[
	\sup_{t \in (0,n^{-1} \log n]}
		\left( \frac{\UU_n(t)^2}{2 t(1-t)} - C_\nu(t) \right) \rightarrow_p - \infty .
\]
The same holds with the supremum over $(0,n^{-1} \log n]$ replaced by $[1 - n^{-1}\log n, 1)$.
\end{lemma}

\begin{proof}
Note that with $d_n = \log n$,
\begin{equation}
	\label{Lemma2Key1}
	\sup_{t \in (0,d_n/n]}
		\left( \frac{\UU_n(t)^2}{2 t(1-t)} - C_\nu(t) \right)
	\le \sup_{t \in (0,d_n/n]}
		\left( \frac{\UU_n(t)^2}{2 t(1-t)} - C(d_n/n) \right)
\end{equation}
since $C_\nu \ge C$ and $C$ is non-increasing. By Fact~\ref{fact3CCHM},
\[
	\sup_{t \in (0,d_n/n]}
		\frac{| \UU_n(t)^2|}{ 2 t(1-t) \log \log \log n}
	\rightarrow_p 1 ,
\]
while
\[
	\frac{C(d_n/n)}{\log \log \log n}
	= \frac{(1 + o(1)) \log\log n}{\log\log\log n}
	\rightarrow \infty .
\]
Thus, the right side of \eqref{Lemma2Key1} can be written as
\begin{align*}
	&\sup_{t \in (0,d_n/n]} \left(
		\frac{\UU_n(t)^2}{2 t(1-t)\log \log \log n}\cdot \log\log \log n - C(d_n/n)
		\right) \\
	&\qquad = \sup_{t \in (0,d_n/n]} \left(
		\frac{\UU_n(t)^2}{2 t(1-t)\log \log \log n} - \frac{C(d_n/n)}{\log \log \log n} \right)
			\log \log \log n \\
	&\qquad \rightarrow_p (1-\infty)\cdot \infty = - \infty .
\end{align*}\\[-5ex]
\end{proof}

\begin{lemma}
\label{lemmaThree}
For any fixed $\nu > 3/4$,
\[
	\sup_{t \in (0,\rho] \cup [1 - \rho,1)}
		\Bigl( \frac{\UU(t)^2}{2 t(1-t)} - C_\nu(t) \Bigr)
	\to - \infty
	\quad\text{almost surely as} \ \rho \searrow 0 .
\]
\end{lemma}

\begin{proof}
Recall that
\[
	T_{\nu} = \sup_{t \in (0,1)} \Bigl( \frac{\UU(t)^2}{2 t(1-t)} - C(t) - \nu D(t) \Bigr)
\]
is finite almost surely for any $\nu > 3/4$. If we choose $\nu' \in (3/4, \nu)$ and write $\nu D(t) = \nu' D(t) + (\nu - \nu')D(t)$, then we see that for any $\rho \in (0,1/2]$,
\begin{align*}
	\sup_{t \in (0,\rho] \cup [1-\rho,1)}
		\left( \frac{\UU(t)^2}{2 t(1-t)} - C(t) - \nu D(t) \right)
	& \le \sup_{t \in (0,\rho] \cup [1-\rho,1)}
		\bigl( T_{\nu'} - (\nu - \nu') D(t) \bigr) \\
	& = T_{\nu'} - (\nu - \nu') D(\rho) ,
\end{align*}
because $D(\cdot)$ is symmetric around $1/2$ and monotone decreasing on $(0,1/2]$. Now the claim follows from $T_{\nu'} < \infty$ almost surely and $D(\rho) \to \infty$ as $\rho \searrow 0$.
\end{proof}

Now we can finish the proof of Theorem~\ref{IntermediateThmPhiDiv}. According to Lemmas~\ref{lemmaTwo} and \ref{lemmaThree}, with $d_n := \log n$,
\[
	\left. \begin{array}{c}
		\tilde{T}_{n,\nu} \\ T_\nu
	\end{array} \!\!\right\}
	= \sup_{t \in [d_n/n, 1 - d_n/n]}
		\left( \frac{1}{2t(1-t)} \left\{ \!\! \begin{array}{c}
			\UU_n(t)^2 \\ \UU(t)^2
		\end{array} \!\! \right\} - C_\nu(t) \right) 
\]
with asymptotic probability one. If we replace the Brownian bridge $\UU$ with the Brownian bridge $\UU^{(n)}$, then Lemma~\ref{lemmaOne} implies that the latter two suprema over $[d_n/n, 1 - d_n/n]$ differ only by $o_p(1)$. Consequently, $\tilde{T}_{n,\nu}$ converges in distribution to $T_\nu$.

\subsection{Proof of Theorem~\ref{thm:PhiDivergenceNull}}

Note first that in case of $s > 0$,
\[
	\sup_{t \in (0,\xi_{n:1})} \bigl( n K_s(\GG_n(t),t) - C_\nu(\GG_n(t),t) \bigr)
	= n K_s(0,\xi_{n:1}) - C_\nu \bigl( \min(\xi_{n:1},1/2) \bigr)
	\to_p - \infty ,
\]
because $K_s(0,t) = t/s + o(t)$ as $t \searrow 0$ and $E(\xi_{n:1}) = 1/(n+1)$. Since $K_s(1,t) = K_s(0,1-t)$, $C_\nu(t) = C_\nu(1-t)$ and $\xi_{n:1} \stackrel{d}{=} 1 - \xi_{n:n}$,
\[
	\sup_{t \in [\xi_{n:n},1)} \bigl( n K_s(\GG_n(t),t) - C_\nu(\GG_n(t),t) \bigr)
	= n K_s(1,\xi_{n:n}) - C_\nu \bigl( \max(\xi_{n:n},1/2) \bigr)
	\to_p - \infty .
\]
Consequently, it suffices to verify Theorem~\ref{thm:PhiDivergenceNull} with the modified test statistic
\[
	T_{n,s,\nu} := \sup_{t \in [\xi_{n:1},\xi_{n:n})}
		\bigl( n K_s(\GG_n(t),t) - C_\nu(\GG_n(t),t) \bigr) ,
\]
provided that we can show that the latter converges in distribution.

In what follows, we show that replacing $s$ with $2$ and $C_\nu(\GG_n(t),t)$ with $C_\nu(t)$ has no effect asymptotically. For these tasks, the following two facts are useful.

\begin{fact}[Linear bounds for $\GG_n$]\strut\\
\label{Fact4}
A. \ By inequality 1, 
\cite{Shorack_Wellner_1986, Shorack_Wellner_2009}, page 415, 
\[
	\sup_{\xi_{n:1} \le t \le 1} \frac{t}{\GG_n(t)} = O_p (1)
	\quad\text{and}\quad
	\sup_{0 \le t < \xi_{n:n}} \frac{1-t}{1-\GG_n(t)} = O_p (1) .
\]
B. \ From Daniels' theorem (Theorem 2, \cite{Shorack_Wellner_1986, Shorack_Wellner_2009}, page 341),
\[
	\sup_{0 < t \le 1} \frac{\GG_n(t)}{t} = O_p (1)
	\quad\text{and}\quad
	\sup_{0 \le t < 1} \frac{1-\GG_n(t)}{1-t} = O_p (1) .
\]
\end{fact}

\begin{fact}
\label{Fact5}
For any sequence of constants $d_n$ with $1 \le d_n \le n$ such that $d_n /n \rightarrow 0$ and $d_n \rightarrow \infty$
\[
	\sup_{d_n/n \le t \le 1} \frac{|\GG_n(t) - t|}{t} = O_p (d_n^{-1/2})
\]
and
\[
	\sup_{0 \le t \le 1- d_n/n} \frac{|\GG_n(t) - t|}{1-t} = O_p (d_n^{-1/2})
\]
(\cite{Wellner_1978}, Lemma~3 and Theorem~1S; \cite{Shorack_Wellner_1986, Shorack_Wellner_2009}, Chapter 10, Section 5, page 424). In fact,
\[
	d_n^{1/2} \sup_{d_n/n \le t \le 1} \frac{|\GG_n(t) - t|}{t}
	\rightarrow_d \sup_{0 \le t \le 1} |\WW(t)| ,
\]
where $\WW$ is a standard Brownian motion, see \cite{Renyi_1969}.
\end{fact}

A particular consequence of Fact~\ref{Fact4} is that
\begin{equation}
\label{eq:Mn1}
	M_{n,1} := \sup_{t \in [\xi_{n:1},\xi_{n:n})} \bigl| \logit(\GG_n(t)) - \logit(t) \bigr|
	= O_p(1) ,
\end{equation}
where $\logit(t) := \log(t/(1-t))$, and Fact~\ref{Fact5} implies that
\begin{equation}
\label{eq:Mn2}
	M_{n,2} := \sup_{t \in [n^{-1} \log n, 1 - n^{-1} \log n]}
		\bigl| \logit(\GG_n(t)) - \logit(t) \bigr|
	= O_p \bigl( (\log n)^{-1/2} \bigr) ,
\end{equation}
with the conventions that $\logit(0) := -\infty$ and $\logit(1) := \infty$. This leads to the following useful bounds:

\begin{lemma}
\label{Lemma4}
For any fixed $s \in \RR$,
\[
	\sup_{t \in [\xi_{n:1},\xi_{n:n})} \frac{K_s(\GG_n(t),t)}{K_2(\GG_n(t),t)}
		= O_p(1)
	\quad\text{and}\quad
	\sup_{t \in [\xi_{n:1},\xi_{n:n})} \bigl( C_\nu(t) - C_\nu(\GG_n(t),t) \bigr)
		= O_p(1) ,
\]
where $K_s(t,t)/K_2(t,t) := 1$. Moreover,
\begin{align*}
	\sup_{t \in [n^{-1} \log n, 1 - n^{-1} \log n]}
			\Bigl| \frac{K_s(\GG_n(t),t)}{K_2(\GG_n(t),t)} - 1 \Bigr|
		&= O_p \bigl( (\log n)^{-1/2} \bigr)
			\quad\text{and} \\
	\sup_{t \in [n^{-1} \log n, 1 - n^{-1} \log n]}
			\bigl( C_\nu(t) - C_\nu(\GG_n(t),t) \bigr)
		&= O_p \bigl( (\log n)^{-1/2} \bigr) ,	
\end{align*}
where $K_s(0,t) = K_s(1,t) := \infty$ in case of $s < 1$.
\end{lemma}

\begin{proof}
With the auxiliary quantities $M_{n,1}$ in \eqref{eq:Mn1} and $M_{n,2}$ in \eqref{eq:Mn2}, it follows from the inequalities (S.14) and Lemma~S.10 that for $\xi_{n:1} \le t < \xi_{n:n}$,
\begin{align*}
	\frac{K_s(\GG_n(t),t)}{K_2(\GG_n(t),t)}
	& \le \exp \bigl( |s-2| M_{n,1} \bigr)
	= O_p(1)
	\quad\text{and}\quad \\
	0 & \le C_\nu(t) - C_\nu(\GG_n(t),t)
	\le (1 + \nu) M_{n,1}
	= O_p(1) .
\end{align*}
Moreover, for $n^{-1} \log n \le t \le 1 - n^{-1} \log n$,
\begin{align*}
	\Bigl| \frac{K_s(\GG_n(t),t)}{K_2(\GG_n(t),t)} - 1 \Bigr|
		\le \exp \bigl( |s - 2| M_{n,2} \bigr) - 1
		&= O_p \bigl( (\log n)^{-1/2}) \bigr)
		\quad\text{and} \\
	0 \le C_\nu(t) - C_\nu(\GG_n(t),t)
		\le (1 + \nu) M_{n,2}
		&= O_p \bigl( (\log n)^{-1/2}) \bigr) .
\end{align*}
(Note that $M_{n,2} = \infty$ if $t < \xi_{n:1}$ or $t \ge \xi_{n:n}$.)
\end{proof}

Now the statement about the (modified) test statistic $T_{n,s,\nu}$ is an immediate consequence of Theorem~\ref{IntermediateThmPhiDiv} and the following lemma.

\begin{lemma}
\label{Lemma5}
For $\nu > 3/4$ and any $s \in \RR$,
\[
	T_{n,s,\nu} = \tilde{T}_{n,\nu} + o_p(1) .
\]
\end{lemma}

\begin{proof}
With $d_n := \log n$, we know that $\xi_{n:n} > 1 - d_n/n$ with asymptotic probability one, and thus it follows from Fact~\ref{fact3CCHM} and Lemma~\ref{Lemma4} that
\begin{align*}
	& \lefteqn{\sup_{t \in [\xi_{n:1},d_n/n]} n K_s(\GG_n(t),t) } \\
	& \le \ \sup_{t \in [\xi_{n:1},1 - d_n/n]} \frac{K_s(\GG_n(t),t)}{K_2(\GG_n(t),t)}
		\sup_{t \in (0,d_n/n]} n K_2(\GG_n(t),t)
	= O_p(\log\log\log n) .
\end{align*}
On the other hand,
\[
	\min_{t \in [\xi_{n:1},d_n/n]} C_\nu(\GG_n(t),t) \ge C(d_n/n) + O_p(1) = (1 + o(1)) \log\log n .
\]
Hence,
\[
	\sup_{t \in [\xi_{n:1},d_n/n]} \bigl( n K_s(\GG_n(t),t) - C_\nu(\GG_n(t),t) \bigr)
	\to_p - \infty ,
\]
and for symmetry reasons,
\[
	\sup_{t \in [1 - d_n/n,\xi_{n:n}]} \bigl( n K_s(\GG_n(t),t) - C_\nu(\GG_n(t),t) \bigr)
	\to_p - \infty .
\]
Since $\tilde{T}_{n,\nu}$ is equal to
\[
	\tilde{T}_{n,\nu}^{\rm restr} = \sup_{t \in [d_n/n, 1 - d_n/n]} \bigl( n K_2(\GG_n(t),t) - C_\nu(t) \bigr)
\]
with asymptotic probability one, it suffices to show that
\[
	T_{n,s,\nu}^{\rm restr} := \sup_{t \in [d_n/n, 1 - d_n/n]} \bigl( n K_s(\GG_n(t),t) - C_\nu(\GG_n(t),t) \bigr)
	= \tilde{T}_{n,\nu}^{\rm restr} + o_p(1) .
\]
To this end, note that $\tilde{T}_{n,\nu}^{\rm restr} \to_d T_\nu$ implies that
\[
	\sup_{t \in [d_n/n, 1 - d_n/n]} n K_2(\GG_n(t),t)
	\le C_\nu(d_n/n) + O_p(1)
	= (1 + o_p(1)) \log \log n .
\]
Consequently,
\begin{align*}
        \lefteqn{
	\bigl| T_{n,s,\nu}^{\rm restr} - \tilde{T}_{n,\nu}^{\rm restr} \bigr| } \\
	&\le \sup_{t \in [d_n/n, 1 - d_n/n]} \bigl| n K_s(\GG_n(t),t) - n K_2(\GG_n(t),t) \bigr|
		+ O_p \bigl( (\log n)^{-1/2}) \bigr) \\
	&\le \sup_{t \in [d_n/n, 1 - d_n/n]} \Bigl| \frac{K_s(\GG_n(t),t)}{K_2(\GG_n(t),t)} - 1 \Bigr|
		\sup_{t \in [d_n/n, 1 - d_n/n]} n K_2(\GG_n(t),t)
		+ O_p \bigl( (\log n)^{-1/2}) \bigr) \\
	&= O_p \bigl( (\log n)^{-1/2} \bigr) (1 + o_p(1)) \log\log n
		= o_p(1) .
\end{align*}\\[-5ex]
\end{proof}

It remains to prove the claim that $\kappa_{n,s,\nu,\alpha} \to \kappa_{\nu,\alpha} > 0$. But this follows immediately from the following lemma.

\begin{lemma}
\label{lem:continuity.limit}
Let $G(r) := P(T_{\nu} \le r)$. Then $G(0) = 0$, and $G$ is continuous and strictly increasing on $[0,\infty)$.
\end{lemma}

To prove this lemma and other results, we make use of the following well-known result.

\begin{fact}[\cite{Borell_1974}, Corollary 2.1; \cite{Gaenssler_Rost_2007}, Lemma~1.1]
\label{Fact6}
The distribution $Q$ of $\UU$ is a log-concave measure on $\mathcal{C}[0,1]$. That means, for Borel sets $\mathcal{B}_0,\mathcal{B}_1 \subset \mathcal{C}[0,1]$ and $\lambda \in (0,1)$,
\[
	\log Q_*((1 - \lambda)\mathcal{B}_0 + \lambda \mathcal{B}_1)
	\ge (1 - \lambda) Q(\mathcal{B}_0) + \lambda Q(\mathcal{B}_1) ,
\]
where $Q_*$ stands for the inner measure induced by $Q$, and $(1 - \lambda) \mathcal{B}_0 + \lambda \mathcal{B}_1 := \{(1 - \lambda) g_0 + \lambda g_1 : g_0 \in \mathcal{B}_0, g_1 \in \mathcal{B}_1\}$.
\end{fact}

From this fact one can deduce the following properties of $\UU$:

\begin{proposition}
\label{prop:log-concavity.Q}
For arbitrary functions $h : [0,1] \to [0,\infty)$ and $h_o : [0,1] \to \RR$,
\[
	G_1(x) := P(|x h_o + \UU| \le h)
\]
is an even, log-concave function of $x \in \RR$. Furthermore, if $h_o \ge 0$, then
\[
	G_2(x) := P \bigl( |\UU| \le \sqrt{h + x h_o} \bigr)
\]
is a non-decreasing and log-concave function of $x \ge 0$.
\end{proposition}

Let $\WW$ be a standard Brownian motion process on $[0,1]$. Then it is well-known
that $\UU(t) := \WW(t) - t \WW(1)$ defines a Brownian bridge process on $[0,1]$.  
The  following self-similarity property of the Brownian bridge process 
$\UU$ seems to be less well-known.

\begin{proposition}
\label{prop:self-similarity.BB}
For fixed numbers $0 \le a < b \le 1$, define a stochastic process $\ZZ_{a,b}$ on $[0,1]$ as follows:
\[
	\ZZ_{a,b}(v) := \UU((1 - v) a + v b) - (1 - v) \UU(a) - v \UU(b) ,
\]
that is, $\ZZ_{a,b}$ describes the interpolation error when replacing $\UU$ on $[a,b]$ with its linear interpolation there. Then the two processes $(\UU(t))_{t \in [0,1] \setminus (a,b)}$ and $\ZZ_{a,b}$ are stochastically independent, and
\[
	\ZZ_{a,b} \stackrel{d}{=} \sqrt{b - a} \, \UU .
\]
\end{proposition}

Proofs of Propositions~\ref{prop:log-concavity.Q} and \ref{prop:self-similarity.BB} are provided in Section~S.4.

\begin{proof}[Proof of Lemma~\ref{lem:continuity.limit}]
Note first that the distribution function 
$r \mapsto G(r) $ coincides with the function $G_2$ in 
Proposition~\ref{prop:log-concavity.Q}, where $h(t) := 2t(1 - t) C_\nu(t)$ and $h_o(t) := 2 t (1 - t)$. In particular, 
$G(r) \le P \bigl( |\UU(1/2)| \le \sqrt{r/2} \bigr)$, and the 
latter bound equals $0$ for $r = 0$ and is strictly smaller than $1$ for any $r \ge 0$.  

By Proposition~\ref{prop:log-concavity.Q}, $G : [0,\infty) \to [0,1]$ is log-concave, 
and since $G(r) < 1 = \lim_{s \to \infty} G(s)$ for all $r \ge 0$, this implies that 
$G$ is continuous and strictly increasing on $(r_o,\infty)$, where 
$r_o := \inf \{r > 0 : G(r) > 0\}$. If we can show that $r_o = 0$, then we know 
that $G$ is, in fact, continuous and strictly increasing on $[0,\infty)$.

To show that $G(r) > 0$ for any $r > 0$, we pick a number 
$\rho \in (0,1/2)$ and write $T_\nu$ as the maximum of the three random variables
\begin{align*}
	T_\nu^{(\rho,1)}
		&:= \max_{t \in [\rho,1-\rho]} \bigl( \UU(t)^2/[2t(1-t)] - C_\nu(t) \bigr) , \\
	T_\nu^{(\rho,2,L)}
		&:= \max_{t \in (0,\rho]} \bigl( \UU(t)^2/[2t(1-t)] - C_\nu(t) \bigr) , \\
	T_\nu^{(\rho,2,R)}
		&:= \max_{t \in [1-\rho,1)} \bigl( \UU(t)^2/[2t(1-t)] - C_\nu(t) \bigr) .
\end{align*}
Then we can write
\begin{align*}
	G(r)
	& = P \bigl( T_\nu^{(\rho,1)} \le r,
		T_\nu^{(\rho,2,L)} \le r, T_\nu^{(\rho,2,R)} \le r \bigr) \\
	& \ge P \Bigl( \max_{t \in [\rho,1-\rho]} |\UU(t)| \le \delta,
		T_\nu^{(\rho,2,L)} \le 0, T_\nu^{(\rho,2,R)} \le 0 \Bigr)
\end{align*}
with $\delta := \sqrt{2\rho(1 - \rho) r} > 0$. 

According to Lemma~\ref{lemmaThree}, we may choose $\rho$ such that $P(T_\nu^{(\rho,2,L)} \le 0) = P(T_\nu^{(\rho,2,R)} \le 0) \ge 1/2$. Now we apply Proposition~\ref{prop:self-similarity.BB} twice, first with $[a,b] = [0,\rho]$, and then with $[a,b] = [1-\rho,1]$. This shows that $\UU$ may be rewritten on $[0,\rho]$ and on $[1 - \rho,1]$ as follows: for $v \in [0,1]$,
\begin{align*}
	\UU(\rho v)
		&= \ v \UU(\rho) + \sqrt{\rho} \, \UU^{(L)}(v) , \\
	\UU(1 - \rho v)
		&= v \UU(1 - \rho) + \sqrt{\rho} \, \UU^{(R)}(v) ,
\end{align*}
where $\UU, \UU^{(L)}, \UU^{(R)}$ are independent Brownian bridge processes. In particular,
\begin{align*}
	P \bigl( & T_\nu^{(\rho,2,L)} \le 0
		\,\big|\, (\UU(t))_{t \in [\rho,1-\rho]} \bigr) \\
	&= P \bigl( \bigl| v \UU(\rho) + \sqrt{\rho} \, \UU^{(L)}(v) \bigr|
		\le \sqrt{2 \rho v (1 - \rho v) C_\nu(\rho v)} \ \text{for all} \ v \in [0,1]
		\,\big|\, (\UU(t))_{t \in [\rho,1-\rho]} \bigr) \\
	&= P \bigl( \bigl| \UU(\rho) v/\sqrt{\rho} + \UU^{(L)}(v) \bigr|
		\le \sqrt{2 v (1 - \rho v) C_\nu(\rho v)} \ \text{for all} \ v \in [0,1]
		\,\big|\, (\UU(t))_{t \in [\rho,1-\rho]} \bigr) \\
	&= G_1(\UU(\rho)) ,
\end{align*}
where $G_1(x) := P( |x h_o + \UU| \le h)$ with $h_o(v) := v / \sqrt{\rho}$ and $h(v) := \sqrt{2v(1 - \rho v) C_\nu(\rho v)}$ for $v \in [0,1]$. Analogously,
\[
	P \bigl( T_\nu^{(\rho,2,R)} \le 0
		\,\big|\, (\UU(t))_{t \in [\rho,1-\rho]} \bigr)
	= G_1(\UU(1 - \rho)) .
\]
According to Proposition~\ref{prop:log-concavity.Q}, $G_1$ is an even, log-concave function on $\RR$. Since $1/2 \le P(T_\nu^{(\rho,2,L)} \le 0) = E[G_1(\UU(\rho))]$, there exists a $\delta_o > 0$ such that $G_1(x) \ge 1/2$ for all $x \in [-\delta_o,\delta_o]$. Consequently,
\[
	G(r)
	\ge E \bigl( 1_{[|\UU| \le \delta \ \text{on} \ [\rho,1-\rho]]} G_1(\UU(\rho)) G_1(\UU(1 - \rho)) \bigr)
	\ge 4^{-1} P \bigl( \|\UU\|_\infty \le \min(\delta,\delta_o) \bigr)
	> 0 .
\]
That $P \bigl( \|\UU\|_\infty \le \lambda) > 0$ for any $\lambda > 0$ follows, for instance, from the expansion
\[
	P \bigl( \|\UU\|_\infty \le \lambda \bigr) 
	= \frac{\sqrt{2 \pi}}{8 \lambda^2} \exp \biggl( - \frac{\pi^2}{8 \lambda^2} \biggr)
		(1 + o(1))
	\quad\text{as} \ \lambda \searrow 0 ;
\]
see \cite{Mogulski_1979} or \cite{Shorack_Wellner_2009}, pp. 526-527. Alternatively, one could use Proposition~\ref{prop:log-concavity.Q} and separability of $\mathcal{C}[0,1]$.
\end{proof}

\section{Proofs for Section~\ref{sec:Implications}}
\label{sec:ProofsSec3}

\subsection{Proofs for Subsection~\ref{subsec:GoF.tests}}

\begin{proof}[Proof of Theorem~\ref{thm:GOF-PwrOneCondition}]
Let $(x_n)_n$ be a sequence in $\RR$ such that $\Delta_n(x_n) \to \infty$. Then for any fixed $\kappa > 0$,
\begin{align}
\label{ineq:Power.1}
	P_{F_n} \bigl[ T_{n,s,\nu}(F_0) &\le \kappa \bigr]
	\le P_{F_n} \bigl[ x_n \not\in [X_{n:1},X_{n:n}) \bigr] \\
\nonumber
	&+ \ P_{F_n} \bigl[ n K_s(\FF_n(x_n),F_0(x_n))
			\le C_\nu(\FF_n(x_n),F_0(x_n)) + \kappa \bigr] ,
\end{align}
where $K_s(u,\cdot) := \infty$ if $s \le 0$ and $u \in \{0,1\}$.

To ensure that the first summand on the right hand side of \eqref{ineq:Power.1} converges to $0$, we show that $x_n$ may be chosen such that $d_n/n \le F_n(x_n) \le 1 - d_n/n$, where $d_n := \log\log n$. To this end we have to analyze the auxiliary function $H_n$ in more detail. Elementary calculus reveals that for $t \in [0,1]$, $(1 + C(t)) t(1 - t)$ is an increasing and $1 + C(t)$ is a decreasing function of $t(1 - t) \in [0,1/4]$. Moreover,
\[
	1 + C(d_n/n) = (1 + o(1)) d_n
	\quad\text{and}\quad
	(d_n/n)(1 - d_n/n) = (1 + o(1)) d_n/n ,
\]
whence
\[
	\min_{t \in [0,1]} H_n(t) \ge (1 + o(1)) n^{-1/2} d_n
	\quad\text{and}\quad
	H_n(d_n/n) \ = \ (2 + o(1)) n^{-1/2} d_n .
\]
In particular,
\[
	|F_n - F_0|(x_n) \ge \Delta_n(x_n) (1 + o(1)) d_n/n .
\]
Now suppose that $F_n(x_n) < d_n/n$. With $\tilde{x}_n := F_n^{-1}(d_n/n)$ we may conclude that
\[
	F_n(\tilde{x}_n) \ge F_n(x_n) > |F_n - F_0|(x_n) - d_n/n \ge \Delta_n(x_n) (1 + o(1)) d_n/n .
\]
In particular, $\max\{ d_n/n, F_n(x_n)\} $ is of order $o(F_n(\tilde{x}_n))$, so
\begin{align*}
	\Delta_n(\tilde{x}_n)
	\ge \frac{\sqrt{n} |F_n - F_0|}{H_n(F_n)}(\tilde{x}_n)
	\ge \frac{(1 + o(1)) \sqrt{n} F_n(\tilde{x}_n)}{(2 + o(1)) n^{-1/2} d_n}
	\ge \ (1/2 + o(1)) \Delta_n(x_n)
	\to \infty .
\end{align*}
Analogously one can show that in case of $F_n(x_n) > 1 - d_n/n$, we may replace 
$x_n$ with $\tilde{x}_n := F_n^{-1}(1 - d_n/n)$ at the cost of reducing $\Delta_n(x_n)$ by a factor of at most $1/2 + o(1)$.

It remains to show that
\begin{equation}
\label{ineq:Power.2}
	P_{F_n} \bigl[ n K_s(\FF_n(x_n),F_0(x_n)) \le C_\nu(\FF_n(x_n),F_0(x_n)) + \kappa \bigr]
	\to 0 .
\end{equation}
By means of the second part of Lemma~S.12, the inequality for $K_s(\FF_n(x_n),F_0(x_n))$ implies that
\begin{align*}
	\sqrt{n} |\FF_n - F_0|(x_n)
	&\le \sqrt{2 (C_\nu(\FF_n,F_0) + \kappa)
		\min \bigl\{ \FF_n(1 - \FF_n), F_0(1 - F_0) \bigr\}}(x_n) \\
		& \qquad + 2(C_\nu(\FF_n,F_0) + \kappa)(x_n)/\sqrt{n} \\
	&\le 2 \max(1 + \nu, \kappa) \min \bigl\{ H_n(\FF_n), H_n(F_0) \bigr\}(x_n) ,
\end{align*}
because $C_\nu(\FF_n,F_0) \le \min \bigl\{ C_\nu(\FF_n), C_\nu(F_0) \bigr\}$, and for the univariate function $C_\nu$, it follows from $D \le C$ that $C_\nu + \kappa \le \max(1+\nu, \kappa) (1 + C)$. Moreover, the assumption that $d_n/n \le F_n(x_n) \le 1 - d_n/n$ implies that
\[
	\frac{h(\FF_n)}{h(F_n)}(x_n) \to_p 1
	\quad\text{for} \ h(t) = t, 1 + C(t), t(1 - t) .
\]
Consequently, \eqref{ineq:Power.2} would be a consequence of
\begin{equation}
\label{ineq:Power.3}
	P_{F_n} \bigl[ \sqrt{n} |\FF_n - F_0|(x_n) \le O_p(1)
		\min \bigl\{ H_n(F_n), H_n(F_0) \bigr\}(x_n) \bigr]
	\to 0 .
\end{equation}

To bound the left-hand side of \eqref{ineq:Power.3} we consider the quantity
\[
	M_n := \max \Bigl\{ \frac{F_0(1 - F_0)}{F_n(1 - F_n)}(x_n), \frac{F_n(1 - F_n)}{F_0(1 - F_0)}(x_n) \Bigr\}
	\ge 1
\]
and distinguish two cases. Suppose first that $M_n \le \Delta_n(x_n)$. Since
\begin{align*}
	& \frac{1 + C(F_n)}{1 + C(F_0)}(x_n) \le 1 \le \frac{F_n(1 - F_n)}{F_0(1 - F_0)}(x_n) \le M_n 
	    \quad\text{or}\quad \\
	& \frac{F_n(1 - F_n)}{F_0(1 - F_0)}(x_n) \le 1 \le \frac{1 + C(F_n)}{1 + C(F_0)}(x_n) \le 1 + \log M_n ,
\end{align*}
the definition of $H_n$ implies that
\[
	\frac{H_n(F_n)}{H_n(F_0)}(x_n) \ \le \Delta_n(x_n)^{1/2} .
\]
Then it follows from $\sqrt{n} (\FF_n - F_n)(x_n) = O_p \bigl( \sqrt{F_n(1 - F_n)}(x_n) \bigr) = O_p \bigl( H_n(F_n(x_n)) \bigr)$ that
\begin{align*}
	P_{F_n} \bigl[ & \sqrt{n} |\FF_n - F_0|(x_n) \le O_p(1)
		\min \bigl\{ H_n(F_n), H_n(F_0) \bigr\}(x_n) \bigr] \\
	&\le P_{F_n} \bigl[ \sqrt{n} |F_n - F_0|(x_n)
		\le O_p(1) \min \bigl\{ H_n(F_n), H_n(F_0) \bigr\}(x_n)
			+ O_p \bigl( H_n(F_n(x_n)) \bigr) \bigr] \\
	&\le P_{F_n} \bigl[ \sqrt{n} |F_n - F_0|(x_n)
		\le O_p \bigl( \Delta_n(x_n)^{1/2} \bigr)
			\min \bigl\{ H_n(F_n), H_n(F_0) \bigr\}(x_n) \bigr] \\
	&= P_{F_n} \bigl[ \Delta_n(x_n) \le O_p \bigl( \Delta_n(x_n)^{1/2} \bigr) \bigr]
		 \to 0 .
\end{align*}
Now suppose that $M_n \ge \Delta_n(x_n)^{1/2}$. Then,
\[
	\frac{|\FF_n - F_0|}{|F_n - F_0|}(x_n)
	\ge 1 - \frac{|\FF_n - F_n|}{|F_n - F_0|}(x_n)
	\ge 1 - \frac{|\FF_n - F_n|}{\bigl| F_n(1 - F_n) - F_0(1 - F_0) \bigr|}(x_n)
	= 1 + O_p(\rho_n)
\]
with
\begin{align*}
	\rho_n &:= \frac{\sqrt{F_n(1 - F_n)}}{\sqrt{n} \bigl| F_n(1 - F_n) - F_0(1 - F_0) \bigr|}(x_n) \\
	&= \frac{F_n(1 - F_n)}{\sqrt{n F_n(1 - F_n)} \bigl| F_n(1 - F_n) - F_0(1 - F_0) \bigr|}(x_n) 
	\le \frac{M_n}{(1 + o(1)) \sqrt{d_n} (M_n - 1)}
		\to 0 .
\end{align*}
Consequently,
\begin{align*}
	P_{F_n} \bigl[ & \sqrt{n} |\FF_n - F_0|(x_n) \le O_p(1)
		\min \bigl\{ H_n(F_n), H_n(F_0) \bigr\}(x_n) \bigr] \\
	&\le P_{F_n} \bigl[ \sqrt{n} |F_n - F_0|(x_n) (1 + o_p(1))
		\le O_p(1) \min \bigl\{ H_n(F_n), H_n(F_0) \bigr\}(x_n) \bigr] \\
	&\le P_{F_n} \bigl[ \Delta_n(x_n) \le O_p(1) \bigr]
		\to 0 .
\end{align*}\\[-5ex]
\end{proof}

\begin{proof}[Proof of Corollary~\ref{cor:GaussMixtureConsistencyCriterion}]
Since $\|F_n - F_0\|_\infty \le \eps_n \to 0$, it suffices to show that \eqref{eq:GOF-PwrOneCondition} is satisfied. In what follows we use frequently the elementary inequalities
\begin{equation}
\label{eq:Mills.ratio}
	\frac{\phi(x)}{x+1} \le \ \Phi(-x) \ \le \ \frac{\phi(x)}{x}
	\quad\text{for} \ x > 0 ,
\end{equation}
where $\phi(x) := \Phi'(x) = \exp(- x^2/2) / \sqrt{2\pi}$. In particular, as $x \to \infty$,
\begin{align*}
	\Phi(-x) \
	&= \ \exp(- x^2/2 + O(\log x)) \quad\text{and} \\
	C(\Phi(x)) \
	&= \ \log \bigl( O(1) + \log(1/\Phi(-x)) \bigr)
		\ = \ 2 \log(x) - \log(2) + o(1) .
\end{align*}

Now consider two sequences $(x_n)_n$ and $(\mu_n)_n$ tending to $\infty$, and let
$F_0 = \Phi$, $F_n = (1 - \eps_n) \Phi + \eps_n \Phi(\cdot - \mu_n)$. Then the inequalities \eqref{eq:Mills.ratio} imply that
\begin{align*}
	[1 + C(F_0(x_n))] F_0(x_n)(1 - F_0(x_n)) \
	&= \ [2 \log(x_n) + O(1)] \Phi(- x_n) (1 + o(1)) \\
	&= \ \exp[- x_n^2 / 2 + O(\log(x_n))] .
\end{align*}
Moreover,
\[
	F_0(x_n) - F_n(x_n)
	\ = \ \eps_n \bigl( \Phi(\mu_n - x_n) - \Phi(- x_n) \bigr)
	\ = \ \eps_n \Phi(\mu_n - x_n) (1 + o(1)) ,
\]
because $\Phi(- x_n) \le \phi(x_n)/x_n$ while
\[
	\Phi(\mu_n - x_n) \ \ge \ \begin{cases}
		1/2
			& \text{if} \ \mu_n \ge x_n , \\
		\displaystyle
		\frac{\phi(x_n - \mu_n)}{x_n - \mu_n + 1}
		\ \ge \ \frac{\phi(x_n) \exp(\mu_n^2/2)}{x_n + 1}
			& \text{if} \ \mu_n < x_n .
	\end{cases}
\]
Consequently, $\Delta_n(x_n) \to \infty$ if
\begin{equation}
\label{eq:critical.ratio}
	\frac{ n \eps_n \Phi(\mu_n - x_n)}
	{n^{1/2} \exp[- x_n^2/4 + O(\log(x_n))] + O(\log(x_n))}
	\ \to \ \infty .
\end{equation}

In part~(a) with $\eps_n = n^{-\beta + o(1)}$ and $\beta \in (1/2,1)$, we imitate the 
arguments of \cite{Donoho_Jin_2004} and consider
\[
	\mu_n \ = \ \sqrt{2 r \log(n)}
	\quad\text{and}\quad
	x_n \ = \ \sqrt{2 q \log(n)}
\]
with $0 < r < q \le 1$. Then by \eqref{eq:Mills.ratio},
\begin{align*}
	n \eps_n \Phi(\mu_n - x_n) \
	&= \ n^{1 - \beta - (\sqrt{q} - \sqrt{r})^2 + o(1)} , \\
	n^{1/2} \exp[- x_n^2/4 + O(\log(x_n))] \
	&= \ n^{1/2 - q/2 + o(1)} , \\
	O(\log(x_n)) \
	&= \ n^{o(1)} ,
\end{align*}
so the left hand side of \eqref{eq:critical.ratio} equals
\[
	\frac{n^{1 - \beta - (\sqrt{q} - \sqrt{r})^2 + o(1)}}
		{n^{1/2 - q/2 + o(1)} + n^{o(1)}}
	\ = \ \frac{n^{1/2 - \beta + q/2 - (\sqrt{q} - \sqrt{r})^2 + o(1)}}
		{1 + n^{(q-1)/2 + o(1)}}
	\ = \ \frac{n^{1/2 - \beta + 2 \sqrt{r}\sqrt{q} - \sqrt{q}^2/2 - r + o(1)}}
		{1 + n^{(q-1)/2 + o(1)}} .
\]
The exponent in the enumerator is maximal in $q \in (r,1]$ if $\sqrt{q} = \min\{2 \sqrt{r},1\}$, i.e.\ $q = \min\{4r,1\}$, and this leads to
\[
	\begin{cases}
		1/2 - \beta + r & \text{if} \ r \le 1/4 , \\
		1 - \beta - (1 - \sqrt{r})^2 & \text{if} \ r \ge 1/4 .
	\end{cases}
\]
Thus when $\beta \in (1/2,3/4)$ we should choose $\beta - 1/2 < r < 1/4$ and $q = 4r$. 
When $\beta \in [3/4,1)$ we should choose $(1 - \sqrt{1 - \beta})^2 < r < 1$ and $q = 1$.

As to part~(b), we consider the more general setting that $\eps_n = n^{-\beta + o(1)}$ for some 
$\beta \in [1/2,3/4)$, where $\pi_n = \sqrt{n} \eps_n \to 0$. Note that this scenario covers also a part of part~(a), so we establish a connection between the two parts. The constraint that $\pi_n \to 0$ is trivial when $\beta > 1/2$ but relevant when $\beta = 1/2$. Now we consider
\[
	\mu_n \ := \ \sqrt{2 \lambda \log(1/\pi_n)}
	\quad\text{and}\quad
	x_n \ := \ \sqrt{2 q \log(1/\pi_n)}
\]
with arbitrary constants $0 < \lambda < q$. Now
\begin{align*}
	n \eps_n \Phi(\mu_n - x_n) \
	&= \ n^{1/2} \pi_n \Phi(\mu_n - x_n) \\
	&= \ n^{1/2} \pi_n^{1 + (\sqrt{q} - \sqrt{\lambda})^2 + o(1)} , \\
	n^{1/2} \exp \bigl( - x_n^2/4 + O(\log(x_n)) \bigr) \
	&= \ n^{1/2} \pi_n^{q/2 + o(1)} , \\
	O(\log(x_n)) \
	&= \ \pi_n^{o(1)} ,
\end{align*}
so the left hand side of \eqref{eq:critical.ratio} equals
\[
	\frac{n^{1/2} \pi_n^{1 + (\sqrt{q} - \sqrt{\lambda})^2 + o(1)}}
		{n^{1/2}\pi_n^{q/2 + o(1)} + \pi_n^{o(1)}}
	\ = \ \frac{\pi_n^{1 + q/2 - 2 \sqrt{q}\sqrt{\lambda} + \lambda + o(1)}}
		{1 + n^{-1/2} \pi_n^{- q/2 + o(1)}}
	\ = \ \frac{\pi_n^{1 + q/2 - 2 \sqrt{q}\sqrt{\lambda} + \lambda + o(1)}}
		{1 + n^{-1/2 + (\beta - 1/2)q/2 + o(1)}} .
\]
The exponent of $\pi_n$ becomes minimal in $q \in (\lambda,\infty)$ if $q = 4\lambda$. Then we obtain
\[
	\frac{\pi_n^{1 - \lambda + o(1)}}{1 + n^{-1/2 + (2\beta - 1) \lambda + o(1)}}
	\ = \ \frac{\pi_n^{1 - \lambda + o(1)}}{1 + \sqrt{n}^{(4\beta - 2) \lambda - 1 + o(1)}} ,
\]
and this converges to $\infty$ if the limiting exponents of $\pi_n$ and $\sqrt{n}$ are negative. This is the case if $1 < \lambda < 1/(4\beta - 2)$. (Note that $4\beta - 2 < 1$ because $\beta < 3/4$.)
\end{proof}

\begin{proof}[Proof of Lemma~\ref{lem:power.contiguous.1}]
Standard LAN theory implies that $P_{F_n}(S_n) \to 0$ for arbitrary events 
$S_n$ depending on $X_1,\ldots,X_n$ such that $P_{F_0}(S_n) \to 0$. Thus for any fixed $0 <\rho <1/2$, $\varphi_n(X_1,\ldots,X_n) \ne \varphi_{n,\rho}(X_1,\ldots,X_n)$ with asymptotic probability zero, both under the null and under the alternative hypothesis. Hence it suffices to show that
\[
	\limsup_{\rho \to 0} \,
		\limsup_{n \to \infty}
			E_{F_n} \varphi_{n,\rho}(X_1,\ldots,X_n)
	\le \alpha .
\]
But $E_{F_n} \varphi_{n,\rho}(X_1,\ldots,X_n)$ does not change if we replace $f_n$ with the modified density
\[
	f_{n,\rho}(x) := \begin{cases}
		f_n(x),  & \text{if} \ x \not\in [x_{\rho},y_{\rho}] \\
		c_{n,\rho} f_0(x),  & \text{if} \ x \in [x_{\rho},y_{\rho}]
	\end{cases}
\]
with
\[
	c_{n,\rho} := \frac{F_n(y_\rho) - F_n(x_\rho)}{1 - 2\rho} .
\]
This follows from the fact that the distribution function $F_{n,\rho}$ of $f_{n,\rho}$ satisfies 
$F_{n,\rho}(x) = F_n(x)$ for $x \not\in [x_{\rho},y_{\rho}]$,
so the distribution of $\{ \FF_n(x)) : x \not\in [x_{\rho},y_{\rho}] \}$ 
under the alternative hypothesis remains unchanged if we replace $f_n$ with $f_{n,\rho}$. But
\[
	\sqrt{n} (c_{n,\rho} - 1) \to
	\delta_\rho := \frac{A(y_\rho) - A(x_\rho)}{1-2\rho} ,
\]
so
\[
	\sqrt{n} (f_{n,\rho}^{1/2} - f_0^{1/2})
	\to \frac{1}{2} a_\rho f_0^{1/2}
	\quad\text{in} \ L_2 (\lambda)
\]
with
\[
	a_\rho(x) = \begin{cases}
		a(x), & \text{if} \ x \not\in [x_{\rho},y_{\rho}] , \\
		\delta_\rho , & \text{if} \ x \in [x_{\rho},y_{\rho}] .
	\end{cases}
\]
Hence the asymptotic power of the test $\varphi_{n,\rho}$ under the alternative is 
bounded by the asymptotic power of the optimal test of $F_0$ versus $F_{n,\rho}$ at level $\alpha$, so
\[
	\limsup_{n \to \infty}
		E_{F_n} \varphi_{n,\rho}(X_1,\ldots,X_n)
	\le \Phi \bigl( \Phi^{-1}(\alpha) + \|a_\rho\|_{L_2(F_0)} \bigr) .
\]
But
\begin{align*}
	\|a_\rho\|_{L_2(F_0)}^2
	&= \int_{(-\infty,x_\rho) \cup (y_\rho,\infty)} a^2 \, dF_0
		+ \bigl( 1- 2\rho \bigr) \delta_\rho^2 \\
	&= \int_{(-\infty, x_\rho) \cup (y_\rho,\infty)} a^2 \, dF_0
		+ \frac{\bigl( A(y_\rho) - A(x_\rho) \bigr)^2}
			{\bigl( 1-2\rho \bigr)}
\end{align*}
converges to $0$ as $\rho \searrow 0$, so 
$\Phi \bigl( \Phi^{-1}(\alpha) + \|a_\rho\|_{L_2(F_0)} \bigr) \to \alpha$ as $\rho \searrow 0$.
\end{proof}

\begin{proof}[Proof of Theorem~\ref{thm:power.contiguous.2}]
Let $\rho \in (0,1/2)$ be fixed. The test statistic $T_{n,s,\nu}$ for the uniform empirical process may be written as the maximum of $T_{n,s,\nu}^{(\rho,1)}$ and $T_{n,s,\nu}^{(\rho,2)}$, where
\begin{align*}
	T_{n,s,\nu}^{(\rho,1)}
	&:= \sup_{t \in \TT_{n,s} \cap [\rho,1-\rho]} \bigl( n K_s(\GG_n(t),t) - C_\nu(\GG_n(t),t) \bigr) , \\
	T_{n,s,\nu}^{(\rho,2)}
	&:= \sup_{t \in \TT_{n,s} \setminus [\rho,1-\rho]} \bigl( n K_s(\GG_n(t),t) - C_\nu(\GG_n(t),t) \bigr) .
\end{align*}
Here $\TT_{n,s} := (0,1)$ if $s > 0$ and $\TT_n := [\xi_{n:1},\xi_{n:n})$ if $s \le 0$. 
A supremum over the empty set is defined to be $-\infty$. 
The proofs of Theorems~\ref{IntermediateThmPhiDiv} and \ref{thm:PhiDivergenceNull} 
can be easily adapted to show that
\[
	T_{n,s,\nu}^{(\rho,1)} \to_d T_{\nu}^{(\rho,1)}
	\quad\text{and}\quad
	T_{n,s,\nu}^{(\rho,2)} \to_d T_{\nu}^{(\rho,2)} := \max \{ T_{\nu}^{(\rho,2,L)}, T_{\nu}^{(\rho,2,R)} \} ,
\]
where the test statistics $T_\nu^{(\rho,1)}$, $T_\nu^{(\rho,2,L)}$ and $T_\nu^{(\rho,2,R)}$ are defined as in the proof of Lemma~\ref{lem:continuity.limit}. In particular, since $C_\nu(1/2) = 0$ and $\UU(1/2) \ne 0$ almost surely,
\begin{align*}
	\liminf_{n \to \infty} P(T_{n,s,\nu}^{(\rho,1)} > 0)
	&= 1 , \\
	\limsup_{n \to \infty} P(T_{n,s,\nu}^{(\rho,2)} \ge 0)
	&\le \pi_0(\rho) := P(T_{\nu}^{(\rho,2)} \ge 0) .
\end{align*}
Note that $\pi_0(\rho) \to 0$ as $\rho \to 0$ by virtue of Lemma~\ref{lemmaThree}.

Now we consider the goodness-of-fit test statistic $T_{n,s,\nu}(F_0)$. 
It is the maximum of $T_{n,s,\nu}^{(\rho,1)}(F_0)$ and $T_{n,s,\nu}^{(\rho,2)}(F_0)$. 
Here $T_{n,s,\nu}^{(\rho,j)}(F_0)$ is defined as $T_{n,s,\nu}^{(\rho,j)}$, 
where $t \in \TT_{n,s}$ is replaced with $x \in \RR$ if $s > 0$ and $x \in [X_{n:1},X_{n:n})$ 
if $s \le 0$, $[\rho,1 - \rho]$ is replaced with $[x_\rho,y_\rho] = [F_0^{-1}(\rho), F_0^{-1}(1 - \rho)]$, 
and $(\GG_n(t),t)$ is replaced with $(\FF_n(x),F_0(x))$. Under the null hypothesis, 
$T_{n,s,\nu}^{(\rho,j)}(F_0)$ has the same distribution as $T_{n,s,\nu}^{(\rho,j)}$ 
for $j = 1,2$. This convergence and standard LAN theory imply that under the alternative hypothesis,
\begin{align*}
	\liminf_{n \to \infty} P_{F_n} \bigl( T_{n,s,\nu}^{(\rho,1)}(F_0) > 0 \bigr)
	&= 1 , \\
	\limsup_{n \to \infty} P_{F_n} \bigl( T_{n,s,\nu}^{(\rho,2)}(F_0) \ge 0 \bigr)
	&\le \pi_A(\rho) := \Phi \bigl( \Phi^{-1}(\pi_0(\rho)) + \|a\|_{L_2(F_0)} \bigr) .
\end{align*}

With standard empirical process theory one can show that under the alternative hypothesis,
\[
	\sqrt{n} (\FF_n - F_0) \to_d \UU \circ F_0 + A
\]
in the space $\ell^\infty(\RR)$ of bounded functions on $\RR$, equipped 
with the supremum norm $\|\cdot\|_\infty$. Moreover, for arbitrary bounded 
functions $h, h_n$ on $\RR$ such that $\|h_n - h\|_\infty \to 0$,
\[
	n K_s(F_0 + n^{-1/2} h_n, F_0) - C_\nu(F_0 + n^{-1/2} h_n, F_0)
	\to h^2 / [2 F_0(1 - F_0)] - C_\nu(F_0)
\]
uniformly on $[x_\rho,y_\rho]$. By virtue of an extended continuous mapping theorem, e.g.\ \cite{vanderVaart_Wellner_1996}, Theorem 1.11.1, page 67, one can conclude that
\[
	T_{n,s,\nu}^{(\rho,1)}(F_0) \to_d T_\nu^{(\rho,1)}(A) ,
\]
where $T_{\nu}^{(\rho,j)}(A)$ is defined as $T_{\nu}^{(\rho,j)}$ with $\UU + A \circ F_0^{-1}$ in place of $\UU$. Finally, note that the distribution $Q_A$ of $\UU + A \circ F_0^{-1}$ is absolutely continuous with respect to the distribution $Q_0$ of $\UU$, where $\log(dQ_A/dQ_0)$ has distribution $N(-\|a\|_{L_2(F_0)}^2/2, \|a\|_{L_2(F_0)}^2)$ under $Q_0$. This follows from \cite{Shorack_Wellner_2009} (Section 4.1 and especially Theorem 4.1.5, page 157), or \cite{vanderVaart_Wellner_1996} (Section 3.10). Consequently,
\[
	P \bigl( T_{\nu}^{(\rho,2)}(A) \ge 0 \bigr)
	\le \pi_A(\rho) .
\]
All in all, we may conclude that
\[
	P_{F_n} \bigl( T_{n,s,\nu}(F_0) \le 0 \bigr)
	\le P_{F_n} \bigl( T_{n,s,\nu}^{(\rho,1)}(F_0) \le 0 \bigr)
	\to 0 ,
\]
and for fixed $r > 0$,
\begin{align*}
	\limsup_{n \to \infty} P_{F_n} \bigl( T_{n,s,\nu}(F_0) \le r \bigr)
	\le &\limsup_{n \to \infty} P_{F_n} \bigl( T_{n,s,\nu}^{(\rho,1)}(F_0) \le r \bigr) \\
	\le &P \bigl( T_{\nu}^{(\rho,1)}(A) \le r \bigr) \\
	\le &P \bigl( T_{\nu}(A) \le r \bigr) + P \bigl( T_{\nu}^{(\rho,2)}(A) > r \bigr) \\
	\le &P \bigl( T_{\nu}(A) \le r \bigr) + \pi_A(\rho) , \\
	\limsup_{n \to \infty} P_{F_n} \bigl( T_{n,s,\nu}(F_0) \ge r \bigr)
	\le &\limsup_{n \to \infty} P_{F_n} \bigl( T_{n,s,\nu}^{(\rho,1)}(F_0) < r \bigr) \\
		&+ \ \limsup_{n \to \infty} P_{F_n} \bigl( T_{n,s,\nu}^{(\rho,2)}(F_0) \ge r \bigr) \\
	\le &P \bigl( T_{\nu}^{(\rho,1)}(A) \ge r \bigr) + \pi_A(\rho) \\
	\le &P \bigl( T_{\nu}(A) \ge r \bigr) + \pi_A(\rho) .
\end{align*}
Since $\pi_A(\rho) \to 0$ as $\rho \searrow 0$, this proves that $T_{n,s,\nu}(F_0)$ converges in distribution to $T_{\nu}(A)$ under the alternative hypothesis.

The  convergence claimed in the second part of the theorem follows from the first part 
together with convergence of the critical values $\kappa_{n,s,\nu,\alpha}$ to $\kappa_{\nu, \alpha}$.
The inequality claimed in the second part is a consequence of Anderson's \cite{Anderson_1955} inequality or Proposition~\ref{prop:log-concavity.Q} with $h_o := A \circ F_0^{-1}$ and $h(t) := \sqrt{2t(1 - t) (C_\nu(t) + \kappa_{\nu,\alpha})}$.

The third part of the theorem follows from the fact that for any $t \in (0,1)$,
\begin{align*}
	P(T_{\nu}(A) > \kappa_{\nu,\alpha})
	&\ge P \Bigl( \frac{(\UU + A \circ F_0^{-1})^2(t)}{2 t(1 - t)}
			> C_\nu(t) + \kappa_{\nu,\alpha} \Bigr) \\
	&\ge \Phi \biggl( \frac{|A(F_0^{-1}(t))|}{\sqrt{t(1-t)}}
		- \sqrt{2 C_\nu(t) + 2 \kappa_{\nu,\alpha}} \biggr) \\
	&= \Phi \biggl( \frac{|A(F_0^{-1}(t))|}{\sqrt{t(1-t)}}
		- \sqrt{2 C(t)} - b_{\nu,\alpha}(t) \biggr) ,
\end{align*}
where $b_{\nu,\alpha} := (2 \nu D + 2 \kappa_{\nu,\alpha}) \big/ \bigl( \sqrt{2 C + 2 \nu D + 2 \kappa_{\nu,\alpha}} + \sqrt{2 C} \bigr)$ is bounded on $(0,1)$.
\end{proof}

\subsection{Proofs for Subsection~\ref{subsec:ConfBands}}

For notational convenience, we suppress the dependence of the confidence bounds on $s$, $\nu$ and $\alpha$ and just write $a_{n,i}^{\rm BJO}$, $a_{n,i}$, $b_{n,i}^{\rm BJO}$ and $b_{n,i}$.

\begin{proof}[Proof of Theorem~\ref{thm:conf.bands.tails}]
Note first that $H_s(u,t) = \gamma H_s(u/\gamma, t/\gamma)$ for arbitrary $u \ge 0$, $t > 0$ and $\gamma > 0$.

Now we prove the claim for the upper bounds $b_{n,i}^{\rm BJO} = 1 - a_{n,n-i}^{\rm BJO}$ and $b_{n,i} = 1 - a_{n,n-i}$. For any integer $i \in [0,n^\delta]$ let
\[
	x_{n,i} := u_{n,i}/\gamma_n = i/\log\log n .
\]
For fixed $\lambda > 0$ let
\[
	\tilde{b}_{n,i}
	:= u_{n,i} + \lambda \gamma_n (B_s(x_{n,i}) - x_{n,i})
	= \gamma_n \bigl( x_{n,i} + \lambda (B_s(x_{n,i}) - x_{n,i}) \bigr)
	> u_{n,i} .
\]
It follows from $x + s \le B_s(x) \le x + 1 + \sqrt{2x + 1}$ that
\[
	\lambda s \gamma_n
	\le \tilde{b}_{n,i}
	\le \lambda \gamma_n B_s(n^\delta/\log\log n)
	= (\lambda + o(1)) n^{\delta-1} .
\]
On the one hand, if $\lambda > 1$, then it follows from the first inequality in (S.15) that
\[
	n K_s(u_{n,i},\tilde{b}_{n,i})
	\ge n H_s(u_{n,i},\tilde{b}_{n,i})
	= n \gamma_n H_s \bigl( x_{n,i}, x_{n,i} + \lambda (B_s(x_{n,i}) - x_{n,i}) \bigr)
	\ge n \gamma_n \lambda ,
\]
because $H_s \bigl( x_{n,i}, x_{n,i} + t(B_s(x_{n,i}) - x_{n,i}) \bigr)$ is convex in $t$ with values $0$ for $t = 0$ and $1$ for $t = 1$. And if $\lambda < 1$, the second inequality in (S.15) implies that
\begin{align*}
	n K_s(u_{n,i},\tilde{b}_{n,i})
	&\le n H_s(u_{n,i},\tilde{b}_{n,i}) / (1 - \tilde{b}_{in})^+ \\
	&= n \gamma_n H_s \bigl( x_{n,i}, x_{n,i} + \lambda (B_s(x_{n,i}) - x_{n,i}) \bigr)
		/ (1 - \tilde{b}_{n,i}) \\
	&\le n \gamma_n \lambda / \bigl( 1 - (\lambda + o(1)) n^{\delta-1} \bigr)
		= n \gamma_n (\lambda + o(1)) .
\end{align*}
On the other hand, $\kappa_{n,s,\alpha}^{\rm BJ} = (1 + o(1)) n\gamma_n$ and
\begin{align*}
	C_\nu(u_{i,n},\tilde{b}_{i,n}) + \kappa_{n,s,\nu,\alpha}
	& = C_\nu(\tilde{b}_{i,n}) + \kappa_{n,s,\nu,\alpha} \\
	& \begin{cases}
		\le C_\nu(\lambda s\gamma_n) + \kappa_{n,s,\nu,\alpha} = (1 + o(1)) n\gamma_n , \\
		\ge C_\nu \bigl( (\lambda + o(1)) n^{\delta-1} \bigr) + \kappa_{n,s,\nu,\alpha} = (1 + o(1)) n\gamma_n .
	\end{cases}
\end{align*}
Consequently, for any fixed $\lambda > 1$ and sufficiently large $n$,
\[
	n K_s(u_{n,i},\tilde{b}_{n,i})
	> \max \bigl\{ C_\nu(u_{n,i},\tilde{b}_{n,i}) + \kappa_{n,s,\nu,\alpha}, \kappa_{n,s,\alpha}^{\rm BJ} \bigr\}
\]
and thus
\[
	\max\{b_{n,i}^{\rm BJO} - u_{n,i}, b_{n,i} - u_{n,i}\}
	\le \lambda \gamma_n (B_s(x_{n,i}) - x_{n,i})
\]
for all integers $i \in [0,n^\delta]$. Likewise, for any fixed $\lambda \in (0,1)$ and sufficiently large $n$,
\[
	n K_s(u_{n,i},\tilde{b}_{n,i})
	< \min \bigl\{ C_\nu(u_{n,i},\tilde{b}_{n,i}) + \kappa_{n,s,\nu,\alpha}, \kappa_{n,s,\alpha}^{\rm BJ} \bigr\}
\]
and thus
\[
	\min\{b_{n,i}^{\rm BJO} - u_{n,i}, b_{ni} - u_{n,i}\}
	\ge \lambda \gamma_n (B_s(x_{n,i}) - x_{n,i})
\]
for all integers $i \in [0,n^\delta]$.

The differences $u_{n,i} - a_{n,i}^{\rm BJO} = b_{n,n-i}^{\rm BJO} - u_{n,n-i}$ and $u_{n,i} - a_{n,i} = b_{n,n-i} - u_{n,n-i}$ can be treated analogously. For each integer $i \in [1,n^\delta]$ and fixed $\lambda > 0$ let $x_{n,i} = u_{n,i}/\gamma_n = i/\log\log n$ as before and
\[
	\tilde{a}_{n,i}
	:= u_{n,i} + \lambda \gamma_n (A_s(x_{n,i}) - x_{n,i})
	= \gamma_n \bigl( x_{n,i} + \lambda (A_s(x_{n,i}) - x_{n,i}) \bigr)
	< u_{n,i} .
\]
On the one hand, if $\lambda > 1$ and $\tilde{a}_{n,i} > 0$, then $A_s(x_{i,n}) > 0$ and
\[
	n K_s(u_{n,i},\tilde{a}_{n,i})
	\ge n H_s(u_{n,i},\tilde{a}_{n,i})
	= n \gamma_n H_s \bigl( x_{n,i}, x_{n,i} + \lambda (A_s(x_{n,i}) - x_{n,i}) \bigr)
	\ge n \gamma_n \lambda ,
\]
because $H_s \bigl( x_{n,i}, x_{n,i} + t(A_s(x_{n,i}) - x_{n,i}) \bigr)$ is convex in $t \in [0,\lambda]$ with values $0$ for $t = 0$ and $1$ for $t = 1$. And if $\lambda < 1$, then
\begin{align*}
	n K_s(u_{n,i},\tilde{a}_{n,i})
	&\le n H_s(u_{n,i},\tilde{a}_{n,i}) / (1 - u_{in}) \\
	&= n \gamma_n H_s \bigl( x_{n,i}, x_{n,i} + \lambda (A_s(x_{n,i}) - x_{n,i}) \bigr)
		/ (1 - u_{n,i}) \\
	&\le n \gamma_n \lambda / \bigl( 1 - n^{\delta-1} \bigr) .
\end{align*}
On the other hand, $\kappa_{n,s,\alpha}^{\rm BJ} = (1 + o(1)) n\gamma_n$ and
\begin{align*}
	C_\nu(u_{i,n},\tilde{a}_{i,n}) + \kappa_{n,s,\nu,\alpha}
	& = C_\nu(u_{i,n}) + \kappa_{n,s,\nu,\alpha} \\
	& \begin{cases}
		\le C_\nu(n^{-1}) + \kappa_{n,s,\nu,\alpha} = (1 + o(1)) n\gamma_n , \\
		\ge C_\nu(\min\{n^{\delta-1}, 1/2\}) + \kappa_{n,s,\nu,\alpha} = (1 + o(1)) n\gamma_n .
	\end{cases}
\end{align*}
Consequently, for any fixed $\lambda > 1$ and sufficiently large $n$,
\[
	\max\{u_{n,i} - a_{n,i}^{\rm BJO}, u_{n,i} - a_{n,i}\}
	\le \lambda \gamma_n (x_{n,i} - A_s(x_{n,i}))
\]
for all integers $i \in [1, n^\delta]$. Likewise, for any fixed $\lambda \in (0,1)$ and sufficiently large $n$,
\[
	\min\{u_{n,i} - a_{n,i}^{\rm BJO}, u_{ni} - a_{n,i}\}
	\ge \lambda \gamma_n (x_{n,i} - A_s(x_{n,i}))
\]
for all integers $i \in [1, n^\delta]$.
\end{proof}

\begin{proof}[Proof of Theorem~\ref{thm:conf.bands.center}]
We only prove the bounds for $a_{n,i}$ and $b_{n,i}$. The bounds for $a_{n,i}^{\rm BJO}$ and $b_{n,i}^{\rm BJO}$ can be derived analogously with obvious modifications. Moreover, since $u_{n,i} - a_{n,i} = b_{n,n-i} - u_{n,n-i}$, it suffices to prove the bounds for $b_{n,i}$ only. For a fixed factor $\lambda > 0$ and any integer $i \in [n^\delta, n - n^\delta]$ let
\[
	\tilde{b}_{n,i} := u_{n,i} + \lambda \sqrt{ 2 \gamma_n(u_{n,i}) u_{n,i}(1 - u_{n,i})} .
\]
Note that
\begin{align*}
	0  \le \frac{\tilde{b}_{n,i} - u_{n,i}}{u_{n,i}(1 - u_{n,i})}
	   & \le \lambda \sqrt{ 2 n^{-1} (C_\nu(n^{\delta-1}) + \kappa_{\nu,\alpha}) n^{1 - \delta} (1 - n^{\delta-1})^{-1} } \\
	   & = O(n^{-\delta/2} (\log\log n)^{1/2}) ,
\end{align*}
whence
\[
	c_n := \max_{n^\delta \le i \le n - n^\delta} \bigl| \logit(\tilde{b}_{n,i}) - \logit(u_{n,i}) \bigr|
	= o(1) .
\]
On the one hand, the inequalities~(S.14) imply that uniformly in $n^\delta \le i \le n - n^\delta$,
\begin{align*}
	n K_s(u_{n,i}, \tilde{b}_{n,i})
		= n K_{1 - s}(\tilde{b}_{n,i},u_{n,i}) 
	& = (1 + o(1)) n K_2(\tilde{b}_{n,i},u_{n,i}) \\
	& = (1 + o(1)) \lambda^2 (C_\nu(u_{n,i}) + \kappa_{\nu,\alpha}) .
\end{align*}
On the other hand, Lemma~S.10 and Theorem~\ref{thm:PhiDivergenceNull} imply that uniformly in $n^\delta \le i \le n - n^\delta$,
\[
	\bigl| C_\nu(u_{n,i},\tilde{b}_{n,i}) + \kappa_{n,s,\nu,\alpha} - C_\nu(u_{n,i}) - \kappa_{\nu,\alpha} \bigr|
	\le (1 + \nu) c_n + |\kappa_{n,s,\nu,\alpha} - \kappa_{\nu,\alpha}|
	= o(1) .
\]
Consequently, for fixed $\lambda > 1$ and sufficiently large $n$,
\[
	n K_s(u_{n,i},\tilde{b}_{n,i}) > C_\nu(u_{n,i},\tilde{b}_{n,i}) + \kappa_{n,s,\nu,\alpha}
\]
and thus
\[
	b_{n,i} - u_{n,i} \le \lambda \sqrt{2 \gamma_n(u_{n,i}) u_{n,i}(1 - u_{n,i})}
\]
for all integers $i \in [n^\delta, n - n^\delta]$. Likewise, for fixed $\lambda \in (0,1)$ and sufficiently large $n$,
\[
	n K_s(u_{n,i},\tilde{b}_{n,i}) < C_\nu(u_{n,i},\tilde{b}_{n,i}) + \kappa_{n,s,\nu,\alpha}
\]
and thus
\[
	b_{n,i} - u_{n,i}
	\ge \lambda \sqrt{2 \gamma_n(u_{n,i}) u_{n,i}(1 - u_{n,i})}
\]
for all integers $i \in [n^\delta, n - n^\delta]$.
\end{proof}

\paragraph{Acknowledgments.}
The authors owe thanks to David Mason for pointing out the relevance of the tools of Cs{\"o}rg{\H{o}} et al.\ \cite{Csorgo_etal_1986} for some of the results presented here. We are also grateful to G\"unther Walther for stimulating conversations about likelihood ratio tests in nonparametric settings and to Rudy Beran for pointing out the interesting results of Bahadur and Savage \cite{Bahadur_Savage_1956}. Constructive comments of two referees and an associate editor are gratefully acknowledged.

\noindent
The first author was supported in part by the Swiss National Science Foundation.\\
The second author was supported in part by NSF Grant DMS-1104832 and NI-AID grant 2R01 AI291968-04.

\addcontentsline{toc}{section}{References}

\clearpage


\appendix
\setcounter{section}{18}
\setcounter{equation}{0}

\section{Supplement}

References within this document start with `S.' or `(S.'. All other references refer to the main paper.

\subsection{Kolmogorov's upper function test}
\label{subsec:KolmogorovTest}

As mentioned in the introduction, inequality~(1.10) is a consequence of Kolmogorov's integral test for ``upper and lower functions'' for Brownian motion. Let $\WW$ denote standard Brownian motion on $[0,\infty)$ starting at $0$, and let $h$ be a positive continuous function on a nonempty interval $(0,b] \subset (0,\infty)$ such that $h \nearrow$ and $t^{-1/2} h(t) \searrow$.

\begin{proposition}   
Let 
\[
	I_h := \int_0^b t^{-3/2} h(t) \exp( - h^2(t)/2t ) dt .
\]
Then 
\[
	P( \WW(t) \ge h(t) \ \ \text{infinitely often as} \ \ t\searrow 0 ) 
	= \begin{cases}
		0, & \text{if} \ I_h < \infty , \\
	    1, & \text{if} \ I_h = \infty .
	\end{cases}
\]
\end{proposition}

If $I_h < \infty$, then $h$ is an ``upper-class function'' for $\WW$, and if $I_h = \infty$, then $h$ is a ``lower-class function'' for $\WW$. In particular, the function 
\[
	h_{\epsilon}(t) = \sqrt{ 2 t \bigl( \log \log(1/t) + (3/2 + \eps) \log\log\log(1/t) \bigr) } ,
	\quad t \in (0, e^{-e}] ,
\]
is an upper class function for $\WW$ if $\epsilon > 0$, and it is a lower class function for $\WW$ if $\epsilon = 0$. See \cite{Erdos_1942} and \cite{Ito_McKean_1974}, pages 33-36.

\subsection{A general non-Gaussian LIL}
\label{subsec:LIL}

Our conditions and results involve the previously defined function $\logit : (0,1) \to \RR$, $\logit(t) = \log(t/(1-t))$. Its inverse is the logistic function $\ell : \RR \to (0,1)$ given by
\[
	\ell(x) := \frac{e^x}{1 + e^x} = \frac{1}{e^{-x} + 1} ,
\]
and
\[
	\ell'(x) = \ell(x) (1 - \ell(x)) = \frac{1}{e^x + e^{-x} + 2} .
\]
We consider stochastic processes $X = (X(t))_{t \in \TT}$ on subsets $\TT$ of $(0,1)$ 
which have locally uniformly sub-exponential tails in the following sense:

\begin{condition}
\label{C}
There exist real constants $M \ge 1$, $\gamma \ge 0$ and a non-increasing function $L : [0,\infty) \to [0,1]$ such that $L(c) = 1 - O(c)$ as $c \searrow 0$, and
\begin{equation}
\label{eq:subexp}
	P \Bigl( \sup_{t \in [\ell(a), \ell(a+c)] \cap \TT} X(t) > \eta \Bigr)
	\le M \exp(- L(c) \eta) \max(1, L(c)\eta)^{-\gamma}
\end{equation}
for arbitrary $a \in \RR$, $c \ge 0$ and $\eta \in \RR$.
\end{condition}

\begin{theorem}
\label{thm:LIL}
Suppose that $X$ satisfies Condition~\ref{C}. 
For arbitrary $\nu > 1 - \gamma/2$ and $L_0 \in (0,1)$, there exists a real constant 
$M_0 \ge 1$ depending only on $M$, $\gamma$, $L(\cdot)$, $\nu$ and $L_0$ such that
\[
	P \Bigl( \sup_{t \in \TT}
		\bigl( X(t) - C_\nu(t) \bigr)
			> \eta \Bigr)
	\le M_0 \exp( - L_0 \eta)
	\quad\text{for arbitrary} \ \eta \ge 0 .
\]
\end{theorem}

\begin{remark}
Suppose that $X$ satisfies Condition~\ref{C}, where $\inf(\TT) = 0$ and $\sup(\TT) = 1$. For any $\nu > 1 - \gamma/2$, the supremum $T_\nu(X)$ of $X - C - \nu D$ over $\TT$ is finite almost surely. But this implies that
\[
	\lim_{t \to \{0,1\}} \bigl( X(t) - C_\nu(t) \bigr)
	= - \infty
\]
almost surely. For if $1 - \gamma/2 < \nu' < \nu$, then
\[
	X(t) - C_\nu(t)
	= X(t) - C(t) - \nu D(t)
	\le T_{\nu'}(X) - (\nu - \nu') D(t) ,
\]
so the claim follows from $T_{\nu'}(X) < \infty$ almost surely and $D(t) \to \infty$ as $t \to \{0,1\}$.
\end{remark}

\begin{remark}
Our definition of the function $D = \log(1 + C^2)$ may look somewhat arbitrary. 
Indeed, we tried various choices, e.g.\ $D = 2 \log(1 + C)$. Theorem~\ref{thm:LIL} is valid 
for any nonnegative function $D$ on $(0,1)$ such that $D(1-\cdot) = D(\cdot)$ and 
$D(t) / \log\log\log(1/t) \to 2$ as $t \searrow 0$. The special choice $D = \log(1 + C^2)$ 
yields a rather uniform distribution of $\argmax_{(0,1)} (X - C_\nu)$ in case of 
$X(t) = \mathbb{U}(t)^2/(2t(1-t))$ and $\nu$ close to one.
\end{remark}

\begin{proof}[Proof of Theorem~\ref{thm:LIL}]
For symmetry reasons it suffices to prove upper bounds for
\[
	P \Bigl( \sup_{\TT \cap [1/2,1)} (X - C_\nu) > \eta \Bigr) .
\]
Let $(a_k)_{k \ge 0}$ be a sequence of real numbers with $a_0 = 0$ such that
\begin{equation}
\label{eq:Partition.a}
	a_k \to \infty
	\quad\text{and}\quad
	0 < \delta_k := a_{k+1} - a_k \to 0
	\quad\text{as} \ k \to \infty .
\end{equation}
Then it follows from $0 \le \logit(t) - \logit(\ell(a_k)) \le \delta_k$ for $t \in [\ell(a_k),\ell(a_{k+1})]$ and Lemma~\ref{lem:CD.Lipschitz} that
\begin{align*}
	\sup_{\TT \cap [\ell(a_k),\ell(a_{k+1})]} (X - C_\nu)
	&\le \sup_{\TT \cap [\ell(a_k),\ell(a_{k+1})]} X
		\, - C_\nu(\ell(a_k))
		+ (1 + \nu) \delta_k \\
	&\le \sup_{\TT \cap [\ell(a_k),\ell(a_{k+1})]} X
		\, - C_\nu(\ell(a_k))
		+ (1 + \nu) \delta_*
\end{align*}
with $\delta_* := \max_{k \ge 0} \delta_k$. Thus Condition~\ref{C} implies that
\begin{align*}
\lefteqn{P \Bigl( \sup_{\TT \cap [1/2,1)} (X - C_\nu) > \eta \Bigr)  
	\le \ \sum_{k \ge 0}
		P \Bigl( \sup_{\TT \cap [\ell(a_k), \ell(a_{k+1})]}
			(X - C_\nu) > \eta \Bigr)  } \\
	&\le \ \sum_{k \ge 0}
		P \Bigl( \sup_{\TT \cap [\ell(a_k), \ell(a_{k+1})]} X
			\, > \eta - (1 + \nu) \delta_*
				+ C(\ell(a_k)) + \nu D(\ell(a_k)) \Bigr) \\
	&\le \ M \exp( (1 + \nu)\delta_* ) L(\delta_*)^{-\gamma} \exp(- \eta L(\delta_*)) \cdot G ,
\end{align*}
where
\begin{align*}
	G \
	&:= \ \sum_{k \ge 0}
			\exp \bigl( - L(\delta_k) C(\ell(a_k))
				- L(\delta_k) \nu D(\ell(a_k)) \bigr)
			\max \bigl( 1, C(\ell(a_k)) - (1 + \nu) \delta_* \bigr)^{-\gamma} \\
	&= \ \sum_{k \ge 0}
		\Bigl( \log \frac{e}{4 \ell'(a_k)} \Bigr)^{-L(\delta_k)}
		\Bigl( 1 + \Bigl( \log \log \frac{e}{4 \ell'(a_k)} \Bigr)^2 \Bigr)^{- \nu L(\delta_k)} \\
	& \qquad \qquad \qquad \cdot \max \Bigl( 1, \log \log \frac{e}{4 \ell'(a_k)} - (1 + \nu) \delta_* \Bigr)^{-\gamma} .
\end{align*}
Now we define
\[
	a_k \ := \ \delta_* A(k)
	\quad\text{with}\quad
	A(s) \ := \ \frac{s}{\log(e + s)}
\]
for some $\delta_* > 0$ such that $L(\delta_*) \ge L_0 \in (0,1)$. Note that $A(\cdot)$ is a continuously differentiable function on $[0,\infty)$ with $A(0) = 0$, limit $A(\infty) = \infty$ and derivative
\[
	A'(s) \ = \ \frac{1}{\log(e + s)}
		\Bigl( 1 - \frac{s}{(e + s) \log(e + s)} \Bigr)
	\ \in \ \Bigl( 0, \frac{1}{\log(e + s)} \Bigr) .
\]
This implies that \eqref{eq:Partition.a} is indeed satisfied with
\[
	\log a_k = \log k + o(\log k)
	\quad\text{and}\quad
	\delta_k
	\ \le \ \frac{\delta_*}{\log(e + k)}
	\ = \ O(1/\log k) \quad\text{as} \ k \to \infty .
\]
Moreover, for any number $a \ge 0$,
\[
	1 \ \le \ \log \frac{e}{4\ell'(a)}
	\ = \ \log \frac{e(e^a + e^{-a} + 2)}{4}
	\ \in \ \bigl( a + \log(e/4), a + 1 \bigr] .
\]
Consequently, as $k \to \infty$,
\begin{align*}
	\Bigl( \log \frac{e}{4 \ell'(a_k)} & \Bigr)^{-L(\delta_k)}
		\Bigl( 1 + \Bigl(\log \log \frac{e}{4 \ell'(a_k)} \Bigr)^2 \Bigr)^{- \nu L(\delta_k)} \\
		& \qquad \max \Bigl( 1, \log \log \frac{e}{4 \ell'(a_k)} - (1 + \nu) \delta_* \Bigr)^{-\gamma} \\
	&= \ O \bigl( a_k^{-L(\delta_k)}
		\log(a_k)^{- 2\nu L(\delta_k) - \gamma} \bigr) \\
	&= \ O \bigl( k^{-L(\delta_k)} (\log k)^{L(\delta_k)}
		(\log k)^{- 2\nu L(\delta_k) - \gamma} \bigr) \\
	&= \ O \bigl( k^{-1 + O(1/\log k)} (\log k)^{- (2\nu - 1) L(\delta_k) - \gamma} \bigr) \\
	&= \ O \bigl( k^{-1} (\log k)^{- (2\nu - 1 + \gamma + o(1))} \bigr) .
\end{align*}
Since $2\nu - 1 + \gamma > 1$, this implies that $G < \infty$. Hence the asserted 
inequality is true with the constant $M_0 = 2 M \exp((1 + \nu) \delta_*) L(\delta_*)^{-\gamma} \cdot G$.
\end{proof}

\paragraph*{Example 1}
Our first example for a process $X$ satisfying Condition~\ref{C} is squared 
and standardized Brownian bridge:

\begin{lemma}
\label{lem:BB}
Let $\TT = (0,1)$ and $X(t) = \mathbb{U}(t)^2/(2t(1 - t))$ with standard Brownian bridge $\mathbb{U}$. 
Then Condition~\ref{C} is satisfied with $M = 2$, $\gamma = 1/2$ and $L(c) = e^{-c}$.
\end{lemma}

\noindent In particular, Lemma~\ref{lem:BB} and Theorem~\ref{thm:LIL} yield inequality (1.6) for any $\nu > 3/4$.

\begin{proof}[Proof of Lemma~\ref{lem:BB}]
To verify Condition~\ref{C} here, recall that if $\mathbb{W} = (\mathbb{W}(t))_{t \ge 0}$ 
is standard Brownian motion, then $(\mathbb{U}(t))_{t \in (0,1)}$ has the same distribution as 
the stochastic process $\bigl( (1 - t) \mathbb{W}(s(t)) \bigr)_{t \in (0,1)}$ with $s(t) := t/(1 - t) = \exp(\logit(t))$. 
Hence for $a \in \RR$ and $c \ge 0$,
\begin{align*}
	\sup_{t \in [\ell(a), \ell(a+c)]} X(t) \
	& \stackrel{d}{=} \ \sup_{t \in [\ell(a),\ell(a+c)]} \frac{(1 - t)^2 \mathbb{W}(s(t))^2}{2t(1 - t)} \\
	&= \ \sup_{t \in [\ell(a),\ell(a+c)]} \frac{\mathbb{W}(s(t))^2}{2 s(t)} \\
	&= \ \sup_{s \in [e^a, e^{a+c}]} \frac{\mathbb{W}(s)^2}{2 s} \\
	&\stackrel{d}{=} \ \sup_{u \in [e^{-c},1]} \frac{\mathbb{W}(u)^2}{2u} \\
	&\le \ \frac{e^c}{2} \max_{u \in [0, 1]} \mathbb{W}(u)^2 .
\end{align*}
Consequently, the probability that $\sup_{t \in [\ell(a), \ell(a+c)]} X(t)$ is at least $\eta \ge 0$ is bounded by
\begin{align*}
	P \Bigl( \max_{u \in [0,1]} |\mathbb{W}(u)| \ge \sqrt{2\eta e^{-c}} \Bigr) \
	&= \ 2 P \Bigl( \max_{u \in [0,1]} \mathbb{W}(u) \ge \sqrt{2\eta e^{-c}} \Bigr) \\
	&= \ 4 P \bigl( \mathbb{W}(1) \ge \sqrt{2\eta e^{-c}} \bigr) \\
	&= \ 4 \bigl( 1 - \Phi \bigl( \sqrt{2\eta e^{-c}} \bigr) \bigr) ,
\end{align*}
where the second last step follows from a standard argument for processes with independent and symmetrically distributed increments, and $\Phi$ denotes the standard Gaussian distribution function. The well-known inequalities $1 - \Phi(x) \le \exp(-x^2/2)/2$ and $1 - \Phi(x) \le \Phi'(x)/x$ for $x \ge 0$ lead to the bound
\[
	P \Bigl( \sup_{t \in [\ell(a), \ell(a+c)]} X(t) \ge \eta \Bigr) \
	\le \ 2 \exp(- e^{-c} \eta) \max(1, e^{-c} \eta)^{-1/2}
\]
for $\eta \ge 0$, and for negative $\eta$, this bound is obviously true.
\end{proof}

\paragraph*{Example 2}
A second example for Theorem~\ref{thm:LIL} is given by
\[
	X_n(t) \ := \ n K(\GG_n(t),t) ,
	\quad t \in \TT = (0,1) ,
\]
with $K = K_1$.

\begin{lemma}
\label{lem:UniformEP.1}
The stochastic process $X_n$ satisfies Condition~\ref{C} with $M = 2$, $\gamma = 0$ and $L(c) = e^{-c}$.
\end{lemma}

Combining this lemma, Theorem~\ref{thm:LIL} and Donsker's Theorem for the uniform empirical process shows that
\[
	\sup_{t \in (0,1)} \bigl( n K(\GG_n(t),t) - C_\nu(t) \bigr)
	\to_d T_\nu
\]
for any fixed $\nu > 1$. We conjecture that Lemma~\ref{lem:UniformEP.1} is true with $\gamma = 1/2$. This conjecture is supported by refined tail inequalities of \cite{Alfers_Dinges_1984} and \cite{Zubkov_Serov_2013} for binomial distributions.

Before proving Lemma~\ref{lem:UniformEP.1}, recall that for $u \in \RR$ and $t \in (0,1)$,
\begin{align*}
	K(u,t) \ 
	 &:= \ \sup_{\lambda \in \RR} \,
		\bigl( \lambda u - \log(1 - t + t e^\lambda) \bigr) \\
	\ & = \ \begin{cases}
		u \log(u/t) + (1 - u) \log [(1-u)/(1-t)]
			& \text{if} \ u \in [0,1] , \\
		\infty
			& \text{else} .
	\end{cases}
\end{align*}
Indeed, \cite{Hoeffding_1963} showed that for a random variable $Y \sim \mathrm{Bin}(n,t)$ and $u \in \RR$,
\begin{align*}
	P(Y \ge n u) \
	&\le \ \exp \Bigl( - n \sup_{\lambda \ge 0} \,
		\bigl( \lambda u - \log(1 - t + t  e^\lambda) \bigr) \Bigr)
		\ = \ \exp(- n K(u,t))
		\quad\text{if} \ u \ge t , \\
	P(Y \le n u) \
	&\le \ \exp \Bigl( - n \sup_{\lambda \le 0} \,
		\bigl( \lambda u - \log(1 - t + t e^\lambda) \bigr) \Bigr)
		\ = \ \exp(- n K(u,t))
		\quad\text{if} \ u \le t .
\end{align*}

\begin{proof}[Proof of Lemma~\ref{lem:UniformEP.1}]
We imitate and modify a martingale argument of \cite{Berk_Jones_1979} which goes back to \cite{Kiefer_1973}. Note first that $\GG_n(t)/t$ is a reverse martingale in $t \in (0,1)$; that means,
\[
	E \bigl( \GG_n(s)/s \,\big|\, (\GG_n(t'))_{t' \ge t} \bigr)
	\ = \ \GG_n(t)/t
	\quad\text{for} \ 0 < s < t < 1 .
\]
Consequently, for $0 < t < t' < 1$ and $0 \le u \le 1$,
\begin{align*}
	P \Bigl( \inf_{s \in [t,t']} \GG_n(s)/s \le u \Bigr) \
	&= \ \inf_{\lambda \le 0} \,
		P \Bigl( \sup_{s \in [t,t']} \exp(\lambda \GG_n(s)/s - \lambda u) \ge 1 \Bigr) \\
	&\le \ \inf_{\lambda \le 0} \,
		E \exp(\lambda \GG_n(t)/t - \lambda u)
\end{align*}
by Doob's inequality for non-negative submartingales. But $n \GG_n(t) \sim \mathrm{Bin}(n,t)$, so
\begin{align*}
	\inf_{\lambda \le 0} \, E \exp(\lambda \GG_n(t)/t - \lambda u) \
	&= \ \inf_{\lambda \le 0} \, E \exp \bigl( \lambda n \GG_n(t) - n\lambda tu \bigr) \\
	&= \ \exp \Bigl( - n \sup_{\lambda \le 0}
		\bigl( \lambda tu - \log(1 - t + t e^\lambda) \bigr) \Bigr) \\
	&= \ \exp(- n K(tu,t)) .
\end{align*}
Thus
\[
	P \Bigl( \inf_{s \in [t,t']} \GG_n(s)/s \le u \Bigr)
	\ \le \ \exp(- n K(tu,t))
	\quad\text{for all} \ u \in [0,1] .
\]
One may rewrite this inequality as
\[
	P \Bigl( \sup_{s \in [t,t']} n K \bigl( t \min\{\GG_n(s)/s, 1\}, t \bigr) \ge \eta \Bigr)
	\ \le \  \exp(- \eta)
	\quad\text{for all} \ \eta \ge 0 .
\]
For if $\eta > - n \log(1 - t)$, the probability on the left hand side equals $0$. 
Otherwise there exists a unique $u = u(t,\eta) \in [0,1]$ such that $n K(tu,t) = \eta$. But then
\[
	n K \bigl( t \min\{\GG_n(s)/s, 1\}, t \bigr) \ge \eta
	\quad\text{if, and only if,}\quad
	\GG_n(s)/s \le u .
\]
Finally, it follows from the inequalities~\eqref{ineq:K1} for $K(\cdot,\cdot)$ that for $t \le s \le t'$,
\[
	K \bigl( \min\{\GG_n(s),s\}, s \bigr)
	\ = \ K \bigl( s \min\{\GG_n(s)/s, 1\}, s \bigr)
	\ \le \ e^c K \bigl( t \min\{\GG_n(s)/s, 1\}, t \bigr)
\]
with $c := \logit(t') - \logit(t)$. Hence
\[
	P \Bigl( \sup_{s \in [t,t']} n K \bigl( \min\{\GG_n(s), s\}, s \bigr) \ge \eta \Bigr)
	\ \le \ \exp(- e^{-c} \eta)
	\quad\text{for all} \ \eta \ge 0 .
\]

Since $\bigl( \GG_n(t) \bigr)_{t \in (0,1)}$ has the same distribution as 
$\bigl( 1 - \GG_n((1 - t) \, -) \bigr)_{t \in (0,1)}$, and because of the symmetry 
relations $K(s,t) = K(1-s,1-t)$ and $\logit(1 - t) = - \logit(t)$, the previous inequality implies further that
\begin{align*}
	P \Bigl(
	& \sup_{s \in [t,t']} n K \bigl( \max\{\GG_n(s), s\}, s \bigr)
			\ge \eta \Bigr) \\
	&= \ P \Bigl(
		\sup_{s \in [t,t']} n K \bigl( \min\{1 - \GG_n(s), 1 - s\}, 1 - s \bigr)
			\ge \eta \Bigr) \\
	&= \ P \Bigl(
		\sup_{s \in [1-t',1-t]} n K \bigl( \min\{\GG_n(s), s\}, s \bigr)
			\ge \eta \Bigr) \\
	&\le \ \exp(- e^{-c} \eta)
		\quad\text{for all} \ \eta \ge 0 .
\end{align*}
Consequently, since $K(\cdot,s) = \max \bigl\{ K(\min\{\cdot,s\},s), K(\max\{\cdot,s\},s) \bigr\}$,
\[
	P \Bigl( \sup_{s \in [t,t']} n K(\GG_n(s),s) \ge \eta \Bigr)
	\ \le \ 2 \exp(- e^{-c} \eta)
	\quad\text{for all} \ \eta \ge 0 .
\]\\[-7ex]
\end{proof}

\paragraph*{Example 3}
Our third and last example concerns a stochastic process on $\TT_n := \{t_{n,i} : i = 1,2,\ldots,n\}$ with $t_{n,i} = i/(n+1)$: 
\[
	\tilde{X}_n(t_{n,i}) \ := \ (n+1) K(t_{n,i}, \xi_{n:i})
\]
with $K = K_1$.

\begin{lemma}
\label{lem:UniformEP.2}
The stochastic process $\tilde{X}_n$ satisfies Condition~\ref{C} with $M = 2$, $\gamma = 0$ and $L(c) = e^{-c}$.
\end{lemma}

Again one could combine this with Theorem~\ref{thm:LIL} and Donsker's theorem for partial sum processes to show that
\[
	\max_{i=1,\ldots,n} \bigl( (n+1) K(t_{n,i},\xi_{n:i}) - C_\nu(t) \bigr)
	\to_d T_\nu
\]
for any $\nu > 1$.

Our proof of Lemma~\ref{lem:UniformEP.2} involves an exponential inequality for Beta distributions from \cite{Duembgen_1998}:

\begin{lemma}
\label{lem:Beta}
Let $s,t \in (0,1)$, and let $Y \sim \mathrm{Beta}(mt, m(1-t))$ for some $m > 0$. Then
\begin{align*}
	P(Y \le s) \
	&\le \ \inf_{\lambda \le 0} \, E \exp(\lambda Y - \lambda s)
		\ \le \ \exp( - m K(t,s))
		\quad\text{if} \ s \le t , \\
	P(Y \ge s) \
	&\le \ \inf_{\lambda \ge 0} \, E \exp(\lambda Y - \lambda s)
		\ \le \ \exp( - m K(t,s))
		\quad\text{if} \ s \ge t .
\end{align*}
\end{lemma}

\begin{proof}[Proof of Lemma~\ref{lem:UniformEP.2}]
We use a well-known representation of uniform order statistics: 
Let $E_1$, $E_2$, \ldots, $E_{n+1}$ be independent random variables with 
standard exponential distribution, i.e.\ $\mathrm{Gamma}(1)$, and let $S_j := \sum_{i=1}^j E_i$. Then
\[
	(\xi_{n:i})_{i=1}^n
	\  \stackrel{d}{=}  \ (S_i/S_{n+1})_{i=1}^n .
\]
In particular, $\xi_{n:i} \sim \mathrm{Beta}(i,n+1-i) = \mathrm{Beta} \bigl( (n+1)t_{n,i},(n+1)(1 - t_{n,i}) \bigr)$ 
and $E U_{n:i} = t_{n,i}$. Furthermore, for $2 \le k \le n+1$, the random vectors $(S_i/S_k)_{i=1}^{k-1}$ 
and $(S_i)_{i=k}^{n+1}$ are stochastically independent. This implies that $(\xi_{n:i}/t_{n,i})_{i=1}^n$ 
is a reverse martingale, because for $1 \le j < k \le n$,
\[
	E \Bigl( \frac{\xi_{n:j}}{t_{n,j}}
		\,\Big|\, (S_i)_{i=k}^{n+1} \Bigr)
	\ = \ E \Bigl( \frac{S_j}{t_{n,j} S_k} \cdot \frac{S_k}{S_{n+1}}
		\,\Big|\, (S_i)_{i=k}^{n+1} \Bigr)
	\ = \ \frac{j}{t_{n,j} k} \cdot \frac{S_k}{S_{n+1}}
	\ = \ \frac{\xi_{n:k}}{t_{nk}} .
\]
Consequently, for $1 \le j \le k \le n$ and $0 < u < 1$, it follows from 
Doob's inequality and Lemma~\ref{lem:Beta} that
\begin{align*}
	P \Bigl( \min_{j \le i \le k} \frac{\xi_{n:i}}{t_{n,i}} \le u \Bigr) \
	&= \ \inf_{\lambda < 0} \,
		P \Bigl( \min_{j \le i \le k}
			\exp \Bigl( \lambda \frac{\xi_{n:i}}{t_{n,i}} - \lambda u \Bigr) \ge 1 \Bigr) \\
	&\le \ \inf_{\lambda < 0} \,
		E \exp \bigl( \lambda \xi_{n:j} - \lambda u t_{n,j} \bigr) \\
	&\le \ \exp \bigl( - (n+1) K(t_{n,j}, t_{n,j} u) \bigr) .
\end{align*}
Again one may reformulate the previous inequalities as follows: For any $\eta > 0$,
\[
	P \Bigl( \max_{j \le i \le k} (n+1)
		K \Bigl( t_{n,j}, t_{n,j} \min \Bigl\{ \frac{\xi_{n:i}}{t_{n,i}}, 1 \Bigr\} \Bigr)
		\ge \eta \Bigr)
	\ \le \ \exp(- \eta) .
\]
But the inequalites~\eqref{ineq:K1} for $K(\cdot,\cdot)$ imply that for $j \le i \le k$,
\[
	K \bigl( t_{n,i}, \min\{\xi_{n:i}, t_{n,i}\} \bigr)
	\ \le \
	e^c K \Bigl( t_{n,j}, t_{n,j} \min \Bigl\{ \frac{\xi_{n:i}}{t_{n,i}}, 1 \Bigr\} \Bigr)
\]
with $c := \logit(t_{nk}) - \logit(t_{n,j})$. Consequently,
\[
	P \Bigl( \max_{j \le i \le k} (n+1)
		K \bigl( t_{n,i}, \min\{\xi_{n:i},t_{n,i}\} \bigr) \ge \eta \Bigr)
	\ \le \ \exp(- e^{-c} \eta)
	\quad\text{for all} \ \eta > 0 .
\]
Since $(1 - \xi_{n:n+1-i})_{i=1}^n$ has the same distribution as $(\xi_{n:i})_{i=1}^n$, 
a symmetry argument as in the proof of Lemma~\ref{lem:UniformEP.1} reveals that
\[
	P \Bigl( \max_{j \le i \le k} (n+1) K(t_{n,i},\xi_{n:i}) \ge \eta \Bigr)
	\ \le \ 2 \exp(- e^{-c} \eta)
	\quad\text{for all} \ \eta > 0 .
\]\\[-7ex]
\end{proof}

\subsection{Auxiliary functions and (in)equalities}
\label{subsec:AuxiliaryFunctions}

\paragraph*{Inequalities involving the logit function}

Recall first that for arbitrary numbers $x > 0$ and $\gamma \in \RR$, the representation $x^\gamma = \exp(\gamma \log x)$ implies that
\[
	\exp \bigl( - |\gamma| |\log x| \bigr) \le x^\gamma \le \exp \bigl( |\gamma| |\log x| \bigr) .
\]
Now we consider arbitrary numbers $t,u \in (0,1)$. Note that either $u/t < 1 < (1 - u)/(1 - t)$ or $u/t \ge 1 \ge (1-u)/(1-t)$. Consequently,
\begin{equation}
\label{eq:LogitUT.1}
	\bigl| \log(u/t) \bigr| + \bigl| \log[(1 - u)/(1 - t)] \bigr|
	= \bigl| \logit(u) - \logit(t) \bigr| ,
\end{equation}
and this implies that
\begin{equation}
\label{ineq:LogitUT.2}
	(u/t)^\gamma, [(1 - u)/(1 - t)]^\gamma
	\ \in \ \bigl[ e^{-|\gamma| c}, e^{|\gamma| c} \bigr]
	\quad\text{with} \ c := \bigl| \logit(u) - \logit(t) \bigr| .
\end{equation}

In the proofs of Theorem~\ref{thm:LIL} and Theorem~2.1, we utilize the following continuity properties of the functions $C, D : (0,1) \to [0,\infty)$.

\begin{lemma}
\label{lem:CD.Lipschitz}
For arbitrary $s,t \in (0,1)$,
\[
	\bigl| D(s) - D(t) \bigr|
	\le \bigl| C(s) - C(t) \bigr|
	\le \bigl| \logit(s) - \logit(t) \bigr| .
\]
\end{lemma}

\begin{proof}
Since $D = \log(1 + C^2)$, the first inequality follows from $d \log(1 + x^2)/dx = 2x/(1 + x^2) \in [0,1]$ for $x \ge 0$. As to the second inequality, if $s(1-s) \le t(1-t)$, then
\begin{align*}
	0 \le C(s) - C(t)
	&= \log \left( \log \Bigl( \frac{e}{4s(1-s)} \Bigr) \Big/ \log \Bigl( \frac{e}{4t(1-t)} \Bigr) \right) \\
	&= \log \left( 1 + \log \Bigl( \frac{t(1-t)}{s(1-s)} \Bigr) \Big/ \log \Bigl( \frac{e}{4t(1-t)} \Bigr) \right) \\
	&\le \log \Bigl( \frac{t(1-t)}{s(1-s)} \Bigr) \\
	&\le \max \Bigl\{ \log \Bigl( \frac{t}{s} \Bigr), \log \Bigl( \frac{1 - t}{1 - s} \Bigr) \Bigr\} \\
	&\le \bigl| \logit(s) - \logit(t) \bigr| ,
\end{align*}
because $\log(t/s) \ge 0 \ge \log((1 - t)/(1-s))$ or $\log(t/s) \le 0 \le \log((1 - t)/(1 - s))$.
\end{proof}

\paragraph*{The divergences $K_s$}
Recall that the divergences $K_s$ can be written as $K_s(u,t) = t \phi_s (u/t) + (1-t) \phi_s [(1-u)/(1-t)]$ with certain auxiliary functions $\phi_s : (0,\infty) \to [0,\infty)$ and their limits $\phi_s(0) := \lim_{x \searrow 0} \phi_s(x) \in (0,\infty]$. In particular,
\[
	K_s(u,t) = K_s(1-u,1-t) .
\]
Precisely, $\phi_s$ is given by $\phi_s(1) = 0 = \phi_s'(1)$ and $\phi_s''(x) = x^{s-2}$.
Any twice continuously differentiable function $f : (0,\infty) \to \RR$ may be written as
\begin{equation}
\label{eq:Taylor2F}
	f(x) = f(1) + f'(1) (x - 1) + \int_1^x (x - u) f''(u) \, du .
\end{equation}
For $\phi_s$ this yields the representation
\begin{equation}
\label{eq:PhiS.1}
	\phi_s(y) = \int_1^y (y - x) x^{s-2} \, dx
\end{equation}
for $y > 0$. Starting from this representation, elementary calculations yield the explicit formulae (3.18) for $\phi_s$ and (1.7) for $K_s$.

Plugging in the representation \eqref{eq:PhiS.1} in the representation of $K_s$ in terms of $\phi_s$ and transforming the two integrals appropriately leads to the representation
\begin{equation}
\label{eq:KS.1}
	K_s(u,t)
	= \int_t^u (u - x) \bigl[ t^{1-s} x^{s-2} + (1 - t)^{1-s} (1 - x)^{s-2} \bigr] \, dx .
\end{equation}
In particular,
\[
	K_2(u,t) = \int_t^u (u - x) [t^{-1} + (1 - t)^{-1}] \, dx = \frac{(u - t)^2}{2t(1-t)} .
\]
Comparing \eqref{eq:KS.1} with \eqref{eq:Taylor2F} reveals that
\begin{eqnarray}
\label{eq:KS.2}
	&& K_s(t,t) = 0, \quad
		\frac{\partial}{\partial u}\Big|_{u=v} K_s(u,t) = 0, 
		\quad\text{and}\quad \\
	&& \frac{\partial^2}{\partial u^2} K_s(u,t)
		= t^{1-s} u^{s-2} + (1 - t)^{1-s} (1 - u)^{s-2} .
\end{eqnarray}
Integrating the latter formula leads to
\begin{equation}
\label{eq:KS.3}
	\frac{\partial}{\partial u} K_s(u,t) = \begin{cases}
		\logit(u) - \logit(t)
			& \text{if} \ s = 1 , \\[1ex]
		\displaystyle
		\frac{(u/t)^{s-1} - [(1 - u)/(1 - t)]^{s-1}}{s-1}
			& \text{if} \ s \ne 1 .
	\end{cases}
\end{equation}

Another interesting identity follows from \eqref{eq:PhiS.1} via the substitution $\tilde{x} = 1/x$:
\begin{equation}
\label{eq:PhiS.2}
	\phi_s(y) = y \phi_{1-s}(1/y)
\end{equation}
for $y > 0$, and this leads to
\begin{equation}
\label{eq:KS.4}
	K_s(u,t) = K_{1-s}(t,u) .
\end{equation}

\paragraph*{Some particular inequalities for $K = K_1$}

For fixed $v \in (0,1)$ and arbitrary $0 < t < t' < 1$,
\begin{equation}
\label{ineq:K1}
	\frac{K(0,t')}{K(0,t)}, \frac{K(t'v,t')}{K(tv,t)} , \frac{K(t',t'v)}{K(t,tv)}
	\ \in \ \Bigl( \frac{t'}{t}, \frac{t'(1-t)}{(1 - t')t} \Bigr) .
\end{equation}

To prove these inequalities, note that on the one hand,
\[
	K(tv,t)
	\ = \ \int_{tv}^{t} \frac{\partial K_0(x,tv)}{\partial x} \, dx
	\ = \ \int_{tv}^{t} \frac{(x - tv)}{x(1 - x)} \, dx
	\ = \ \int_v^1 \frac{t(y-u)}{y(1 - ty)} \, dy .
\]
These formulae remain true if we replace $v$ with $0$. On the other hand,
\[
	K(t,tv)
	\ = \ \int_{tv}^{t} (t - x) \frac{\partial^2}{\partial x^2} K(x,tv) \, dx
	\ = \ \int_{tv}^{t} \frac{(t - x)}{x(1 - x)} \, dx
	\ = \ \int_v^1 \frac{t(1 - y)}{y(1 - ty)} \, dy .
\]
But for any $y \in (0,1)$,
\[
	\frac{\partial}{\partial t} \log \frac{t}{1 - ty}
	\ = \ \frac{1}{t(1 - ty)}
	\ \in \ \Bigl( \frac{1}{t}, \frac{1}{t(1-t)} \Bigr)
	\ = \ \bigl( \log'(t), \logit'(t) \bigr) .
\]
Thus for $0 < t < t' < 1$,
\[
	\frac{t'}{1 - t' y} \Big/ \frac{t}{1 - ty}
	\ \in \ \Bigl( \frac{t'}{t}, \frac{t'(1 - t)}{(1 - t') t} \Bigr) ,
\]
and this entails the asserted inequalities for the three ratios $K(0,t')/K(0,t)$, $K(t'v,t')/K(tv,t)$ and $K(t',t'v)/K(t,tv)$.

\paragraph*{Relating $K_s$ and $K_2$}

Starting from \eqref{eq:KS.1}, we may write
\begin{align*}
	K_s(u,t)
	&= \int_t^u (u - x) \bigl[ t^{-1} (x/t)^{s-2} + (1 - t)^{-1} [(1 - x)/(1 - t)]^{s-2} \bigr] \, dx \\
	&= \int_u^t (x - u) \bigl[ t^{-1} (x/t)^{s-2} + (1 - t)^{-1} [(1 - x)/(1 - t)]^{s-2} \bigr] \, dx .
\end{align*}
Note that either $t < u$ and $u/t \ge x/t \ge 1 \ge (1 - x)/(1 - t) \ge (1 - u)/(1 - t)$, or $t \ge u$ and $u/t \le x/t \le 1 \le (1 - x)/(1 - t) \le (1 - u)/(1 - t)$. Hence, it follows from these representations of $K_s(u,t)$ and the inequalities \eqref{ineq:LogitUT.2} that
\begin{equation}
\label{ineq:KSK2}
	\frac{K_s(u,t)}{K_2(u,t)} \in \bigl[ e^{-|s-2| c}, e^{|s-2| c} \bigr]
	\quad\text{with} \ c := \bigl| \logit(u) - \logit(t) \bigr| ,
\end{equation}
where $K_s(t,t)/K_2(t,t) := 1$.

\paragraph*{Some bounds for $\phi_s$ and $K_s$}

In what follows, we restrict our attention to parameters $s \in [-1,2]$. The next lemma provides lower bounds for $\phi_s$.

\begin{lemma}
\label{lem:CoolS1}
Let $s \in [-1,2]$. Then
\[
	\phi_s(1 + x)
	\ge \frac{x^2}{2(1 + ax)}
		\quad \text{for} \ x > -1 ,
\]
where $a := (2 - s)/3 \in [0,1]$.
\end{lemma}

Lemma~\ref{lem:CoolS1} implies useful bounds for $K_s$.

\begin{lemma}
\label{lem:CoolS2}
Let $s \in [-1,2]$. Then for $t,u \in (0,1)$,
\[
	K_s(u,t) \ \ge \ \frac{\delta^2}{2(t + a\delta)(1 - t - a\delta)} ,
\]
where $\delta := u - t \in (-t,1-t)$ and $a := (2 - s)/3 \in [0,1]$. Moreover, for any $\gamma > 0$, the inequality $K_s(u,t) \le \gamma$ implies that
\[
	|\delta| \ \le \ \begin{cases}
		\sqrt{2\gamma \,t(1 - t)} \, + 2 |1 - 2\,t| a \gamma , \\
		\sqrt{2\gamma u(1 - u)} + 2 |1 - 2u| (1 - a) \gamma .
	\end{cases}
\]
\end{lemma}

\begin{proof}[Proof of Lemma~\ref{lem:CoolS1}]
The asserted inequality reads $\phi_s(1 + x) \ge h_a(x)$ for $x > -1$ with the auxiliary function $h_a(x) := 2^{-1} x^2 /(1 + ax)$. Elementary calculations reveal that $h_a(0) = 0 = h_a'(0)$ and $h_a''(x) = (1 + ax)^{-3}$. On the other hand, $\phi_s(1) = 0 = \phi_s'(1)$ and $\phi_s''(1 + x) = (1 + x)^{s-2} = (1 + x)^{-3a}$. Consequently, it suffices to show that $\phi_s''(1 + \cdot) \ge h_a''$, that is,
\[
	(1 + x)^{-3a} \ge (1 + ax)^{-3}
\]
for $x > -1$. This is equivalent to the inequality
\[
	- a \log(1 + x) \ge - \log(1 + ax) .
\]
But this inequality follows from convexity of $-\log$, because
\begin{align*}
	-\log(1 + ax)
	& = -\log[a \cdot (1 + x) + (1 - a) \cdot 1] \\
	&  \le -a \log(1 + x) - (1 - a) \log(1) 
	 = -a \log(1 + x) .
\end{align*}\\[-5ex]
\end{proof}

\begin{proof}[Proof of Lemma~\ref{lem:CoolS2}]
It follows from Lemma~\ref{lem:CoolS1} that
\begin{align*}
	K_s(u,t)
	&= t \phi_s(1 + \delta/t) + (1 - t) \phi_s[1 - \delta/(1 - t)] \\
	&\ge \ \frac{t (\delta/t)^2}{2(1 + a \delta/t)}
		+ \frac{(1 - t) [\delta/(1-t)]^2}{2 (1 - a \delta/(1 - t))} \\
	&= \ \frac{\delta}{2(t + a \delta)}
			+ \frac{\delta^2}{2 (1 - t - a \delta)}
		= \frac{\delta^2}{2(t + a\delta)(1 - t - a\delta)} .
\end{align*}
As a consequence, the inequality $K_s(u,t) \le \gamma$ implies that
\[
	\delta^2
	\le 2 \gamma (t + a\delta)(1 - t - a\delta)
	\le 2\gamma t(1 - t) + 2 \delta (1 - 2t) a \gamma .
\]
With $b := a (1 - 2t)$, this leads to $\delta^2 - 2\delta b\gamma \le 2 \gamma t(1 - t)$, that is,
\[
	(\delta - b \gamma)^2 \le 2\gamma t(1 - t) + b^2 \gamma^2 .
\]
Consequently,
\[
	|\delta|
	\le |b| \gamma + \sqrt{2\gamma t(1 - t) + b^2 \gamma^2}
	\le \sqrt{2\gamma t(1 - t)} + 2 |b| \gamma
	= \sqrt{2\gamma t(1 - t)} + 2 |1 - 2t| a \gamma , 
\]
because $\sqrt{x+y} \le \sqrt{x} + \sqrt{y}$ for $x,y \ge 0$. The second inequality for $|\delta|$ follows from the first one and the identity \eqref{eq:KS.4}: Since $K_s(u,t) = K_{1-s}(t,u)$, and since $(2 - (1 - s))/3 = (s+1)/3 = 1 - a$, it follows from $K_s(u,t) \le \gamma$ that
\[
	|\delta| \le \sqrt{2\gamma u(1 - u)} + 2 |1 - 2u| (1 - a) \gamma .
\]\\[-5ex]
\end{proof}

\paragraph*{Approximating $K_s$ close to $(0,0)$}

The following bounds show that $K_s(u,t)$ can be approximated by a simpler function if $u,t$ are close to $0$: For $s \in [-1,2]$ and $u,t \in (0,1)$,
\begin{equation}
\label{eq:Ks.Hs}
	t \phi_s(u/t) \le K_s(u,t) \le t \phi_s(u,t) / (1 - \max\{u,t\}) .
\end{equation}
If $s \in (0,2]$, then \eqref{eq:Ks.Hs} is even true for $u = 0$ and reads as $t/s \le K_s(0,t) \le (t/s) / (1 - t)$. To verify \eqref{eq:Ks.Hs}, recall that $K_s(u,t)$ is the sum of the nonnegative terms $t \phi_s(u/t)$ and $(1 - t) \phi_s[(1-u)/(1 - t)]$. If $u < t$, then
\[
	t \phi_s(u/t)
	= t \int_{u/t}^1 (r - u/t) r^{s-2} \, dr
	\ge t \int_{u/t}^1 (r - u/t) \, dr
	= (u - t)^2/(2t) ,
\]
because $r \le 1$ and $s - 2 \le 0$, whereas
\begin{align*}
	(1 - t) \phi_s[(1-u)/(1-t)]
	&= (1 - t) \int_1^{(1-u)/(1-t)} [(1 - u)/(1 - t) - r] r^{s-2} \, dr \\
	&\le (1 - t) \int_1^{(1-u)/(1-t)} [(1 - u)/(1 - t) - r] \, dr \\
	&= (u - t)^2/[2(1-t)]
		= (u - t)^2/(2t) \cdot t/(1-t) ,
\end{align*}
because $r \ge 1$. If $t < u$, we use the identity~\eqref{eq:PhiS.2} to verify that
\[
	t \phi_s(u/t) = u \phi_{1-s}(t/u) \ge (u - t)^2/(2u)
\]
and
\[
	(1-t) \phi_s[(1-u)/(1-t)] = (1-u) \phi_{1-s}[(1-t)/(1-u)] \le (u - t)^2/(2u) \cdot u/(1 - u) ,
\]
because $(1 - s) - 2 = -s-1 \le 0$.

The next lemma summarizes some properties of the function $(x,y) \mapsto y \phi_s(x/y)$ which appears in \eqref{eq:Ks.Hs}.

\begin{lemma}
\label{lem:Hs.As.Bs}
For $s \in [-1,2]$ and $x,y > 0$ let
\[
	H_s(x,y) := y \phi_s(x/y) = x \phi_{1-s}(y/x) .
\]
This defines a continuous, convex function $H_s : (0,\infty) \times (0,\infty) \to [0,\infty)$. For $x,\lambda > 0$, $H_s(x,\lambda x) = x \phi_{1-s}(\lambda)$, and $H_s(x,x) = 0$. In case of $s > 0$, the function $H_s$ can be extended continuously to $[0,\infty) \times (0,\infty)$ via $H_s(0,y) := y/s$, and in case of $0 < s < 1$, it can be extended continuously to $[0,\infty) \times [0,\infty)$ via $H_s(x,0) := x/(1 - s)$.

For $x \ge 0$ let
\begin{align*}
	a_s(x) &:= \begin{cases}
		0
			& \text{if} \ x = 0 , \\
		\inf \{ y \in (0,x) : H_s(x,y) \le 1 \}
			& \text{else} ,
		\end{cases} \\
	b_s(x) &:= \begin{cases}
		s^+
			& \text{if} \ x = 0 , \\
		\max\{ y > x : H_s(x,y) \le 1 \}
			& \text{else} .
		\end{cases}
\end{align*}
This defines continuous functions $a_s, b_s : [0,\infty) \to [0,\infty)$ where $a_s$ is convex with $a_s(x) = 0$ if and only if $x \le (1 - s)^+$, and $b_s$ is concave. Moreover, for fixed $x \ge 0$, $a_s(x)$ and $b_s(x)$ are non-decreasing in $s \in [-1,2]$ and satisfy the inequalities
\begin{align*}
	x + \tilde{a} - \sqrt{2x + \tilde{a}^2} \le a_s(x) & \le x + 1 - \sqrt{2x + 1} , \\
	x + \max \bigl\{ s, \sqrt{2x} \bigr\} \le b_s(x) & \le x + \tilde{a} + \sqrt{2x + \tilde{a}^2} ,
\end{align*}
where $\tilde{a} := (1+s)/3 \in [0,1]$.
\end{lemma}

This lemma implies that $a_s(x)/x \to 0$ and $b_s(x)/x \to \infty$ as $x \searrow 0$, whereas $a_s(x) = x - \sqrt{2x} + O(1)$ and $b_s(x) = x + \sqrt{2x} + O(1)$ as $x \to \infty$.

\begin{remark}
Since $K_s(u,t) = H_s(u,z) + H_s(1 - u, 1 - t)$, Lemma~\ref{lem:Hs.As.Bs} implies that $K_s$ is a convex function on $(0,1) \times (0,1)$ with $K_s(t,t) = 0$ for all $t \in (0,1)$. Joint convexity of the functions $(u,v) \mapsto K_s (u,v)$ is a very special case of \cite{Simon_2011}, Theorem~16.3.
\end{remark}

\begin{proof}[Proof of Lemma~\ref{lem:Hs.As.Bs}]
Convexity of $H_s$ follows from the fact that for $x,y > 0$, the Hessian matrix of $H_s$ at $(x,y)$ equals
\[
	x^{s-1} y^{-s} \begin{bmatrix} y/x, & -1 \\ -1, & x/y \end{bmatrix} ,
\]
which is positive semidefinite.

For $x > 0$, it follows from the formula $H_s(x,y) = x \phi_{1-s}(y/x)$ and $\phi_{1-s} : [1,\infty) \to [0,\infty)$ being increasing and bijective that $b_s(x)$ is the unique number $y \in (x,\infty)$ such that $H_s(x,y) = 1$. More precisely, for $y > x$, $b_s(x) \le y$ is equivalent to $H_s(x,y) \ge 1$, and $b_s(x) \ge y$ is equivalent to $H_s(x,y) \le 1$.

If $s \le 0$, then for any fixed $y > 0$, $H_s(x,y) = y \phi_s(x/y) \to \infty$ as $x \searrow 0$, whence $b_s(x) \to 0$ as $x \searrow 0$. If $s > 0$, then $H_s(x,s) = s \phi_s(x/s)$ is strictly decreasing in $x \in [0,s]$ with $H_s(0,s) = 1$, whence $b_s(x) \ge s$ for all $x \ge 0$. On the other hand, for any $y > s$, $H_s(x,y) = y \phi_s(x/y) \to y/s > 1$ as $x \searrow 0$, whence $b_s(x) \to s$ as $x \searrow 0$. This shows that $b_s$ is continuous at $0$.

Convexity of $H_s$ implies that $b_s$ is concave and thus continuous on $(0,\infty)$. Together with continuity at $0$, this implies that $b_s$ is continuous and concave on $[0,\infty)$.

For $x > 0$ and $y \in [0,x]$, it follows from $\phi_{1-s} : [0,1] \to [0,1/(1 - s)^+]$ being decreasing and bijective that $a_s(x) = 0$ if $x \le (1 - s)^+$, and for $x > (1 - s)^+$, $a_s(x)$ is the unique number $y \in (0,x)$ such that $H_s(x,y) = 1$. More precisely, for $y \in (0,x)$, $a_s(x) \ge y$ is equivalent to $H_s(x,y) \ge 1$, and $a_s(x) \le y$ is equivalent to $H_s(x,y) \le 1$. Convexity of $H_s$ implies that $a_s$ is convex too, and since $0 \le a_s(x) < x$ for all $x > 0$, $a_s$ is a convex and continuous function on $[0,\infty)$.

By continuity, it suffices to verify the remaining claims for $x > 0$. It follows from Lemma~\ref{lem:CoolS1} that for $x,y > 0$,
\[
	H_s(x,y) = y \phi_s(x/y)
	\ge \frac{y (x/y - 1)^2}{2(1 - a + ax/y)}
	= \frac{(x - y)^2}{2(\tilde{a}y + ax)} ,
\]
where $a = (2 - s)/3 \in [0,1]$ and $\tilde{a} = 1 - a = (1+s)/3$. Consequently, the inequality $H_s(x,y) \le 1$ implies that $(y - x)^2 \le 2(\tilde{a}y + ax)$, and this is equivalent to $(y - x - \tilde{a})^2 \le 2x + \tilde{a}^2$, that is,
\[
	a_s(x) \ge x + \tilde{a} - \sqrt{2x + \tilde{a}^2}
	\quad\text{and}\quad
	b_s(x) \le x + \tilde{a} + \sqrt{2x + \tilde{a}^2} .
\]

For $0 < x < y$, $H_s(x,y) = y \int_{x/y}^1 (r - x/y) r^{s-2} \, dr$ is monotone decreasing in $s \in [-1,2]$. By construction of $b_s(x)$, this entails that $b_s(x)$ is monotone increasing in $s \in [-1,2]$. Consequently, $b_s(x) \ge b_{-1}(x) = x + \sqrt{2x}$, because
\[
	H_{-1}(x,y) = x \phi_2(y/x) = (y - x)^2/(2x) = 1
	\quad\text{if and only if}\quad
	y = x \pm \sqrt{2x} .
\]
Furthermore, if $s > 0$, then $H_s(0,s) = 1$, and $H_s(x,x + \sqrt{2x}) \le 1$ for all $x > 0$. For $x_o = s^2/2$, $x_o + \sqrt{2x_o} = x_o + s$. By convexity of $H_s$,
\[
	H_s(x, x + s) \le (1 - x/x_o) H_s(0,s) + (x/x_o) H_s(x_o,x_o + s) \le 1
\]
for $0 \le x \le x_o$, whence $b_s(x) \ge x + s$ for $0 \le x \le x_o$. Since $x + \sqrt{2x} \ge x + s$ if and only if $x \ge x_o$, this shows that $b_s(x) \ge x + \max \bigl\{ s, \sqrt{2x} \bigr\}$.

For $0 < y < x$, $H_s(x,y) = y \int_1^{x/y} (x/y - r) r^{s-2} \, dr$ is monotone increasing in $s \in [-1,2]$, so $a_s(x)$ is monotone increasing by its construction. Consequently $a_s(x) \le a_2(x) = x+1 - \sqrt{2x+1}$, because
\[
	H_2(x,y) = y \phi_2(x/y) = (y - x)^2/(2y) = 1
	\quad\text{if and only if}\quad
	y = x +1 \pm \sqrt{2x + 1} .	
\]\\[-5ex]
\end{proof}

\subsection{Further proofs for Section~2}
\label{subsec:LogConcavityQ_SelfSimilarityBB}

\begin{proof}[Proof of Proposition~4.13]
Log-concavity of $G_1$ follows from the facts that $G_1(x) = Q(\mathcal{B}_1(x))$ with the closed set $\mathcal{B}_1(x) := \bigl\{ g \in \mathcal{C}[0,1] : |x h_o + g| \le h \bigr\}$, and that $(1 - \lambda) \mathcal{B}_1(x_0) + \lambda \mathcal{B}_1(x_1) \subset \mathcal{B}_1((1 - \lambda) x_0 + \lambda x_1)$ for $x_0, x_1 \in \RR$ and $\lambda \in (0,1)$. Indeed, if $g_0 \in \mathcal{B}_1(x_0)$ and $g_1 \in \mathcal{B}_1(x_1)$, then
\[
	\bigl| (1 - \lambda) x_0 h_o + \lambda x_1 h_o + (1 - \lambda) g_0 + \lambda g_1 \bigr|
	\le (1 - \lambda) |x_0 h_o + g_0| + \lambda |x_1 h_o + g_1|
	\le h .
\]
Similarly, $G_2(x) = Q(\mathcal{B}_2(x))$ with $\mathcal{B}_2(x) := \bigl\{ g \in \mathcal{C}[0,1] : |g| \le \sqrt{h + x h_o} \bigr\}$, and for $x_0,x_1 \ge 0$ and $\lambda \in (0,1)$, $(1 - \lambda) \mathcal{B}_2(x_0) + \lambda \mathcal{B}_2(x_1) \subset \mathcal{B}_2((1 - \lambda) x_0 + \lambda x_1)$. Indeed, if $g_0 \in \mathcal{B}_1(x_0)$ and $g_1 \in \mathcal{B}_2(x_1)$, then
\begin{align*}
	|(1 - \lambda) g_0 + \lambda g_1|
	\le (1 - \lambda) |g_0| + \lambda |g_1|
	& \le (1 - \lambda) \sqrt{h + x_0 h_o} + \lambda \sqrt{h + x_1 h_o} \\
	& \le \sqrt{h + ((1 - \lambda) x_0 + \lambda x_1) h_o} ,
\end{align*}
where the last inequality is a consequence of $\sqrt{\cdot}$ being concave.

That $G_1$ is an even function follows from $Q$ being symmetric around $0 \in \mathcal{C}[0,1]$. That $G_2$ is non-decreasing follows from $\mathcal{B}_2(x_1) \subset \mathcal{B}_2(x_2)$ for $0 \le x_1 \le x_2$.
\end{proof}

\begin{proof}[Proof of Proposition~4.14]
Note that $\UU$ and $\ZZ_{a,b}$ have pointwise expectation $0$ and are jointly Gaussian, because $\ZZ_{a,b}$ is a linear function of $\UU$. Recall that the covariance function of $\UU$ is given by $E \bigl( \UU(r) \UU(t) \bigr) = r (1 - t)$ for $0 \le r \le t \le 1$. With elementary calculations one can show that
\[
	E \bigl( \UU(t) \ZZ_{a,b}(v) \bigr) = 0
	\quad\text{for} \ t \in [0,1] \setminus (a,b) \ \text{and} \ v \in [0,1] ,
\]
and this implies stochastic independence of $(\UU(t))_{t \in [0,1] \setminus (a,b)}$ and $\ZZ_{a,b}$. Furthermore, tedious but elementary calculations reveal that
\[
	E \bigl( \ZZ_{a,b}(v) \ZZ_{a,b}(w) \bigr) = (b - a) v(1 - v)
	\quad\text{for} \ 0 \le v \le w \le 1 ,
\]
and this shows that $\ZZ_{a,b} \stackrel{d}{=} \sqrt{b - a} \, \UU$.
\end{proof}

\subsection{Proof of Theorem~3.10}
\label{subsec:BahadurSavageProof} 

By symmetry, it suffices to prove the claim about $B_n$. By monotonicity of $B_n$
\[
	P_F \Bigl( \inf_{x \in \RR} B_n (x) < \epsilon \Bigr)
	= \sup_{x \in \RR, \delta \in (0,\epsilon) } P_F ( B_n (x) < \delta ) .
\]
Hence it suffices to show that $P_F (B_n (x) < \delta ) \le (1-\epsilon)^{-n} \alpha$ for any single point $x \in \RR$ and $\delta \in (0,\epsilon)$.  To this end, consider $F_{\epsilon, \mu} := (1-\epsilon )F + \epsilon F( \cdot - \mu)$ for our given $\epsilon$ and some $\mu \in \RR$. Note that ${\cal L}_{F_{\epsilon, \mu}}(X_1,X_2,\ldots,X_n)$ describes the distribution of 
\[
	(Y_1 + Z_1 \mu, Y_2 + Z_2 \mu, \ldots, Y_n + Z_n \mu)
\]
with $2n$ independent random variables $Y_1,\ldots,Y_n \sim F$ and $Z_1, Z_2, \ldots , Z_n \sim \mathrm{Bin}(1, \epsilon)$. In particular, for any event $S_n \subset \RR^n$, 
\begin{align*}
	P_{F_{\epsilon, \mu}} \bigl( (X_1,\ldots,X_n) \in S_n \bigr)
	&=  P \bigl( (Y_1 + Z_1 \mu, \ldots, Y_n + Z_n \mu) \in S_n \bigr) \\
	&\ge  P \bigl( (Y_1,\ldots,Y_n) \in S_n , \ Z_1 = \cdots = Z_n = 0 ) \\ 
	&=  (1-\epsilon)^n P_F \bigl( (X_1,\ldots,X_n)  \in S_n \bigr) .
\end{align*}
Consequently, since $F_{\epsilon, \mu} \in {\cal F}$ too, we may conclude from 
\[
	P_{F_{\epsilon, \mu}} (A_n \le F_{\epsilon, \mu} \le B_n \ \mbox{on} \ \RR ) \ge 1-\alpha
\]
that
\begin{align*} 
	\alpha
	&\ge  P_{F_{\epsilon, \mu}} (B_n (x) < F_{\epsilon, \mu} (x) ) \\
	&\ge  (1 -\epsilon)^n P_F (B_n (x) < (1-\epsilon) F(x) + \epsilon  F(x- \mu)) \\
	&\ge (1-\epsilon)^n P_F ( B_n (x) < \epsilon F(x-\mu) ) .
\end{align*}
But for sufficiently small (negative) $\mu$, the value $\epsilon F(x-\mu)$ is greater 
than or equal to $\delta$. Then we may conclude that $\alpha \ge (1-\epsilon)^n P_F (B_n (x) < \delta)$. \hfill $\Box$

\subsection{Duality between goodness-of-fit tests and confidence bands}
\label{subsec:Duality}

\paragraph*{Continuous distribution functions}
All goodness-of-fit tests considered in this paper are of the following type. For a continuous distribution function $F$, the test statistic $T_n(F) = T_n(F,(X_i)_{i=1}^n)$ equals
\begin{equation}
\label{eq:Tn.type1}
	T_n(F) = \sup_{x \in [X_{n:1},X_{n:n-1})} \Gamma_n(\FF_n(x), F(x))
\end{equation}
or
\begin{equation}
\label{eq:Tn.type2}
	T_n(F) = \sup_{x\colon 0 < F(x) < 1} \Gamma_n(\FF_n(x), F(x))
\end{equation}
with $\Gamma_n : [0,1] \times [0,1] \to (-\infty,\infty]$ such that for any fixed $u \in [0,1]$, the function $\Gamma_n(u,\cdot)$ is continuous, decreasing on $[0,u]$ and increasing on $[u,1]$. This implies that $T_n(F)$ in \eqref{eq:Tn.type1} can be written as
\begin{equation}
\label{eq:Tn.type1b}
	T_n(F) = \max_{1 \le i < n}
		\max \bigl\{ \Gamma_n(i/n,F(X_{n:i})), \Gamma_n(i/n,F(X_{n:i+1})) \bigr\} ,
\end{equation}
while $T_n(F)$ in \eqref{eq:Tn.type2} equals
\begin{equation}
\label{eq:Tn.type2b}
	T_n(F) = \max_{1 \le i \le n}
		\max \bigl\{ \Gamma_n((i-1)/n,F(X_{n:i})), \Gamma_n(i/n,F(X_{n:i})) \bigr\} .
\end{equation}
In particular, if $F$ is the distribution function of the observations $X_i$, then $T_n(F)$ has the same distribution as
\[
	T_n = \max_{1 \le i < n}
		\max \bigl\{ \Gamma_n(i/n,\xi_{n:i}), \Gamma_n(i/n,\xi_{n:i+1}) \bigr\} ,
\]
or
\[
	T_n = \max_{1 \le i \le n}
		\max \bigl\{ \Gamma_n((i-1)/n,\xi_{n:i}), \Gamma_n(i/n,\xi_{n:i}) \bigr\} ,
\]
respectively, because $(F(X_{n:i})_{i=1}^n$ has the same distribution as $(\xi_{n:i})_{i=1}^n$.  For any critical value $\kappa \in \RR$, the inequality $T_n(F) \le \kappa$ is equivalent to
\begin{equation}
\label{eq:CBand}
	F(x) \in [a_{n,i}(\kappa), b_{n,i}(\kappa)]
	\quad\text{for} \ x \in [X_{n:i}, X_{n:i+1}) \ \text{and} \ 0 \le i \le n
\end{equation}
with certain constants $a_{n,i}(\kappa), b_{n,i}(\kappa) \in [0,1]$ such that $a_{n,0}(\kappa) = 0$ and $b_{n,n}(\kappa) = 1$. Specifically, if $T_n(F)$ is given by \eqref{eq:Tn.type1}, then $a_{n,n}(\kappa) = a_{n,n-1}(\kappa)$, $b_{n,0}(\kappa) = b_{n,1}(\kappa)$, and for $1 \le i < n$,
\begin{align*}
	a_{n,i}(\kappa) &= \min \bigl\{ t \in [0,i/n] \colon \Gamma_n(i/n,t) \le \kappa \bigr\} , \\
	b_{n,i}(\kappa) &= \max \bigl\{ t \in [i/n,1] \colon \Gamma_n(i/n,t) \le \kappa \bigr\} .
\end{align*}
If $T_n(F)$ is given by \eqref{eq:Tn.type2}, then
\begin{align*}
	a_{n,i}(\kappa) &= \min \bigl\{ t \in [0,i/n] \colon \Gamma_n(i/n,t) \le \kappa \bigr\}
		\quad\text{for} \ 1 \le i \le n , \\
	b_{n,i}(\kappa) &= \max \bigl\{ t \in [i/n,1] \colon \Gamma_n(i/n,t) \le \kappa \bigr\}
		\quad\text{for} \ 0 \le i < n .
\end{align*}
If $\Gamma_n$ satisfies the symmetry property that $\Gamma_n(u,t) = \Gamma_n(1-u,1-t)$ for all $u,t \in [0,1]$, then
\[
	a_{n,i}(\kappa) \ = \ 1 - b_{n,n-i}(\kappa)
	\quad\text{for} \ 0 \le i \le n .
\]

To compute the probability $P_F(T_n(F) \le \kappa) = P(T_n \le \kappa)$ numerically, one can use the dual representation \eqref{eq:CBand}, applied to the uniform distribution on $[0,1]$, to verify that
\begin{equation}
\label{eq:Duality}
	P(T_n \le \kappa)
	= P \bigl( a_{n,i}(\kappa) \le \xi_{n:i} \le b_{n,i-1}(\kappa)
		\ \text{for} \ 1 \le i \le n \bigr) .
\end{equation}
If for all relevant $u$, $\Gamma_n(u,t)$ is strictly decreasing on $[0,u]$ and strictly increasing on $[u,1]$, then the bounds $a_{n,i}(\kappa)$ and $b_{n,i}(\kappa)$ are continuous in $\kappa$, whence the distribution function of $T_n$ is continuous.

\paragraph*{Confidence bands for arbitrary distribution functions}
Suppose that we have chosen numbers $0 \le a_{n,i,\alpha} < b_{n,i,\alpha} \le 1$, $0 \le i \le n$, with $a_{n,0,\alpha} = 0$ and $b_{n,n,\alpha} = 1$ such that $P(a_{n,i,\alpha} \le \xi_{n:i} \le b_{n,i-1,\alpha} \ \text{for} \ 1 \le i \le n) \ge 1 - \alpha$. This leads to the confidence band $(A_{n,\alpha}, B_{n,\alpha})$ given by
\[
	\bigl[ A_{n,\alpha}(x), B_{n,\alpha}(x) \bigr]
	:= [a_{n,i,\alpha}, b_{n,i,\alpha}]
	\quad\text{for} \ x \in [X_{n:i}, X_{n:i+1}) \ \text{and} \ 0 \le i \le n .
\]
Indeed, this confidence band satisfies inequality (1.1),
\[
	P_F(A_{n,\alpha} \le F \le B_{n,\alpha} \ \text{on} \ \RR)
	\ge 1 - \alpha ,
\]
even if the underlying distribution function $F$ is not continuous. To verify this, note that $(X_{n:i})_{i=1}^n$ has the same distribution as $(F^{-1}(\xi_{n:i}))_{i=1}^n$ with $F^{-1}(u) = \min\{x \in \RR\colon F(x) \ge u\}$ for $0 < u < 1$. Moreover, $F(F^{-1}(\xi_{n:i})-) \le \xi_{n:i} \le F(F^{-1}(\xi_{n:i}))$ for $0 \le i \le n+1$. Consequently, $A_{n,\alpha} \le F \le B_{n,\alpha}$ on $\RR$ whenever $[\xi_{n:i},\xi_{n:i+1}] \subset [a_{n,i,\alpha}, b_{n,i,\alpha}]$ for $0 \le i \le n$, and the latter inclusions are equivalent to $a_{n,i,\alpha} \le \xi_{n:i} \le b_{n,i-1,\alpha}$ for $1 \le i \le n$.

\subsection{Critical values for various goodness-of-fit tests}
\label{subsec:CritValues}

Tables~\ref{tab:CritValuesDW.restr} and \ref{tab:CritValuesDW} contain $(1 - \alpha)$-quantiles of the statistics
\begin{equation}
\label{eq:Tns1DW.restr}
	T_{n,s,1} := \sup_{t \in [\xi_{n:1},\xi_{n:n-1})}
		\bigl[ n K_s(\GG_n(t),t) - C_1(\GG_n(t),t) \bigr]
\end{equation}
and
\begin{equation}
\label{eq:Tns1DW}
	T_{n,s,1} := \sup_{t \in (0,1)}
		\bigl[ n K_s(\GG_n(t),t) - C_1(\GG_n(t),t) \bigr] ,
\end{equation}
respectively, for various sample sizes $n$ and test levels $\alpha$. The parameters $s$ for the divergences $K_s$ are in $\{j/10 \colon -10 \le j \le 9\}$ and $\{j/10\colon 0 < j \le 20\}$, respectively. Thus, the critical values $\kappa_{n,s,1,\alpha}$ in the main paper are the quantiles in Table~\ref{tab:CritValuesDW.restr} for $s \le 0$ only and all quantiles in Table~\ref{tab:CritValuesDW}.

Note the big difference between the quantiles for $T_{n,s,1}$ in \eqref{eq:Tns1DW.restr} and for $T_{n,s,1}$ in \eqref{eq:Tns1DW} if $s > 0$ is small. This is not surprising, because the full supremum differs from the restricted supremum by the two terms $n K_s(0,\xi_{n:1}) \ge n \xi_{n:1}/s - C_\nu(\min\{\xi_{n:1},0.5\})$ and $n K_s(1,\xi_{n:n}) \ge n (1 - \xi_{n:n})/s - C_\nu(\max\{\xi_{n:n},0.5\})$, see the beginning of the proof of Theorem~2.1. Taking the full supremum has the advantage that the upper confidence bound for $F(x)$ is strictly smaller on $(-\infty, X_{n:1})$ than at $X_{n:1}$, just as the bound of Berk-Jones-Owen, so we might not want to always restrict the supremum.

In a similar fashion, Tables~\ref{tab:CritValuesBJ.restr} and \ref{tab:CritValuesBJ} contain $(1 - \alpha)$-quantiles of
\begin{equation}
\label{eq:TnsBJ.restr}
	T_{n,s}^{\rm BJ} := \sup_{t \in [\xi_{n:1},\xi_{n:n-1})} n K_s(\GG_n(t),t)
\end{equation}
and
\begin{equation}
\label{eq:TnsBJ}
	T_{n,s}^{\rm BJ} := \sup_{t \in (0,1)} n K_s(\GG_n(t),t) ,
\end{equation}
respectively.

Finally, Table~\ref{tab:CritValuesSP} contains critical values for the goodness-of-fit statistic
\begin{equation}
\label{eq:TnSP.supp}
	T_n^{\rm SP} = \sup_{t \in [\xi_{n:1},\xi_{n:n})}
		\frac{\sqrt{n} \bigl| \GG_n(t) - t \bigr|}
		     {\sqrt{\GG_n(1 - \GG_n)(t) h(t)}}
\end{equation}
of \cite{Stepanova_Pavlenko_2018}, where $h(t) = \log(1/[t(1-t)])$. These critical values are larger than the asymptotic ones provided by \cite{Orasch_Pouliot_2004} and used by \cite{Stepanova_Pavlenko_2018}. Table~\ref{tab:CovProbSPOP} shows that even for rather large sample sizes $n$, using the asymptotic critical values would imply too small coverage probabilities.

All these critical values and coverage probabilities have been computed numerically via the dual representation \eqref{eq:Duality} and a variant of No{\'e}'s \cite{Noe_1972} recursion; we do not rely on asymptotic theory. The critical values have been rounded up to three digits. The algorithm is essentially the same as the one of \cite{Owen_1995}, but our variant of No{\'e}'s recursion works with log-probabilities rather than probabilities. As confirmed by extensive Monte Carlo experiments, this improves numerical accuracy substantially. A description and complete computer code in R \cite{R_2019} can be found on the first author's web site https://github.com/duembgen-lutz/ConfidenceBands.

\begin{table}
\small
\[
	\arraycolsep=3pt
	\begin{array}{|l||r|r|r|r|r|r|}
	\hline
			 & \multicolumn{6}{|l|}{n} \\[-1.0ex]
		   s & 100 & 250 & 500 & 1000 & 2000 & 4000 \\
	\hline\hline
		-1.0 &  2.109 &  2.130 &  2.133 &  2.131 &  2.126 &  2.120 \\
			 &  6.718 &  6.545 &  6.372 &  6.203 &  6.051 &  5.918 \\
			 &  9.690 &  9.529 &  9.315 &  9.087 &  8.868 &  8.667 \\
			 & 18.769 & 19.009 & 18.920 & 18.745 & 18.544 & 18.343 \\
	\hline
		-0.9 &  2.066 &  2.088 &  2.092 &  2.091 &  2.087 &  2.082 \\
			 &  6.303 &  6.140 &  5.984 &  5.834 &  5.699 &  5.584 \\
			 &  8.953 &  8.773 &  8.559 &  8.338 &  8.129 &  7.941 \\
			 & 16.978 & 17.110 & 16.983 & 16.787 & 16.575 & 16.368 \\
	\hline
		-0.8 &  2.026 &  2.049 &  2.053 &  2.053 &  2.051 &  2.047 \\
			 &  5.936 &  5.788 &  5.649 &  5.517 &  5.400 &  5.300 \\
			 &  8.302 &  8.112 &  7.905 &  7.696 &  7.503 &  7.332 \\
			 & 15.404 & 15.445 & 15.284 & 15.071 & 14.850 & 14.637 \\
	\hline
		-0.7 &  1.989 &  2.012 &  2.017 &  2.018 &  2.017 &  2.014 \\
			 &  5.613 &  5.481 &  5.360 &  5.246 &  5.145 &  5.059 \\
			 &  7.729 &  7.538 &  7.344 &  7.152 &  6.978 &  6.826 \\
			 & 14.021 & 13.985 & 13.796 & 13.569 & 13.340 & 13.124 \\
	\hline
		-0.6 &  1.954 &  1.977 &  1.984 &  1.986 &  1.985 &  1.983 \\
			 &  5.329 &  5.215 &  5.111 &  5.013 &  4.927 &  4.854 \\
			 &  7.226 &  7.043 &  6.866 &  6.694 &  6.541 &  6.409 \\
			 & 12.807 & 12.708 & 12.498 & 12.260 & 12.026 & 11.808 \\
	\hline
		-0.5 &  1.921 &  1.945 &  1.953 &  1.955 &  1.955 &  1.955 \\
			 &  5.080 &  4.984 &  4.896 &  4.812 &  4.740 &  4.678 \\
			 &  6.787 &  6.619 &  6.461 &  6.311 &  6.179 &  6.066 \\
			 & 11.743 & 11.595 & 11.371 & 11.128 & 10.894 & 10.679 \\
	\hline
		-0.4 &  1.891 &  1.916 &  1.924 &  1.927 &  1.928 &  1.928 \\
			 &  4.861 &  4.783 &  4.709 &  4.639 &  4.578 &  4.526 \\
			 &  6.405 &  6.255 &  6.118 &  5.990 &  5.877 &  5.783 \\
			 & 10.814 & 10.631 & 10.401 & 10.161 &  9.934 &  9.729 \\
	\hline
		-0.3 &  1.864 &  1.888 &  1.897 &  1.901 &  1.903 &  1.904 \\
			 &  4.670 &  4.608 &  4.548 &  4.490 &  4.439 &  4.396 \\
			 &  6.075 &  5.946 &  5.829 &  5.721 &  5.627 &  5.548 \\
			 & 10.006 &  9.804 &  9.578 &  9.349 &  9.138 &  8.951 \\
	\hline
		-0.2 &  1.838 &  1.863 &  1.872 &  1.877 &  1.880 &  1.882 \\
			 &  4.503 &  4.457 &  4.408 &  4.361 &  4.320 &  4.285 \\
			 &  5.789 &  5.683 &  5.586 &  5.496 &  5.419 &  5.354 \\
			 &  9.307 &  9.101 &  8.888 &  8.679 &  8.492 &  8.329 \\
	\hline
		-0.1 &  1.815 &  1.840 &  1.849 &  1.855 &  1.859 &  1.861 \\
			 &  4.358 &  4.325 &  4.287 &  4.250 &  4.217 &  4.189 \\
			 &  5.544 &  5.460 &  5.381 &  5.308 &  5.245 &  5.193 \\
			 &  8.707 &  8.511 &  8.320 &  8.138 &  7.977 &  7.841 \\
	\hline
	\end{array}
	\quad
	\begin{array}{|l||r|r|r|r|r|r|}
	\hline
			& \multicolumn{6}{|l|}{n} \\[-1.0ex]
		  s & 100 & 250 & 500 & 1000 & 2000 & 4000 \\
	\hline\hline
		0.0 & 1.794 & 1.819 & 1.829 & 1.835 & 1.840 & 1.843 \\
			& 4.231 & 4.212 & 4.183 & 4.155 & 4.129 & 4.107 \\
			& 5.334 & 5.271 & 5.209 & 5.151 & 5.101 & 5.060 \\
			& 8.197 & 8.022 & 7.858 & 7.704 & 7.572 & 7.462 \\
	\hline
		0.1 & 1.775 & 1.800 & 1.810 & 1.817 & 1.822 & 1.827 \\
			& 4.122 & 4.114 & 4.094 & 4.073 & 4.054 & 4.038 \\
			& 5.155 & 5.111 & 5.064 & 5.020 & 4.982 & 4.950 \\
			& 7.767 & 7.620 & 7.485 & 7.362 & 7.256 & 7.169 \\
	\hline
		0.2 & 1.758 & 1.783 & 1.794 & 1.801 & 1.807 & 1.812 \\
			& 4.028 & 4.030 & 4.018 & 4.004 & 3.991 & 3.980 \\
			& 5.003 & 4.977 & 4.944 & 4.911 & 4.883 & 4.860 \\
			& 7.408 & 7.294 & 7.188 & 7.092 & 7.011 & 6.945 \\
	\hline
		0.3 & 1.744 & 1.768 & 1.779 & 1.787 & 1.794 & 1.799 \\
			& 3.949 & 3.959 & 3.953 & 3.945 & 3.937 & 3.931 \\
			& 4.876 & 4.866 & 4.844 & 4.822 & 4.802 & 4.787 \\
			& 7.112 & 7.031 & 6.953 & 6.882 & 6.822 & 6.774 \\
	\hline
		0.4 & 1.732 & 1.756 & 1.767 & 1.775 & 1.782 & 1.788 \\
			& 3.882 & 3.899 & 3.900 & 3.897 & 3.893 & 3.891 \\
			& 4.770 & 4.774 & 4.763 & 4.749 & 4.737 & 4.728 \\
			& 6.871 & 6.823 & 6.769 & 6.719 & 6.678 & 6.645 \\
	\hline
		0.5 & 1.722 & 1.745 & 1.756 & 1.765 & 1.772 & 1.799 \\
			& 3.827 & 3.851 & 3.856 & 3.857 & 3.858 & 3.858 \\
			& 4.685 & 4.700 & 4.697 & 4.691 & 4.685 & 4.681 \\
			& 6.679 & 6.659 & 6.626 & 6.595 & 6.569 & 6.549 \\
	\hline
		0.6 & 1.714 & 1.737 & 1.748 & 1.757 & 1.765 & 1.771 \\
			& 3.784 & 3.812 & 3.821 & 3.826 & 3.830 & 3.833 \\
			& 4.618 & 4.641 & 4.645 & 4.646 & 4.645 & 4.645 \\
			& 6.530 & 6.534 & 6.519 & 6.503 & 6.489 & 6.479 \\
	\hline
		0.7 & 1.710 & 1.732 & 1.742 & 1.751 & 1.759 & 1.766 \\
			& 3.753 & 3.783 & 3.795 & 3.802 & 3.809 & 3.814 \\
			& 4.568 & 4.598 & 4.607 & 4.612 & 4.616 & 4.619 \\
			& 6.420 & 6.442 & 6.441 & 6.436 & 6.432 & 6.429 \\
	\hline
 		0.8 & 1.709 & 1.729 & 1.740 & 1.748 & 1.756 & 1.763 \\
			& 3.734 & 3.765 & 3.778 & 3.787 & 3.795 & 3.802 \\
			& 4.537 & 4.569 & 4.581 & 4.589 & 4.596 & 4.602 \\
			& 6.346 & 6.380 & 6.388 & 6.392 & 6.394 & 6.397 \\
	\hline
		0.9 & 1.715 & 1.732 & 1.741 & 1.749 & 1.756 & 1.763 \\
			& 3.731 & 3.759 & 3.772 & 3.781 & 3.789 & 3.796 \\
			& 4.527 & 4.558 & 4.570 & 4.579 & 4.586 & 4.593 \\
			& 6.313 & 6.349 & 6.361 & 6.368 & 6.374 & 6.380 \\
	\hline
	\end{array}
\]
\caption{$(1 - \alpha)$-quantiles of $T_{n,s,1}$ in \eqref{eq:Tns1DW.restr} for $\alpha = 0.5, 0.1, 0.05, 0.01$.}
\label{tab:CritValuesDW.restr}
\end{table}

\begin{table}
\small
\[
	\arraycolsep=3pt
	\begin{array}{|l||r|r|r|r|r|r|}
	\hline
			& \multicolumn{6}{|l|}{n} \\[-1.0ex]
		  s & 100 & 250 & 500 & 1000 & 2000 & 4000 \\
	\hline\hline
		0.1 &  9.785 &  9.419 &  9.182 &  8.972 &  8.786 &  8.619 \\
			& 27.325 & 27.060 & 26.834 & 26.615 & 26.411 & 26.224 \\
			& 34.306 & 34.140 & 33.942 & 33.732 & 33.529 & 33.340 \\
			& 50.094 & 50.263 & 50.166 & 49.999 & 49.811 & 49.625 \\
	\hline
		0.2 &  4.136 &  3.908 &  3.770 &  3.656 &  3.560 &  3.478 \\
			& 12.692 & 12.304 & 12.038 & 11.798 & 11.584 & 11.393 \\
			& 16.253 & 15.893 & 15.630 & 15.387 & 15.168 & 14.971 \\
			& 24.302 & 24.062 & 23.826 & 23.590 & 23.367 & 23.163 \\
	\hline
		0.3 &  2.828 &  2.712 &  2.643 &  2.586 &  2.539 &  2.500 \\
			&  7.919 &  7.532 &  7.282 &  7.067 &  6.881 &  6.721 \\
			& 10.278 &  9.866 &  9.589 &  9.344 &  9.127 &  8.937 \\
			& 15.712 & 15.337 & 15.055 & 14.796 & 14.562 & 14.353 \\
	\hline
		0.4 &  2.336 &  2.266 &  2.225 &  2.193 &  2.166 &  2.144 \\
			&  5.823 &  5.543 &  5.376 &  5.239 &  5.126 &  5.033 \\
			&  7.468 &  7.108 &  6.882 &  6.693 &  6.535 &  6.401 \\
			& 11.469 & 11.044 & 10.750 & 10.491 & 10.263 & 10.063 \\
	\hline
		0.5 &  2.090 &  2.046 &  2.021 &  2.002 &  1.986 &  1.974 \\
			&  4.844 &  4.671 &  4.572 &  4.493 &  4.430 &  4.379 \\
			&  6.064 &  5.821 &  5.678 &  5.563 &  5.471 &  5.396 \\
			&  9.084 &  8.702 &  8.458 &  8.254 &  8.084 &  7.943 \\
	\hline
		0.6 &  1.951 &  1.923 &  1.908 &  1.896 &  1.888 &  1.882 \\
			&  4.347 &  4.246 &  4.188 &  4.144 &  4.109 &  4.083 \\
			&  5.349 &  5.203 &  5.121 &  5.056 &  5.006 &  4.967 \\
			&  7.750 &  7.487 &  7.331 &  7.208 &  7.110 &  7.033 \\
	\hline
		0.7 &  1.866 &  1.849 &  1.841 &  1.835 &  1.832 &  1.830 \\
			&  4.076 &  4.019 &  3.989 &  3.966 &  3.949 &  3.936 \\
			&  4.967 &  4.887 &  4.843 &  4.810 &  4.785 &  4.766 \\
			&  7.032 &  6.883 &  6.799 &  6.735 &  6.687 &  6.650 \\
	\hline
		0.8 &  1.815 &  1.805 &  1.801 &  1.800 &  1.799 &  1.800 \\
			&  3.923 &  3.895 &  3.881 &  3.871 &  3.865 &  3.862 \\
			&  4.758 &  4.720 &  4.700 &  4.685 &  4.675 &  4.669 \\
			&  6.652 &  6.583 &  6.545 &  6.518 &  6.498 &  6.484 \\
	\hline
		0.9 &  1.787 &  1.782 &  1.780 &  1.781 &  1.782 &  1.785 \\
			&  3.842 &  3.832 &  3.827 &  3.825 &  3.824 &  3.826 \\
			&  4.650 &  4.636 &  4.629 &  4.625 &  4.624 &  4.624 \\
			&  6.464 &  6.441 &  6.429 &  6.421 &  6.416 &  6.414 \\
	\hline
		1.0 &  1.780 &  1.776 &  1.776 &  1.777 &  1.779 &  1.782 \\
			&  3.824 &  3.817 &  3.815 &  3.815 &  3.816 &  3.819 \\
			&  4.624 &  4.616 &  4.613 &  4.612 &  4.613 &  4.615 \\
			&  6.415 &  6.406 &  6.401 &  6.398 &  6.397 &  6.398 \\
	\hline
	\end{array}
	\ \
	\begin{array}{|l||r|r|r|r|r|r|}
	\hline
			& \multicolumn{6}{|l|}{n} \\[-1.0ex]
		  s & 100 & 250 & 500 & 1000 & 2000 & 4000 \\
	\hline\hline
		1.1 &  1.787 &  1.785 &  1.785 &  1.786 &  1.789 &  1.791 \\
			&  3.872 &  3.861 &  3.856 &  3.852 &  3.851 &  3.851 \\
			&  4.700 &  4.683 &  4.673 &  4.667 &  4.664 &  4.662 \\
			&  6.594 &  6.556 &  6.534 &  6.519 &  6.507 &  6.500 \\
	\hline
		1.2 &  1.805 &  1.804 &  1.804 &  1.805 &  1.807 &  1.810 \\
			&  3.984 &  3.963 &  3.950 &  3.941 &  3.935 &  3.931 \\
			&  4.888 &  4.850 &  4.828 &  4.811 &  4.798 &  4.789 \\
			&  7.160 &  7.050 &  6.987 &  6.938 &  6.899 &  6.869 \\
	\hline
		1.3 &  1.831 &  1.831 &  1.831 &  1.832 &  1.834 &  1.836 \\
			&  4.157 &  4.120 &  4.098 &  4.081 &  4.068 &  4.058 \\
			&  5.202 &  5.131 &  5.090 &  5.057 &  5.031 &  5.010 \\
			&  8.398 &  8.161 &  8.023 &  7.912 &  7.821 &  7.746 \\
	\hline
		1.4 &  1.863 &  1.864 &  1.864 &  1.865 &  1.866 &  1.867 \\
			&  4.396 &  4.338 &  4.303 &  4.275 &  4.253 &  4.235 \\
			&  5.675 &  5.556 &  5.487 &  5.431 &  5.386 &  5.350 \\
			& 10.901 & 10.534 & 10.306 & 10.113 &  9.946 &  9.802 \\
	\hline
		1.5 &  1.901 &  1.903 &  1.903 &  1.903 &  1.903 &  1.904 \\
			&  4.711 &  4.625 &  4.574 &  4.532 &  4.497 &  4.469 \\
			&  6.376 &  6.189 &  6.079 &  5.991 &  5.918 &  5.859 \\
			& 15.201 & 14.812 & 14.566 & 14.352 & 14.163 & 13.993 \\
	\hline
		1.6 &  1.944 &  1.946 &  1.946 &  1.945 &  1.945 &  1.945 \\
			&  5.127 &  5.002 &  4.928 &  4.867 &  4.817 &  4.776 \\
			&  7.427 &  7.153 &  6.988 &  6.853 &  6.741 &  6.647 \\
			& 21.701 & 21.319 & 21.076 & 20.865 & 20.677 & 20.509 \\
	\hline
		1.7 &  1.992 &  1.994 &  1.993 &  1.992 &  1.990 &  1.990 \\
			&  5.678 &  5.502 &  5.397 &  5.311 &  5.241 &  5.182 \\
			&  9.001 &  8.646 &  8.424 &  8.236 &  8.075 &  7.937 \\
			& 31.292 & 30.914 & 30.674 & 30.464 & 30.278 & 30.111 \\
	\hline
		1.8 &  2.044 &  2.045 &  2.044 &  2.042 &  2.040 &  2.038 \\
			&  6.420 &  6.180 &  6.035 &  5.916 &  5.817 &  5.734 \\
			& 11.255 & 10.864 & 10.614 & 10.397 & 10.206 & 10.038 \\
			& 45.476 & 45.101 & 44.862 & 44.654 & 44.468 & 44.302 \\
	\hline
 		1.9 &  2.100 &  2.101 &  2.099 &  2.096 &  2.093 &  2.090 \\
			&  7.426 &  7.118 &  6.926 &  6.766 &  6.631 &  6.517 \\
			& 14.323 & 13.929 & 13.677 & 13.458 & 13.263 & 13.090 \\
			& 66.584 & 66.212 & 65.974 & 65.766 & 65.582 & 65.416 \\
	\hline
		2.0 &  2.160 &  2.161 &  2.158 &  2.154 &  2.150 &  2.146 \\
			&  8.777 &  8.414 &  8.182 &  7.983 &  7.811 &  7.662 \\
			& 18.383 & 17.995 & 17.747 & 17.530 & 17.338 & 17.167 \\
			& 98.206 & 97.837 & 97.600 & 97.393 & 97.209 & 97.044 \\
	\hline
	\end{array}
\]
\caption{$(1 - \alpha)$-quantiles of $T_{n,s,1}$ in \eqref{eq:Tns1DW} for $\alpha = 0.5, 0.1, 0.05, 0.01$.}
\label{tab:CritValuesDW}
\end{table}

\begin{table}
\small
\[
	\arraycolsep=3pt
	\begin{array}{|l||r|r|r|r|r|r|}
	\hline
			 & \multicolumn{6}{|l|}{n} \\[-1.0ex]
		   s & 100 & 250 & 500 & 1000 & 2000 & 4000 \\
	\hline\hline
		-1.0 &  2.800 &  3.037 &  3.186 &  3.316 &  3.431 &  3.534 \\
			 &  7.955 &  8.258 &  8.379 &  8.452 &  8.497 &  8.529 \\
			 & 11.012 & 11.420 & 11.567 & 11.643 & 11.684 & 11.706 \\
			 & 20.055 & 20.927 & 21.229 & 21.382 & 21.459 & 21.497 \\
	\hline
		-0.9 &  2.744 &  2.980 &  3.129 &  3.259 &  3.375 &  3.479 \\
			 &  7.492 &  7.783 &  7.902 &  7.977 &  8.027 &  8.063 \\
			 & 10.242 & 10.617 & 10.753 & 10.826 & 10.866 & 10.889 \\
			 & 18.251 & 19.015 & 19.279 & 19.413 & 19.481 & 19.515 \\
	\hline
		-0.8 &  2.692 &  2.927 &  3.076 &  3.206 &  3.322 &  3.427 \\
			 &  7.075 &  7.355 &  7.474 &  7.552 &  7.607 &  7.649 \\
			 &  9.552 &  9.899 & 10.028 & 10.099 & 10.140 & 10.165 \\
			 & 16.661 & 17.332 & 17.565 & 17.683 & 17.742 & 17.772 \\
	\hline
		-0.7 &  2.645 &  2.878 &  3.026 &  3.157 &  3.273 &  3.378 \\
			 &  6.699 &  6.971 &  7.092 &  7.174 &  7.235 &  7.283 \\
			 &  8.934 &  9.258 &  9.382 &  9.453 &  9.496 &  9.524 \\
			 & 15.259 & 15.851 & 16.056 & 16.161 & 16.214 & 16.241 \\
	\hline
		-0.6 &  2.600 &  2.832 &  2.980 &  3.111 &  3.227 &  3.332 \\
			 &  6.361 &  6.627 &  6.750 &  6.837 &  6.904 &  6.958 \\
			 &  8.382 &  8.687 &  8.809 &  8.881 &  8.928 &  8.960 \\
			 & 14.022 & 14.546 & 14.728 & 14.822 & 14.870 & 14.894 \\
	\hline
		-0.5 &  2.560 &  2.790 &  2.938 &  3.068 &  3.184 &  3.290 \\
			 &  6.057 &  6.318 &  6.445 &  6.538 &  6.611 &  6.672 \\
			 &  7.889 &  8.180 &  8.300 &  8.376 &  8.428 &  8.466 \\
			 & 12.930 & 13.396 & 13.560 & 13.645 & 13.689 & 13.712 \\
	\hline
		-0.4 &  2.523 &  2.752 &  2.898 &  3.028 &  3.145 &  3.250 \\
			 &  5.784 &  6.042 &  6.173 &  6.272 &  6.351 &  6.419 \\
			 &  7.449 &  7.728 &  7.850 &  7.931 &  7.988 &  8.033 \\
			 & 11.966 & 12.380 & 12.533 & 12.611 & 12.653 & 12.676 \\
	\hline
		-0.3 &  2.489 &  2.716 &  2.862 &  2.992 &  3.108 &  3.213 \\
			 &  5.540 &  5.798 &  5.932 &  6.036 &  6.122 &  6.196 \\
			 &  7.058 &  7.330 &  7.454 &  7.540 &  7.605 &  7.657 \\
			 & 11.114 & 11.494 & 11.632 & 11.706 & 11.748 & 11.772 \\
	\hline
		-0.2 &  2.459 &  2.684 &  2.829 &  2.958 &  3.074 &  3.179 \\
			 &  5.323 &  5.580 &  5.718 &  5.827 &  5.919 &  5.999 \\
			 &  6.710 &  6.978 &  7.105 &  7.198 &  7.270 &  7.330 \\
			 & 10.366 & 10.713 & 10.843 & 10.917 & 10.960 & 10.988 \\
	\hline
		-0.1 &  2.432 &  2.655 &  2.799 &  2.928 &  3.043 &  3.148 \\
			 &  5.129 &  5.386 &  5.528 &  5.643 &  5.741 &  5.826 \\
			 &  6.403 &  6.668 &  6.800 &  6.899 &  6.979 &  7.047 \\
			 &  9.708 & 10.030 & 10.156 & 10.232 & 10.280 & 10.313 \\
	\hline
	\end{array}
	\quad
	\begin{array}{|l||r|r|r|r|r|r|}
	\hline
			& \multicolumn{6}{|l|}{n} \\[-1.0ex]
		  s & 100 & 250 & 500 & 1000 & 2000 & 4000 \\
	\hline\hline
		0.0 & 2.408 & 2.629 & 2.772 & 2.900 & 3.015 & 3.120 \\
			& 4.958 & 5.216 & 5.362 & 5.481 & 5.584 & 5.674 \\
			& 6.133 & 6.397 & 6.533 & 6.640 & 6.727 & 6.803 \\
			& 9.132 & 9.436 & 9.561 & 9.641 & 9.696 & 9.737 \\
	\hline
		0.1 & 2.387 & 2.606 & 2.749 & 2.876 & 2.990 & 3.095 \\
			& 4.808 & 5.066 & 5.215 & 5.340 & 5.447 & 5.542 \\
			& 5.897 & 6.160 & 6.302 & 6.416 & 6.510 & 6.593 \\
			& 8.631 & 8.922 & 9.050 & 9.136 & 9.200 & 9.250 \\
	\hline
		0.2 & 2.370 & 2.587 & 2.728 & 2.855 & 2.969 & 3.073 \\
			& 4.678 & 4.936 & 5.088 & 5.216 & 5.328 & 5.427 \\
			& 5.692 & 5.956 & 6.103 & 6.223 & 6.324 & 6.414 \\
			& 8.197 & 8.481 & 8.613 & 8.709 & 8.783 & 8.844 \\
	\hline
		0.3 & 2.358 & 2.571 & 2.711 & 2.837 & 2.950 & 3.054 \\
			& 4.566 & 4.825 & 4.979 & 5.111 & 5.225 & 5.328 \\
			& 5.517 & 5.782 & 5.933 & 6.059 & 6.166 & 6.261 \\
			& 7.826 & 8.106 & 8.245 & 8.350 & 8.435 & 8.507 \\
	\hline
		0.4 & 2.349 & 2.560 & 2.699 & 2.823 & 2.936 & 3.038 \\
			& 4.473 & 4.731 & 4.887 & 5.021 & 5.139 & 5.244 \\
			& 5.369 & 5.635 & 5.790 & 5.921 & 6.034 & 6.134 \\
			& 7.513 & 7.792 & 7.939 & 8.054 & 8.150 & 8.232 \\
	\hline
		0.5 & 2.345 & 2.553 & 2.690 & 2.813 & 2.925 & 3.027 \\
			& 4.399 & 4.654 & 4.812 & 4.948 & 5.068 & 5.175 \\
			& 5.249 & 5.514 & 5.673 & 5.808 & 5.925 & 6.030 \\
			& 7.255 & 7.535 & 7.689 & 7.814 & 7.919 & 8.012 \\
	\hline
		0.6 & 2.344 & 2.550 & 2.686 & 2.808 & 2.919 & 3.021 \\
			& 4.343 & 4.596 & 4.754 & 4.891 & 5.012 & 5.121 \\
			& 5.157 & 5.420 & 5.581 & 5.718 & 5.839 & 5.947 \\
			& 7.050 & 7.330 & 7.491 & 7.624 & 7.738 & 7.839 \\
	\hline
		0.7 & 2.357 & 2.557 & 2.689 & 2.809 & 2.918 & 3.019 \\
			& 4.309 & 4.558 & 4.714 & 4.851 & 4.973 & 5.083 \\
			& 5.094 & 5.353 & 5.514 & 5.653 & 5.776 & 5.887 \\
			& 6.899 & 7.178 & 7.343 & 7.482 & 7.604 & 7.711 \\
	\hline
		0.8 & 2.380 & 2.576 & 2.705 & 2.823 & 2.930 & 3.028 \\
			& 4.301 & 4.543 & 4.697 & 4.833 & 4.954 & 5.063 \\
			& 5.066 & 5.319 & 5.478 & 5.617 & 5.740 & 5.851 \\
			& 6.806 & 7.080 & 7.247 & 7.390 & 7.515 & 7.627 \\
	\hline
		0.9 & 2.408 & 2.602 & 2.730 & 2.846 & 2.952 & 3.049 \\
			& 4.339 & 4.568 & 4.715 & 4.846 & 4.964 & 5.071 \\
			& 5.090 & 5.330 & 5.483 & 5.618 & 5.739 & 5.848 \\
			& 6.792 & 7.052 & 7.214 & 7.356 & 7.481 & 7.593 \\
	\hline
	\end{array}
\]
\caption{$(1 - \alpha)$-quantiles of $T_{n,s}^{\rm BJ}$ in \eqref{eq:TnsBJ.restr} for $\alpha = 0.5, 0.1, 0.05, 0.01$.}
\label{tab:CritValuesBJ.restr}
\end{table}

\begin{table}
\small
\[
	\arraycolsep=3pt
	\begin{array}{|l||r|r|r|r|r|r|}
	\hline
			& \multicolumn{6}{|l|}{n} \\[-1.0ex]
		  s & 100 & 250 & 500 & 1000 & 2000 & 4000 \\
	\hline\hline
		0.1 & 12.248 & 12.271 & 12.279 & 12.283 & 12.285 & 12.286 \\
			& 29.327 & 29.549 & 29.623 & 29.661 & 29.679 & 29.689 \\
			& 36.177 & 36.527 & 36.644 & 36.709 & 36.732 & 36.747 \\
			& 51.722 & 52.460 & 52.708 & 52.708 & 52.896 & 52.927 \\
	\hline
		0.2 &  6.235 &  6.273 &  6.296 &  6.316 &  6.335 &  6.353 \\
			& 14.689 & 14.788 & 14.822 & 14.838 & 14.847 & 14.851 \\
			& 18.123 & 18.279 & 18.331 & 18.357 & 18.370 & 18.377 \\
			& 25.930 & 26.258 & 26.369 & 26.424 & 26.452 & 26.466 \\
	\hline
		0.3 &  4.424 &  4.506 &  4.563 &  4.615 &  4.663 &  4.709 \\
			&  9.845 &  9.911 &  9.936 &  9.950 &  9.959 &  9.965 \\
			& 12.122 & 12.217 & 12.249 & 12.266 & 12.275 & 12.280 \\
			& 17.336 & 17.529 & 17.594 & 17.627 & 17.643 & 17.651 \\
	\hline
		0.4 &  3.633 &  3.747 &  3.825 &  3.897 &  3.963 &  4.026 \\
			&  7.511 &  7.584 &  7.621 &  7.649 &  7.674 &  7.695 \\
			&  9.181 &  9.260 &  9.292 &  9.313 &  9.328 &  9.341 \\
			& 13.063 & 13.193 & 13.238 & 13.261 & 13.274 & 13.280 \\
	\hline
		0.5 &  3.211 &  3.344 &  3.434 &  3.517 &  3.594 &  3.666 \\
			&  6.227 &  6.330 &  6.391 &  6.444 &  6.492 &  6.537 \\
			&  7.518 &  7.613 &  7.663 &  7.704 &  7.739 &  7.772 \\
			& 10.560 & 10.668 & 10.711 & 10.737 & 10.756 & 10.770 \\
	\hline
		0.6 &  2.959 &  3.103 &  3.201 &  3.291 &  3.374 &  3.452 \\
			&  5.477 &  5.611 &  5.697 &  5.773 &  5.842 &  5.906 \\
			&  6.527 &  6.653 &  6.729 &  6.796 &  6.856 &  6.911 \\
			&  8.999 &  9.117 &  9.177 &  9.223 &  9.262 &  9.297 \\
	\hline
		0.7 &  2.800 &  2.951 &  3.054 &  3.148 &  3.236 &  3.317 \\
			&  5.024 &  5.184 &  5.288 &  5.381 &  5.467 &  5.545 \\
			&  5.926 &  6.083 &  6.183 &  6.271 &  6.352 &  6.425 \\
			&  8.021 &  8.170 &  8.257 &  8.332 &  8.399 &  8.459 \\
	\hline
		0.8 &  2.708 &  2.862 &  2.966 &  3.063 &  3.153 &  3.236 \\
			&  4.754 &  4.931 &  5.048 &  5.153 &  5.248 &  5.336 \\
			&  5.568 &  5.748 &  5.865 &  5.969 &  6.064 &  6.151 \\
			&  7.431 &  7.613 &  7.727 &  7.828 &  7.918 &  8.000 \\
	\hline
		0.9 &  2.656 &  2.813 &  2.921 &  3.019 &  3.111 &  3.195 \\
			&  4.618 &  4.803 &  4.925 &  5.036 &  5.137 &  5.230 \\
			&  5.384 &  5.576 &  5.702 &  5.815 &  5.918 &  6.012 \\
			&  7.120 &  7.324 &  7.455 &  7.571 &  7.676 &  7.771 \\
	\hline
		1.0 &  2.629 &  2.791 &  2.901 &  3.002 &  3.095 &  3.181 \\
			&  4.609 &  4.793 &  4.916 &  5.027 &  5.129 &  5.222 \\
			&  5.377 &  5.566 &  5.691 &  5.804 &  5.907 &  6.001 \\
			&  7.103 &  7.300 &  7.429 &  7.545 &  7.650 &  7.746 \\
	\hline
	\end{array}
	\ \
	\begin{array}{|l||r|r|r|r|r|r|}
	\hline
			& \multicolumn{6}{|l|}{n} \\[-1.0ex]
		  s & 100 & 250 & 500 & 1000 & 2000 & 4000 \\
	\hline\hline
		1.1 &  2.620 &  2.786 &  2.898 &  3.002 &  3.097 &  3.185 \\
			&  4.698 &  4.879 &  5.000 &  5.109 &  5.209 &  5.301 \\
			&  5.536 &  5.715 &  5.834 &  5.941 &  6.039 &  6.129 \\
			&  7.518 &  7.677 &  7.785 &  7.882 &  7.971 &  8.053 \\
	\hline
		1.2 &  2.623 &  2.794 &  2.909 &  3.015 &  3.112 &  3.201 \\
			&  4.871 &  5.046 &  5.163 &  5.268 &  5.364 &  5.453 \\
			&  5.850 &  6.012 &  6.120 &  6.218 &  6.308 &  6.390 \\
			&  8.492 &  8.593 &  8.663 &  8.727 &  8.787 &  8.843 \\
	\hline
		1.3 &  2.638 &  2.812 &  2.930 &  3.038 &  3.137 &  3.228 \\
			&  5.125 &  5.291 &  5.401 &  5.500 &  5.591 &  5.674 \\
			&  6.336 &  6.473 &  6.565 &  6.649 &  6.726 &  6.797 \\
			& 10.284 & 10.323 & 10.349 & 10.374 & 10.397 & 10.419 \\
	\hline
		1.4 &  2.661 &  2.840 &  2.960 &  3.070 &  3.171 &  3.263 \\
			&  5.470 &  5.621 &  5.721 &  5.812 &  5.895 &  5.971 \\
			&  7.035 &  7.140 &  7.211 &  7.275 &  7.335 &  7.390 \\
			& 13.224 & 13.232 & 13.237 & 13.240 & 13.243 & 13.246 \\
	\hline
		1.5 &  2.693 &  2.875 &  2.998 &  3.110 &  3.212 &  3.305 \\
			&  5.922 &  6.053 &  6.140 &  6.219 &  6.291 &  6.358 \\
			&  8.016 &  8.086 &  8.132 &  8.173 &  8.211 &  8.247 \\
			& 17.669 & 17.672 & 17.673 & 17.674 & 17.674 & 17.674 \\
	\hline
		1.6 &  2.732 &  2.918 &  3.042 &  3.156 &  3.259 &  3.354 \\
			&  6.505 &  6.613 &  6.684 &  6.748 &  6.806 &  6.861 \\
			&  9.364 &  9.404 &  9.428 &  9.449 &  9.467 &  9.485 \\
			& 24.209 & 24.211 & 24.212 & 24.213 & 24.213 & 24.213 \\
	\hline
		1.7 &  2.778 &  2.967 &  3.093 &  3.208 &  3.312 &  3.408 \\
			&  7.253 &  7.336 &  7.388 &  7.435 &  7.478 &  7.518 \\
			& 11.178 & 11.200 & 11.211 & 11.219 & 11.226 & 11.231 \\
			& 33.817 & 33.820 & 33.821 & 33.822 & 33.822 & 33.822 \\
	\hline
		1.8 &  2.831 &  3.022 &  3.150 &  3.266 &  3.371 &  3.467 \\
			&  8.205 &  8.265 &  8.300 &  8.330 &  8.357 &  8.382 \\
			& 13.576 & 13.591 & 13.597 & 13.600 & 13.602 & 13.604 \\
			& 48.012 & 48.016 & 48.017 & 48.018 & 48.018 & 48.018 \\
	\hline
		1.9 &  2.891 &  3.084 &  3.214 &  3.330 &  3.436 &  3.533 \\
			&  9.407 &  9.449 &  9.471 &  9.488 &  9.503 &  9.515 \\
			& 16.716 & 16.729 & 16.733 & 16.735 & 16.736 & 16.737 \\
			& 69.125 & 69.131 & 69.133 & 69.134 & 69.134 & 69.135 \\
	\hline
		2.0 &  2.958 &  3.153 &  3.283 &  3.400 &  3.506 &  3.603 \\
			& 10.914 & 10.945 & 10.959 & 10.968 & 10.975 & 10.980 \\
			& 20.815 & 20.827 & 20.831 & 20.833 & 20.834 & 20.835 \\
			& 100.76 & 100.76 & 100.77 & 100.77 & 100.77 & 100.77 \\
	\hline
	\end{array}
\]
\caption{$(1 - \alpha)$-quantiles of $T_{n,s}^{\rm BJ}$ in \eqref{eq:TnsBJ} for $\alpha = 0.5, 0.1, 0.05, 0.01$.}
\label{tab:CritValuesBJ}
\end{table}

\begin{table}
\[
	\begin{array}{|r|r|r|r|r|r|r|}
	\hline
		\multicolumn{7}{|l|}{n} \\[-1.0ex]
		100 & 250 & 500 & 1000 & 2000 & 4000 & 8000 \\
	\hline
		2.892 & 2.914 & 2.919 & 2.919 & 2.916 & 2.912 & 2.907 \\
		4.286 & 4.282 & 4.270 & 4.256 & 4.244 & 4.233 & 4.224 \\
		4.768 & 4.758 & 4.742 & 4.726 & 4.712 & 4.701 & 4.691 \\
		5.780 & 5.754 & 5.728 & 5.704 & 5.684 & 5.668 & 5.655 \\
	\hline
	\end{array}
\]
\caption{$(1 - \alpha)$-quantiles of $T_n^{\rm SP}$ in \eqref{eq:TnSP.supp} for $\alpha = 0.5, 0.1, 0.05, 0.01$.}
\label{tab:CritValuesSP}
\end{table}

\begin{table}
\[
	\begin{array}{|l||r|r|r|r|r|r|r|r|}
	\hline
			& \multicolumn{8}{|l|}{n} \\[-1.0ex]
		\kappa & 100    & 250    & 500    & 1000   & 2000   & 4000   & 8000   & \infty^* \\
	\hline
		2.80   & 0.4586 & 0.4473 & 0.4438 & 0.4428 & 0.4433 & 0.4446 & 0.4464 & 0.50 \\
		4.12   & 0.8748 & 0.8751 & 0.8770 & 0.8792 & 0.8811 & 0.8829 & 0.8843 & 0.90 \\
		4.57   & 0.9331 & 0.9339 & 0.9353 & 0.9367 & 0.9380 & 0.9390 & 0.9399 & 0.95 \\
		5.53   & 0.9849 & 0.9855 & 0.9860 & 0.9865 & 0.9869 & 0.9873 & 0.9875 & 0.99 \\
	\hline
	\end{array}
\]
\caption{True coverage probabilities of the confidence bands of \cite{Stepanova_Pavlenko_2018} with the quantiles of \cite{Orasch_Pouliot_2004}, rounded to four digits. $^*$Intended limits.}
\label{tab:CovProbSPOP}
\end{table}

\subsection{Additional numerical examples}
\label{subsec:Additional.examples}

\begin{figure}
\includegraphics[width=0.9\textwidth]{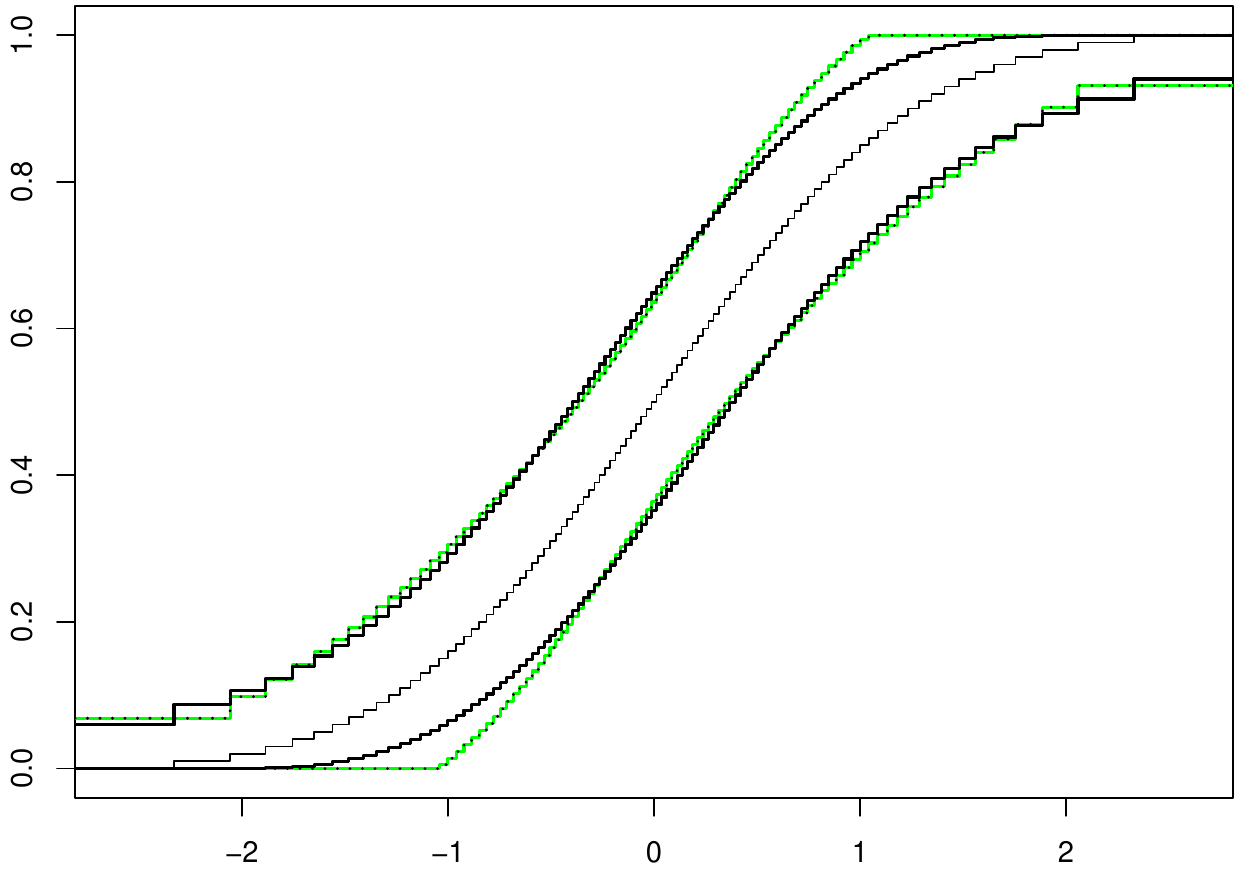}

\includegraphics[width=0.9\textwidth]{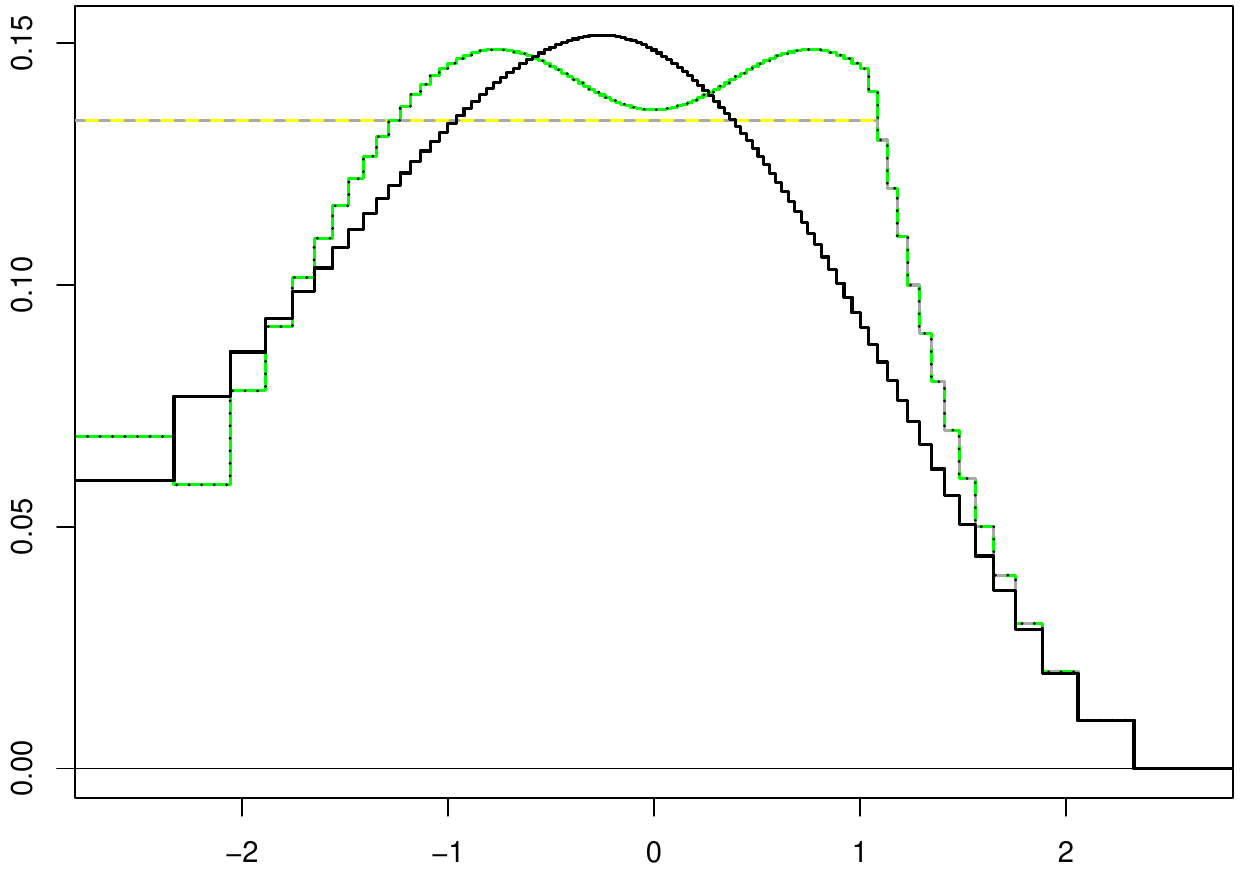}
\caption{$95\%$-confidence bands for $n = 100$. Upper panel: $(A_{n,1,1,\alpha}, B_{n,1,1,\alpha})$ (solid) and $(A_{n,\alpha}^{\rm SP}, B_{n,\alpha}^{\rm SP})$ (green, dotted). Lower panel: centered upper bounds $B_{n,1,1,\alpha} - \FF_{n}$ (solid), $B_{n,\alpha}^{\rm SP} - \FF_{n}$ (green, dotted) and $B_{n,\alpha}^{\rm KS} - \FF_{n}$ (dashed).}
\label{fig:AsEfficiency.1add}
\end{figure}

\begin{figure}
\includegraphics[width=0.9\textwidth]{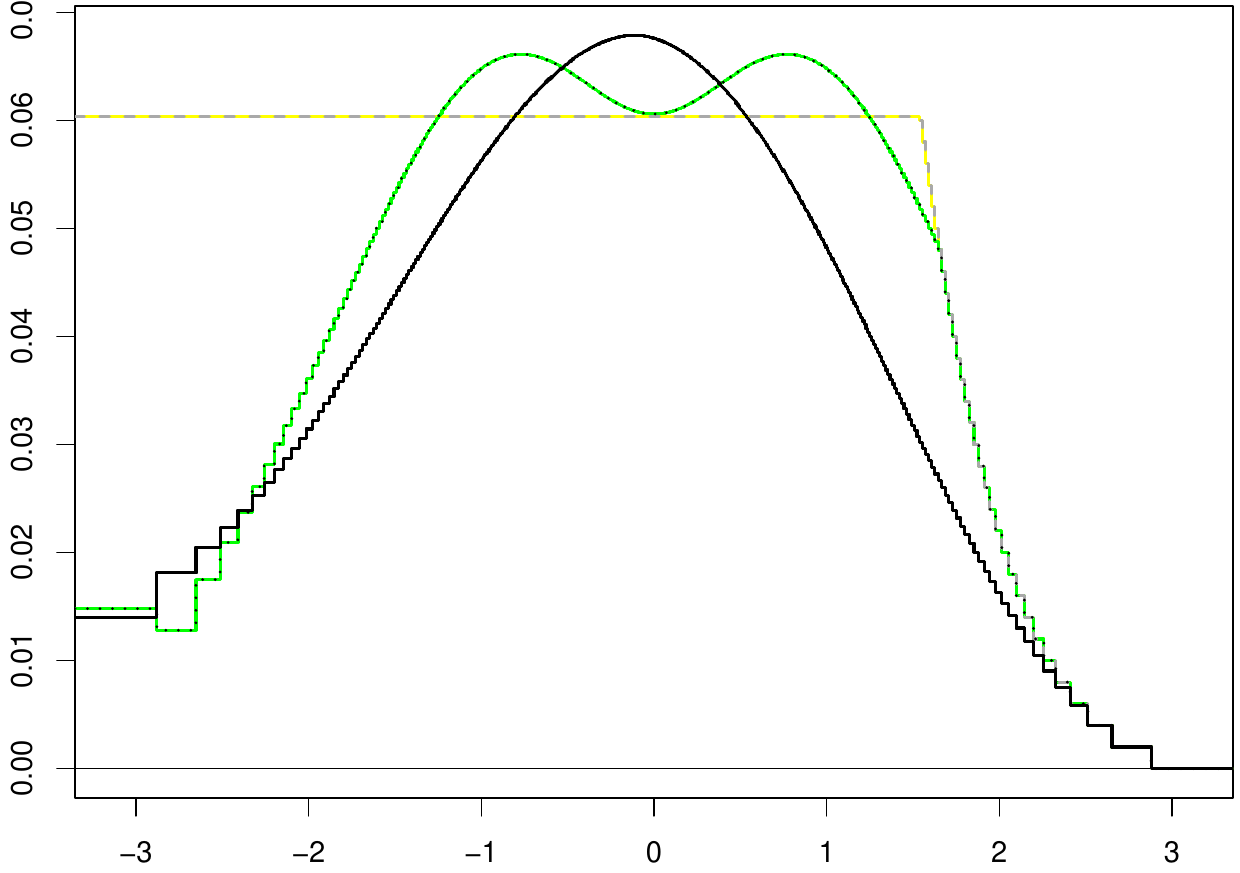}

\includegraphics[width=0.9\textwidth]{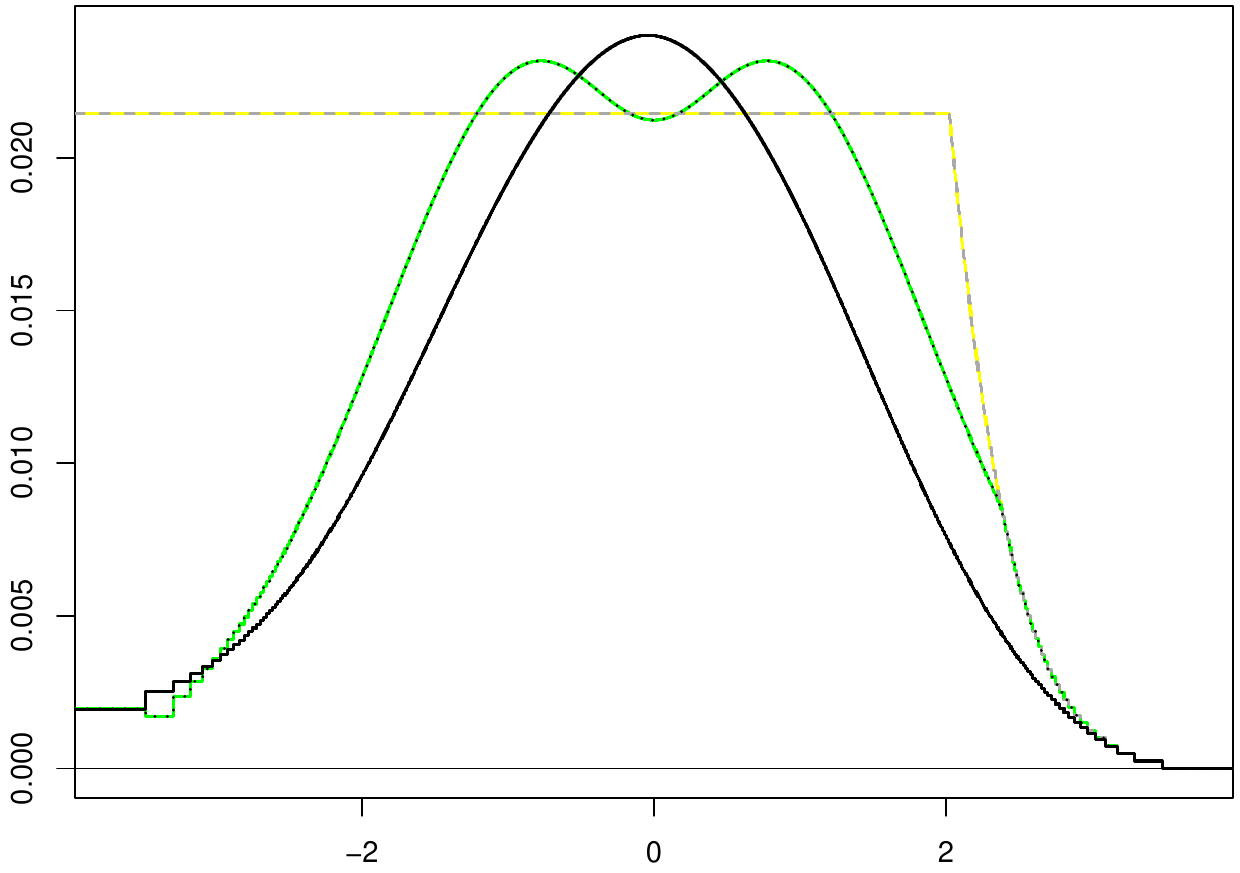}
\caption{Centered upper $95\%$-confidence bounds $B_{n,1,1,\alpha} - \FF_{n}$ (solid), $B_{n,\alpha}^{\rm SP} - \FF_{n}$ (green, dotted) and $B_{n,\alpha}^{\rm KS} - \FF_{n}$ (yellow, dashed) for $n = 500$ (upper panel) and $n = 4000$ (lower panel).}
\label{fig:AsEfficiency.2add}
\end{figure}

In Example~3.10, we compared the new $95\%$-confidence bands $(A_{n,1,1,\alpha}, B_{n,1,1,\alpha})$ with the confidence bands $(A_{n,\alpha}^{\rm KS}, B_{n,\alpha}^{\rm KS})$ and $(A_{n,1,\alpha}^{\rm BJO}, B_{n,1,\alpha}^{\rm BJO})$. In Figures~\ref{fig:AsEfficiency.1add} and \ref{fig:AsEfficiency.2add}, we compare the new bands with the $95\%$-confidence bands $(A_{n,\alpha}^{\rm SP}, B_{n,\alpha}^{\rm SP})$ of Stepanova and Pavlenko \cite{Stepanova_Pavlenko_2018}. The latter have been computed with the nonasymptotic critical values in Section~\ref{subsec:CritValues}. As predicted by our Remark~3.8, the band $(A_{n,\alpha}^{\rm SP}, B_{n,\alpha}^{\rm SP})$ is wider than $(A_{n,1,1,\alpha}, B_{n,1,1,\alpha})$ in the boundary regions, except for a rather small region in the left (resp.\ right) tail where $B_{n,\alpha}^{\rm SP} < B_{n,1,1,\alpha}$ (resp.\ $A_{n,\alpha}^{\rm SP} > A_{n,1,1,\alpha}$). An explanation for this is the fact that the test statistic $T_n^{\rm SP}$ corresponds to the divergences $K_s(\cdot,\cdot)$ with $s = -1$, see also Remark~3.6.

\begin{figure}
\includegraphics[width=0.9\textwidth]{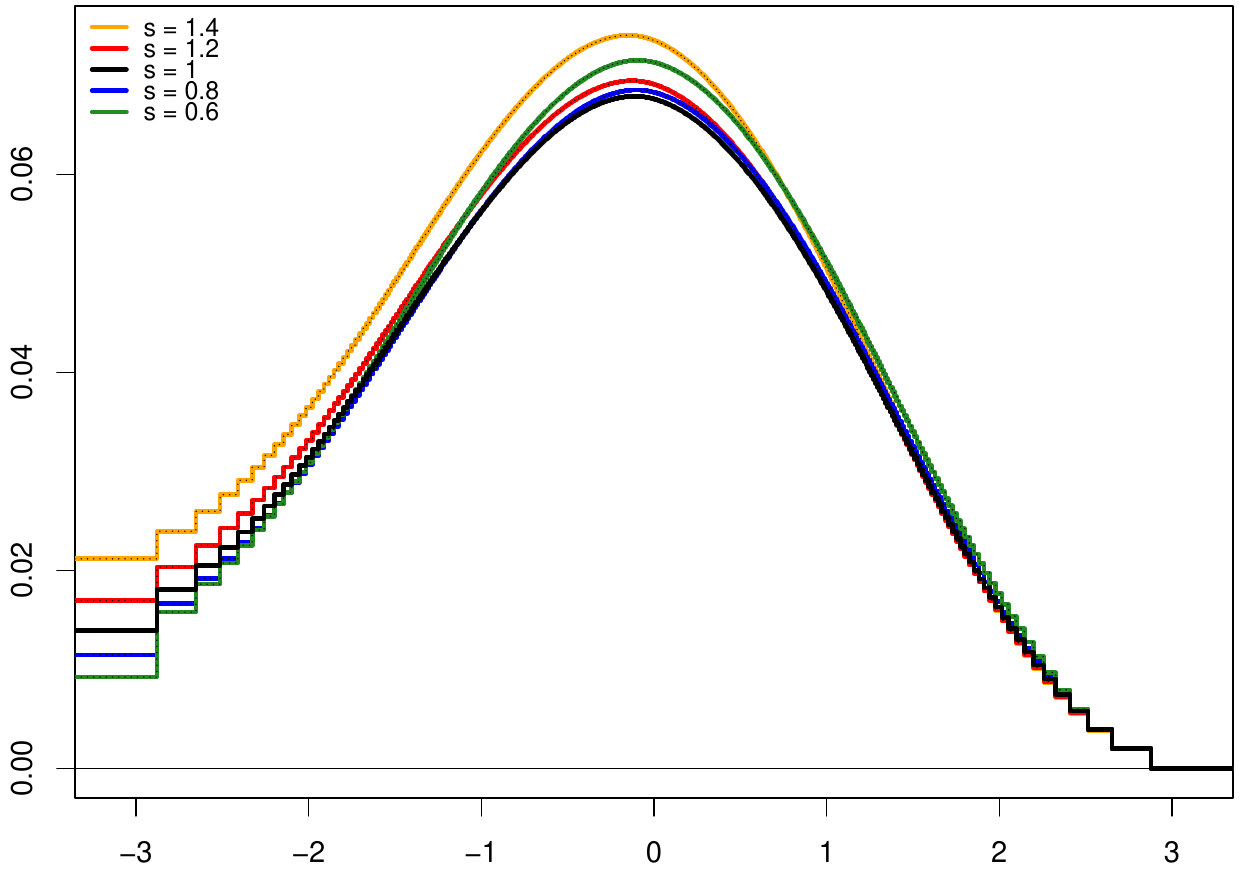}

\includegraphics[width=0.9\textwidth]{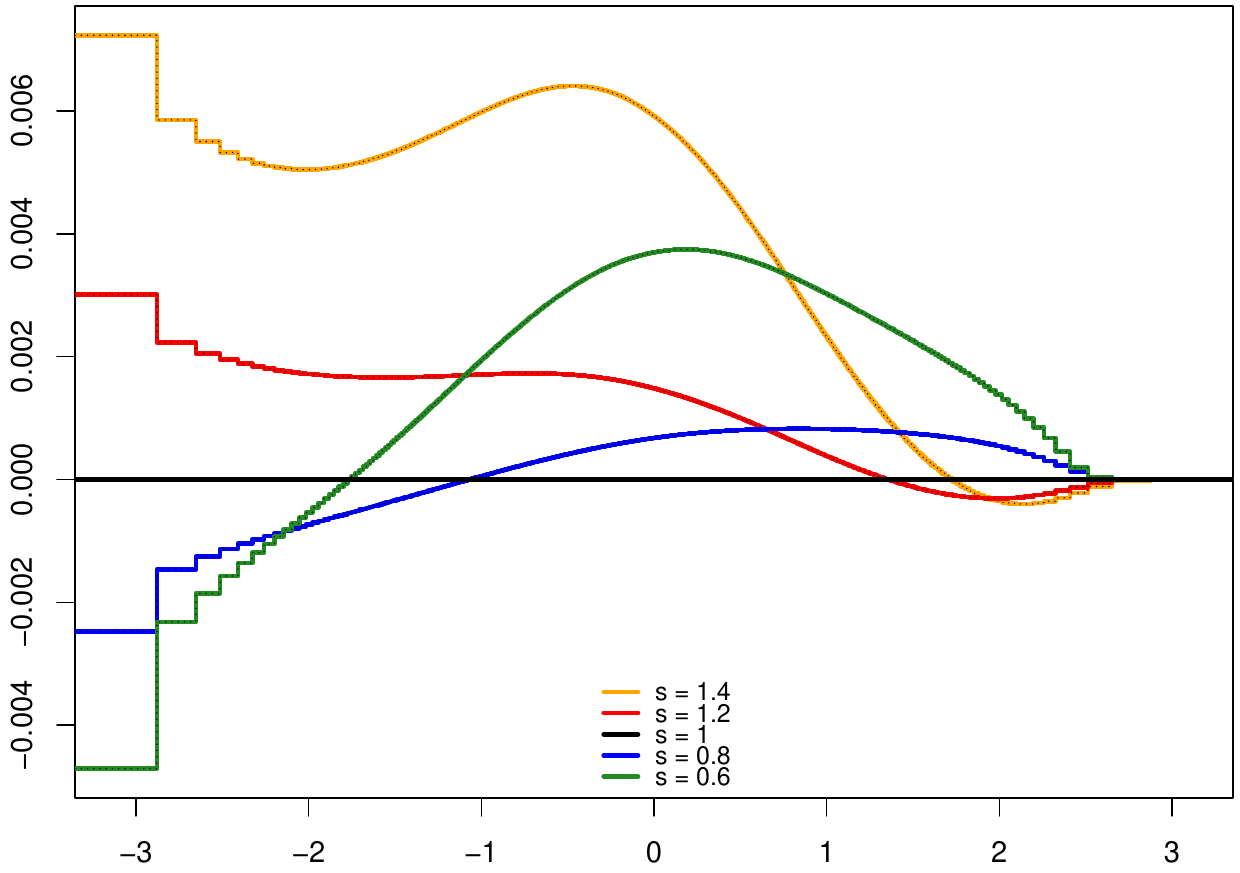}
\caption{Upper $95\%$-confidence bounds for $n = 500$ and $s \in \{0.6,0.8,1,1.2,1.4\}$. Upper panel: centered bounds $B_{n,s,1,\alpha} - \FF_n$. Lower panel: differences $B_{n,s,1,\alpha} - B_{n,1,1,\alpha}$.}
\label{fig:UBSadd500}
\end{figure}

\begin{figure}
\includegraphics[width=0.9\textwidth]{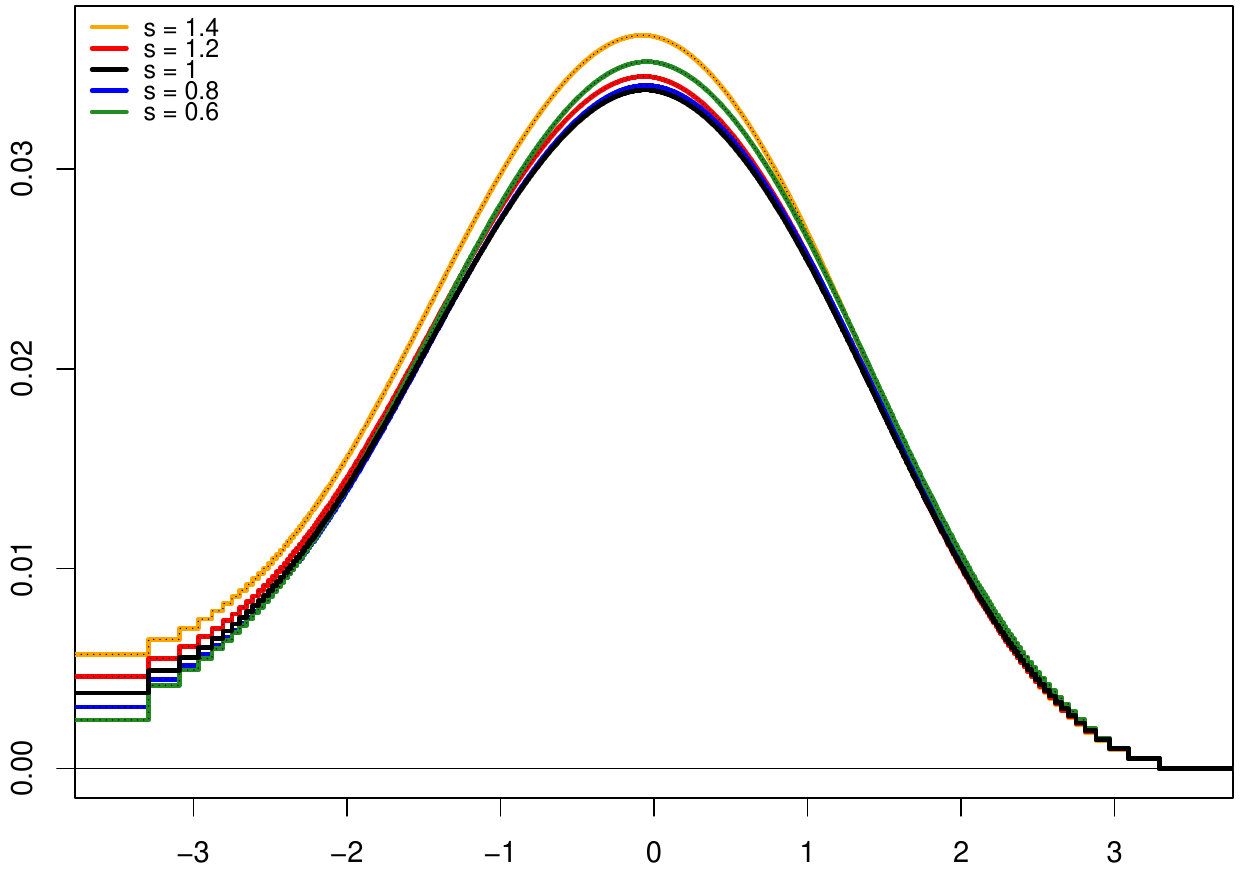}

\includegraphics[width=0.9\textwidth]{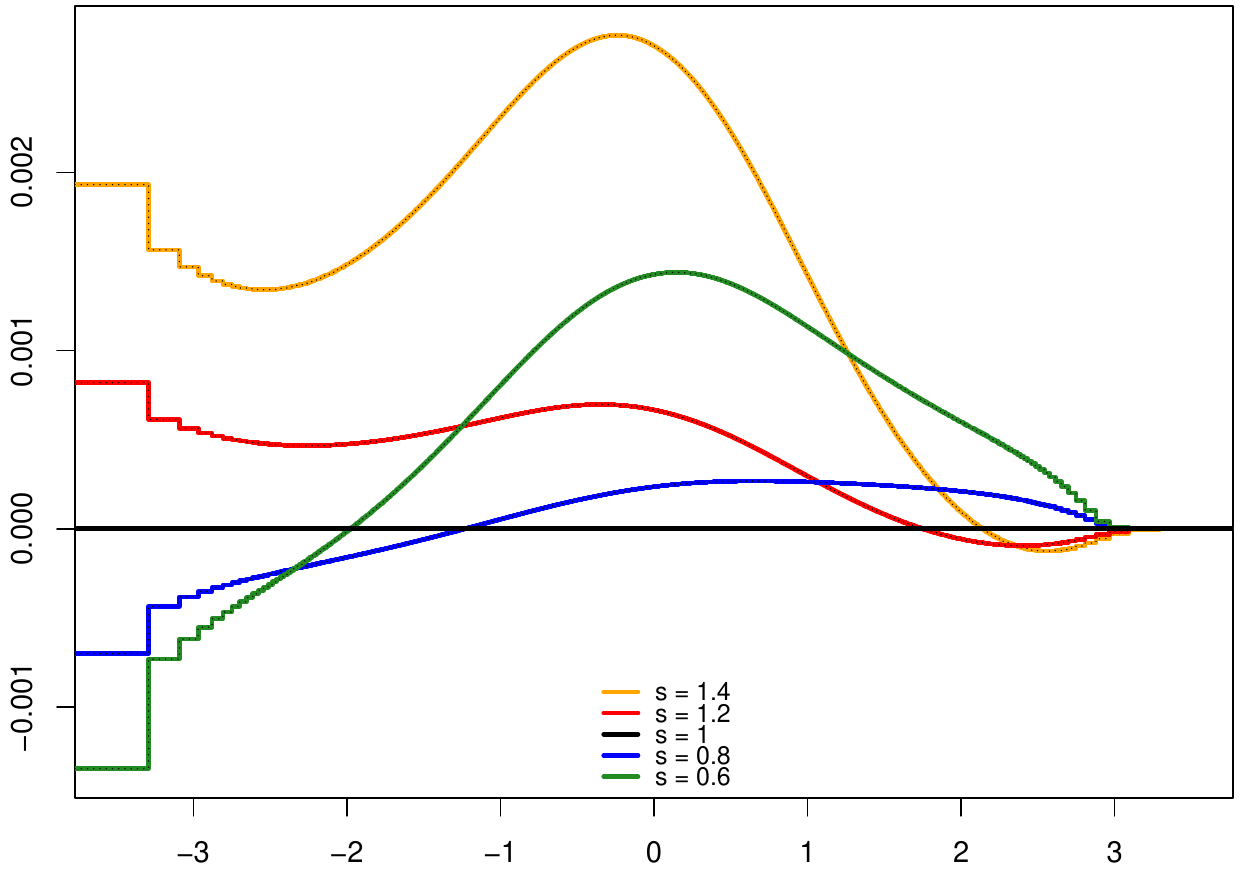}
\caption{Upper $95\%$-confidence bounds for $n = 2000$ and $s \in \{0.6,0.8,1,1.2,1.4\}$. Upper panel: centered bounds $B_{n,s,1,\alpha} - \FF_n$. Lower panel: differences $B_{n,s,1,\alpha} - B_{n,1,1,\alpha}$.}
\label{fig:UBSadd2000}
\end{figure}

In Example~3.11, we illustrated the impact of $s$ on the confidence bands $(A_{n,s,1,\alpha},B_{n,s,1,\alpha})$ by comparing these bands for $n = 500$, $\alpha = 0.05$ and $s \in \{0.6, 1, 1.4\}$. Figure~\ref{fig:UBSadd500} provides these comparisons for the same $n$ and $\alpha$ but $s \in \{0.6,0.8,1,1.2,1.4\}$. Figure~\ref{fig:UBSadd2000} shows analogous pictures for $n = 2000$.

\end{document}